# The polynomial $X^2 + Y^4$ captures its primes

By John Friedlander and Henryk Iwaniec*

*To Cherry and to Kasia*

**Table of Contents**

*JF was supported in part by NSERC grant A5123 and HI was supported in part by NSF grant DMS-9500797.

## 1. Introduction and statement of results

The prime numbers which are of the form $a^2 + b^2$ are characterized in a beautiful theorem of Fermat. It is easy to see that no prime $p = 4n-1$ can be so written and Fermat proved that all $p = 4n+1$ can be. Today we know that for a general binary quadratic form $\phi(a,b) = \alpha a^2 + \beta ab + \gamma b^2$ which is irreducible the primes represented are characterized by congruence and class group conditions. Therefore $\phi$ represents a positive density of primes provided it satisfies a few local conditions. In fact a general quadratic irreducible polynomial in two variables is known [Iw] to represent the expected order of primes (these are not characterized in any simple fashion). Polynomials in one variable are naturally more difficult and only the case of linear polynomials is settled, due to Dirichlet.

In this paper we prove that there are infinitely many primes of the form $a^2 + b^4$, in fact getting the asymptotic formula. Our main result is

THEOREM 1. *We have*

$$(1.1) \qquad \sum\sum_{a^2+b^4 \leqslant x} \Lambda(a^2 + b^4) = 4\pi^{-1}\kappa x^{\frac{3}{4}} \left\{1 + O\left(\frac{\log\log x}{\log x}\right)\right\}$$

*where $a$, $b$ run over positive integers and*

$$(1.2) \qquad \kappa = \int_0^1 (1 - t^4)^{\frac{1}{2}} \, dt = \Gamma\left(\tfrac{1}{4}\right)^2 / 6\sqrt{2\pi} \, .$$

Here of course, $\Lambda$ denotes the von Mangoldt function and $\Gamma$ the Euler gamma function. The factor $4/\pi$ is meaningful; it comes from the product (2.17) which in our case is computed in (4.8). Also the elliptic integral (1.2) arises naturally from the counting (with multiplicity included) of the integers $n \leqslant x$, $n = a^2 + b^4$ (see (3.15) and take $d = 1$). In view of these computations one can interpret $4/\pi \log x$ as the "probability" of such an integer being prime. By comparing (1.1) with the asymptotic formula in the case of $a^2 + b^2$ (change $x^{\frac{3}{4}}$ to $x$ and $t^4$ to $t^2$), we see that the probability of an integer $a^2 + b^2$ being prime is the same when we are told that $b$ is a square as it is when we are told that $b$ is not a square. In contrast to the examples given above which involved sets of primes of order $x(\log x)^{-1}$ and $x(\log x)^{-3/2}$, the one given here is much thinner.

Our work was inspired by results of Fouvry and Iwaniec [FI] wherein they proved the asymptotic formula

$$(1.3) \qquad \sum\sum_{a^2+b^2 \leqslant x} \Lambda(a^2 + b^2)\Lambda(b) = \sigma x \left\{1 + O\left((\log x)^{-A}\right)\right\}$$

with a positive constant $\sigma$ which gives the primes of the form $a^2 + b^2$ with $b$ prime.

Theorem 1 admits a number of refinements. It follows immediately from our proof that the expected asymptotic formula holds when the variables $a, b$ are restricted to any fixed arithmetic progressions, and moreover that the distribution of such points is uniform within any non-pathological planar domain. We expect, but did not check, that the methods carry over to the prime values of $\phi(a, b^2)$ for $\phi$ a quite general binary quadratic form. The method fails however to produce primes of the type $\phi(a, b^2)$ where $\phi$ is a non-homogeneous quadratic polynomial.

One may look at the equation

$$p = a^2 + b^4 \qquad (1.4)$$

in two different ways. First, starting from the sequence of Fermat primes $p = a^2 + b^2$ one may try to select those for which $b$ is square. We take the alternative approach of beginning with the integers

$$n = a^2 + b^4 \qquad (1.5)$$

and using the sieve to select primes. In the first case one would begin with a rather dense set but would then have to select a very thin subset. In our approach we begin with a very thin set but one which is sufficiently regular in behaviour for us to detect primes.

In its classical format the sieve is unable to detect primes for a very intrinsic reason, first pointed out by Selberg [Se] and known as the parity problem. The asymptotic sieve of Bombieri [Bo1], [FI1] clearly exhibits this problem. We base our proof on a new version of the sieve [FI3], which should be regarded as a development of Bombieri's sieve and was designed specifically to break this barrier and to simultaneously treat thinner sets of primes. This paper, [FI3], represents an indispensible part of the proof of Theorem 1. Originally we had intended to include it within the current paper but, expecting it to trigger other applications, we have split it off. Here, in Section 2, we briefly summarize the necessary results from that paper.

Any sieve requires good estimates for the remainder term in counting the numbers (1.5) divisible by a given integer $d$. Such an estimate is required also by our sieve and for our problem a best possible estimate of this type was provided in [FI] as a subtle deduction from the Davenport-Halberstam Theorem [DH]. It was this particular result of [FI] which most directly motivated the current work. In Section 3 we give, for completeness, that part of their work in a form immediately applicable to our problem. We also briefly describe at the end of that section how the other standard sieve assumptions listed in Section 2 follow easily for the particular sequence considered here.

In departing from the classical sieve, we introduce (see (2.11)) an additional axiom which overcomes the parity problem. As a result of this we are now required also to verify estimates for sums of the type

$$\sum_{\ell} \Big| \sum_{m} \beta(m) a_{\ell m} \Big| \tag{1.6}$$

where $\beta(m)$ is very much like the Möbius function and $a_n$ is the number of representations of (1.5) for given $n = \ell m$. The estimates for these bilinear forms constitute the major part of the paper and several of them are of interest on their own.

For example we describe an interesting by-product of one part of this work. Given a Fermat prime $p$ we define its spin $\sigma_p$ to be the Jacobi symbol $\left(\frac{s}{r}\right)$ where $p = r^2 + s^2$ is the unique representation in positive integers with $r$ odd. We show the equidistribution of the positive and negative spins $\sigma_p$. Actually we obtain this in a strong form, specifically:

THEOREM 2.    *We have*

$$\sum\sum_{r^2+s^2=p\leqslant x} \left(\frac{s}{r}\right) \ll x^{\frac{76}{77}} \tag{1.7}$$

*where $r, s$ run over positive integers with $r$ odd and $\left(\frac{s}{r}\right)$ is the Jacobi symbol.*

*Remarks.* The primes in (1.7) are not directly related to those in (1.4). As in the case of Theorem 1 the bound (1.7) holds without change when $\pi = r + is$ restricted to a fixed sector and in any fixed arithmetic progression. The exponent $\frac{76}{77}$ can be reduced by refining our estimates for the relevant bilinear forms (see Theorem $2^\psi$ in Section 26 for a more general statement and further remarks).

In studying bilinear forms of type (1.6) we are led, following some preliminary technical reductions in Sections 4 and 5, to the lattice point problem of counting points in an arithmetic progression inside the "biquadratic ellipse"

$$b_1^4 - 2\gamma b_1^2 b_2^2 + b_2^4 \leqslant x \tag{1.8}$$

for a parameter $0 < \gamma < 1$. The counting is accomplished in Sections 6–9 by a rather delicate harmonic analysis necessitated by the degree of uniformity required. The modulus $\Delta$ of the progression is very large and there are not many lattice points compared to the area of the region, at least for a given value of the parameter. It is in this counting that we exploit the great regularity in the distribution of the squares and after this step the problem of the thinness of the sifting set is gone.

There remains the task of summing the resulting main terms, that is those coming from the zero frequency in the harmonic analysis, over the relevant values of the parameter $\gamma$. The structure of these main terms is arithmetic in nature and there is some cancellation to be found in their sum, albeit requiring for its detection techniques more subtle than those needed for the nonzero frequencies. This sum is given by a bilinear form (not to be confused with (1.6)) which involves roots of quadratic congruences, again to modulus $\Delta$, which are then, as is familiar, expressed in terms of the Jacobi symbol and arithmetic progressions, this time with moduli $d$ running through the divisors of $\Delta$. Decomposing in Section 10 the relevant sum in accordance with the size of the divisors $d$ we find that we need very different techniques to deal with the divisors in different ranges.

For all but the smallest and largest ranges the relevant sum may be treated by rather general mean-value theorems of Barban-Davenport-Halberstam type. That is we need to estimate Jacobi-twisted sums on average over all residue classes and their moduli. Although, as in other theorems of this type, the results pertain for linear forms with very general coefficients, because of the rather hybrid nature of our sum (the real characters over progressions are mixed with the multiplicative inverse) new ideas are required. The goal is achieved in three steps; see Sections 11, 12, 13, their combination in Section 14 and application in Section 15.

In Section 16 we treat the smallest moduli. We require what is in essence an equidistribution result on Gaussian primes in sectors and residue classes. Now the shape of our coefficients is crucial; the cancellation will come from their resemblance to the Möbius function. The machinery for this result was developed by Hecke [He]. However, greater uniformity in the conductor is required than could have been done by him at a time prior to the famous estimate of Siegel [Si]. Siegel's work deals with $L$-functions of real Dirichlet characters rather than Grossencharacters, but today it is a routine matter to extend his argument to our case. Here we employ an elegant argument of Goldfeld [Go]. This analogue of the Siegel-Walfisz bound is applied to our problem as in the original framework and the implied constants are not computable.

There remains only the treatment of the largest moduli. We regard this as perhaps the most interesting part and hence we save it for last. In Section 17 we make some preliminary reductions and state our final goal, Proposition 17.2, for these sums. In Section 18 we show how this proposition, when combined with our earlier results, completes the proof of the main Theorem 1.

It has been familiar since the time of Dirichlet that, in dealing with various ranges in a divisor problem it is often profitable to replace large divisors by smaller ones by means of the involution $d \to |\Delta|/d$. Already this was required here in Section 12 for the final two steps in the treatment of the mid-sized moduli. An interesting feature in our case is that, due to the presence of the Jacobi symbol, the law of quadratic reciprocity plays a crucial role in this involution and an extra Jacobi symbol (of the type occurring in Theorem 2) emerges in the transformed sum (see Lemma 17.1). This extra symbol (see (20.1)) is essentially treated as a function of one complex variable and as such it is reminiscent of the Kubota symbol. This "Jacobi-Kubota symbol" later creates in Section 23, by summation over all Gaussian integers of given norm, a function on the positive integers to which we refer as a "quadratic eigenvalue".

Because the mean-value theorems of Sections 11–13 hold for such general coefficients the appearance of the Jacobi-Kubota symbol does not affect the arguments of those sections so we are able to cover completely the range of mid-sized moduli. When we again apply the Dirichlet involution, this time to transform the largest moduli, we now arrive in the same range of small moduli which have just been treated in Section 16. Now however the presence of the Jacobi-Kubota symbol destroys the previous argument, that is the theory of Hecke $L$-functions is not applicable here.

In the solution of this final part of our problem a prominent role is played by the real characters in the Gaussian domain. Dirichlet [Di] was first to treat these as an extension of the Legendre symbol. In this paper we require this Dirichlet symbol for all primary numbers, not just primes, in the same way the Jacobi symbol generalizes that of Legendre. These are introduced in Section 19. They enter our study via a kind of theta multiplier rule for the multiplication of the Jacobi-Kubota symbol, a rule we establish in Section 20.

Of particular interest are the results of Sections 21–22 concerned with general bilinear forms with the Dirichlet symbol and special linear forms with the Jacobi-Kubota symbol. This time a cancellation is received from the sign changes of these symbols rather than from those of the Möbius function which also makes an appearance arising as coefficients from our particular sieve theory. Originally, in the estimation of both of the above forms we used the Burgess bound for short character sums (thus appealing indirectly to the Riemann Hypothesis for curves, that is the Hasse-Weil Theorem). This allowed us to obtain results which in some cases are stronger than those presented here. After several attempts to simplify the original arguments we ended up with the current treatment for bilinear forms producing satisfactory results in wider ranges. Because of the wider ranges in the bilinear forms we were able to accept linear form estimates which are less uniform in the involved parameters, and consequently were able to dispense with the Burgess bound, replacing it (see Section 22) with the more elementary Polyá-Vinogradov inequality. Should we have combined the original and the present arguments then a substantial quantitative sharpening of Theorem 2 would follow.

Our estimates for the bilinear form with the Dirichlet symbol and for the special linear form with the Jacobi-Kubota symbol are then in Section 23, via the multiplier rule, transformed into corresponding results for forms in quadratic eigenvalues.

Our final job is to transform (in Sections 25 and 26) these linear and bilinear forms in the quadratic eigenvalues into sums supported on the primes (which completes Theorem 2) or sums weighted by Möbius type functions (which completes Proposition 17.2 and hence Theorem 1). There are by now a number of known combinatorial identities which can be used to achieve such a goal. The identity we introduce in Section 24 has some novel features. In particular, it enables us to reduce rather quickly from Möbius-type functions to primes and hence allows us to achieve two goals at once.

The statement of Theorem 1 may be re-interpreted in terms of the elliptic curve

$$(1.9) \qquad E: y^2 = x^3 - x \ .$$

This curve, the congruent number curve, has complex multiplication by $\mathbb{Z}[i]$ and the corresponding Hasse-Weil $L$-function

$$L_E(s) = \sum_{n=1}^{\infty} \lambda_n n^{-s}$$

is the Mellin transform of a theta series

$$f(z) = \sum_{n=1}^{\infty} \lambda_n e(nz)$$

which is a cusp form of weight two on $\Gamma_0(32)$ and is an eigenfunction of all the Hecke operators $T_p f = \lambda_p f$ . Precisely, the eigenvalues are given by

$$\lambda_n = \sum_{w\overline{w}=n}^{\wedge} w$$

where $\wedge$ restricts the summation to $w \equiv 1 (\mathrm{mod}\, 2(1+i))$, that is $w$ is primary. Hence $\lambda_p = \pi + \overline{\pi}$ if $p = \pi\overline{\pi}$ with $\pi$ primary. In particular, if $p = a^2 + b^4$, with $4 \mid a$, then $\pi = b^2 + ia$ is primary so that $\lambda_p = 2b^2$. Thus Theorem 1 gives the asymptotic formula for primes for which the Hecke eigenvalue is twice a square. Using Jacobsthal sums for these primes one expresses this property as

$$-\sum_{0<x<p/2} \left(\frac{x^3 - x}{p}\right) = \text{ square.}$$

The primes of type $p = a^2 + b^4$ give points of infinite order on the quartic twists

$$E_p : y^2 = x^3 - px \ ,$$

namely $(x, y) = (-b^2, ab)$. That this is not a torsion point follows from the Lutz-Nagell criterion. We thank Andrew Granville for pointing this out to us. The parity conjecture asserts in this case that the rank of $E_p$ is odd if $a \equiv 0 (\mathrm{mod}\, 4)$ and even if $a \equiv 2 (\mathrm{mod}\, 4)$. Recent results concerning points on quartic twists have been established by Stewart and Top [ST] improving and generalizing earlier work of Gouvea and Mazur [GM].

Further interesting connections to elliptic curves hold for primes of the form $27a^2 + 4b^6$ and there is some hope to produce such primes using our arguments in the domain $\mathbb{Z}[\zeta_3]$.

The results of this paper have been announced together with a very brief sketch of the main ideas of the proof in the paper [FI2] in the Proceedings of the National Academy of Sciences of the USA. We close here by repeating the last sentences of that paper: "Although the proofs of our results are rather lengthy and complicated we are able to avoid much of the high-powered technology frequently used in modern analytic number theory such as the bounds of Weil and Deligne. We also do not appeal to the theory of automorphic functions although experts will, in several places, detect it bubbling just beneath the surface."

*Acknowledgements.* We thank the Institute for Advanced Study for providing us with excellent conditions during the early stages of this work beginning in December 1995. HI thanks the University of Toronto for their hospitality during several short visits. We also enjoyed the hospitality of Carleton University during the CNTA conference in August 1996. We thank E. Fouvry for his encouragement. Finally we thank the referee as well as E. Fouvry, A. Granville, D. Shiu, and especially M. Watkins, for pointing out a number of minor inaccuracies.

## 2. Asymptotic sieve for primes

In this section we state a result of [FI3] in a form which is suitable for the proof of the main theorem. Let $\mathcal{A} = (a_n)$ be a sequence of real, nonnegative numbers for $n = 1, 2, 3, ...$ Our objective is an asymptotic formula for

$$S(x) = \sum_{p \leqslant x} a_p \log p$$

subject to various hypotheses familiar from sieve theory.

Let $x$ be a given number, sufficiently large in terms of $\mathcal{A}$. Put

$$A(x) = \sum_{n \leqslant x} a_n .$$

We assume the crude bounds

(2.1) $$A(x) \gg A(\sqrt{x})(\log x)^2,$$

(2.2) $$A(x) \gg x^{\frac{1}{3}} \left( \sum_{n \leqslant x} a_n^2 \right)^{\frac{1}{2}} .$$

For any $d \geqslant 1$ we write

$$(2.3) \qquad A_d(x) = \sum_{\substack{n \leqslant x \\ n \equiv 0 \pmod{d}}} a_n = g(d)A(x) + r_d(x)$$

where $g$ is a nice multiplicative function and $r_d(x)$ may be regarded as an error term which is small on average. These must of course be made more specific.

We assume that $g$ has the following properties:

$$(2.4) \qquad 0 \leqslant g(p^2) \leqslant g(p) < 1 ,$$

$$(2.5) \qquad g(p) \ll p^{-1} ,$$

and

$$(2.6) \qquad g(p^2) \ll p^{-2} .$$

Furthermore, for all $y \geqslant 2$,

$$(2.7) \qquad \sum_{p \leqslant y} g(p) = \log\log y + c + O\left((\log y)^{-10}\right) ,$$

where $c$ is a constant depending only on $g$.

We assume another crude bound

$$(2.8) \qquad A_d(x) \ll d^{-1} \tau(d)^8 A(x) \qquad \text{uniformly in } d \leqslant x^{\frac{1}{3}} .$$

We assume that the error terms satisfy

$$(2.9) \qquad \sum_{d \leqslant DL^2}{}^3 |r_d(t)| \leqslant A(x) L^{-2}$$

uniformly in $t \leqslant x$, for some $D$ in the range

$$(2.10) \qquad x^{\frac{2}{3}} < D < x .$$

Here the superscript 3 in (2.9) restricts the summation to cubefree moduli and $L = (\log x)^{2^{24}}$.

We require an estimate for bilinear forms of the type

$$(2.11) \qquad \sum_m \Big| \sum_{\substack{N < n \leqslant 2N \\ mn \leqslant x \\ (n, m\Pi)=1}} \beta(n) a_{mn} \Big| \leqslant A(x)(\log x)^{-2^{26}}$$

where the coefficients are given by

$$(2.12) \qquad \beta(n) = \beta(n, C) = \mu(n) \sum_{c | n, c \leqslant C} \mu(c) .$$

This is required for every $C$ with

$$(2.13) \qquad 1 \leqslant C \leqslant xD^{-1} ,$$

and for every $N$ with

(2.14) $$\Delta^{-1}\sqrt{D} < N < \delta^{-1}\sqrt{x} \ ,$$

for some $\Delta \geqslant \delta \geqslant 2$. Here $\Pi$ is the product of all primes $p < P$ with $P$ which can be chosen at will subject to

(2.15) $$2 \leqslant P \leqslant \Delta^{1/2^{35} \log \log x} \ .$$

PROPOSITION 2.1. *Assuming the above hypotheses, we have*

(2.16) $$S(x) = HA(x) \left\{ 1 + O\left(\frac{\log \delta}{\log \Delta}\right) \right\}$$

*where $H$ is the positive constant given by the convergent product*

(2.17) $$H = \prod_p (1 - g(p)) \left(1 - \tfrac{1}{p}\right)^{-1} \ ,$$

*and the implied constant depends only on the function $g$.*

In practice $\delta$ is a large power of $\log x$ and $\Delta$ is a small power of $x$. For most sequences all of the above hypotheses are easy to verify with the exception of (2.9) and (2.11). The hypothesis (2.9) is a traditional one while (2.11) is quite new in sieve theory.

We conclude this section by giving some technical results on the divisor function which will find repeated application in this paper.

LEMMA 2.2. *Fix $k \geqslant 2$. Any $n \geqslant 1$ has a divisor $d \leqslant n^{1/k}$ such that*

$$\tau(n) \leqslant (2\tau(d))^{\frac{k \log k}{\log 2}} \ ,$$

*and, in case $n$ is squarefree, then we may strengthen this to $\tau(n) \leqslant (2\tau(d))^k$. For any $n \geqslant 1$ we also have*

$$\tau(n) \leqslant 9 \sum_{d|n, d \leqslant n^{\frac{1}{3}}} \tau(d) \ .$$

The first two of these three statements are also from [FI3] (see Lemmata 1 and 2 there for the proofs). To prove the last of these we note that

$$\tau_3(n) \leqslant 3 \sum_{d|n, d \leqslant n^{\frac{1}{3}}} \tau(\frac{n}{d}) \ ,$$

and hence by Cauchy's inequality

$$t(n) = \tau_3(n)^2 \Big(\sum_{d|n} \tau(\frac{n}{d})^2 \tau(d)^{-1}\Big)^{-1} \leqslant 9 \sum_{d|n, d \leqslant n^{\frac{1}{3}}} \tau(d) \ .$$

On the other hand we have $t(n) \geqslant \tau(n)$ which, due to multiplicativity, can be checked by verifying on prime powers. This completes the proof of the lemma.

## 3. The sieve remainder term

In this section we verify the hypothesis (2.9) by arguments of [FI]. Given an arithmetic function $\mathfrak{z} : \mathbb{Z} \to \mathbb{C}$ we consider the sequence $\mathcal{A} = (a_n) : \mathbb{N} \to \mathbb{C}$ with

$$(3.1) \qquad a_n = \sum\sum_{a^2+b^2=n} \mathfrak{z}(b)$$

where $a$ and $b$ are integers, not necessarily positive. In our particular sequence $\mathfrak{z}$ will be supported on squares. Note that this use of $a, b$ changes from now on that in the introduction. We have

$$A_d(x) = \sum_{\substack{0<n\leqslant x \\ n\equiv 0 (\mathrm{mod}\, d)}} a_n = \sum\sum_{\substack{0<a^2+b^2\leqslant x \\ a^2+b^2\equiv 0(\mathrm{mod}\, d)}} \mathfrak{z}(b) \ .$$

We expect that $A_d(x)$ is well approximated by

$$M_d(x) = \frac{1}{d}\sum\sum_{0<a^2+b^2\leqslant x} \mathfrak{z}(b)\rho(b;d)$$

where $\rho(b;d)$ denotes the number of solutions $\alpha(\mathrm{mod}\, d)$ to the congruence

$$\alpha^2 + b^2 \equiv 0 \ (\mathrm{mod}\, d) \ .$$

For $b = 1$ we denote $\rho(1;d) = \rho(d)$; it is the multiplicative function such that

$$\rho(p^\alpha) = 1 + \chi_4(p)$$

except that $\rho(2^\alpha) = 0$ if $\alpha \geqslant 2$. Here $\chi_4$ is the character of conductor four. Thus if $4 \nmid d$

$$\rho(d) = \prod_{p|d}(1+\chi_4(p)) = \sum\nolimits^\flat_{\nu|d} \chi_4(\nu)$$

and $\rho(d) = 0$ if $4|d$. The notation $\sum^\flat$ indicates a summation over squarefree integers. For any $b$ we have

$$(3.2) \qquad \rho(b;d) = (b, d_2)\rho\left(d/(b^2,d)\right)$$

where $d_2$ denotes the largest square divisor of $d$, that is $d = d_1 d_2^2$ with $d_1$ squarefree.

LEMMA 3.1. *Suppose $\mathfrak{z}(b)$ is supported on squares and $|\mathfrak{z}(b)| \leqslant 2$. Then*

$$\sum_{d \leqslant D} |A_d(x) - M_d(x)| \ll D^{\frac{1}{4}} x^{\frac{9}{16} + \varepsilon} \tag{3.3}$$

*for any $D \geqslant 1$ and $\varepsilon > 0$, the implied constant depending only on $\varepsilon$.*

*Remarks.* This result is a modification of Lemma 4 of [FI] for our particular sequence $\mathcal{A} = (a_n)$ supported on numbers $n = a^2 + c^4$. Of course, in [FI] the authors had no reason to consider such a thin sequence so their version did not take advantage of the lacunarity of the squares.

In our case we have the individual bounds

$$\sum_{d \leqslant x} A_d(x) \ll x^{\frac{3}{4}} (\log x)^2, \tag{3.4}$$

$$\sum_{d \leqslant x} M_d(x) \ll x^{\frac{3}{4}} (\log x)^2. \tag{3.5}$$

These are derived as follows:

$$\sum_d A_d(x) \leqslant \sum\sum_{0 < a^2+b^2 \leqslant x} |\mathfrak{z}(b)| \tau(a^2 + b^2)$$

$$\leqslant 16\sqrt{x} \sum_{0 \leqslant b \leqslant \sqrt{x}} |\mathfrak{z}(b)| \sum_{d \leqslant \sqrt{x}} \rho(b;d) d^{-1}.$$

To estimate the inner sum we use the bounds $\rho(b;d) \leqslant d_2\rho(d) \leqslant \rho(d_1)\rho(d_2)d_2$, for $d$ odd, $\rho(b;d) \leqslant 4\sqrt{d}$ for $d$ a power of 2, and

$$\sum_{d \leqslant x} \rho(d) d^{-1} \ll \log x.$$

Hence we obtain (3.4) while (3.5) is derived similarly. In view of (3.4) and (3.5) our estimate (3.3) is trivial if $D > x^{3/4}$, as expected. Therefore we can assume that $D \leqslant x^{3/4}$.

The proof of Lemma 3.1 requires an application of harmonic analysis and it rests on the fact that there is an exceptional well-spacing property of the rationals $\nu/d \pmod{1}$ with $\nu$ ranging over the roots of

$$\nu^2 + 1 \equiv 0 \pmod{d}.$$

These roots correspond to the primitive representations of the modulus as the sum of two squares

$$d = r^2 + s^2 \quad \text{with} \quad (r,s) = 1.$$

By choosing $-s < r \leqslant s$ we see that each such representation gives the unique root defined by $\nu s \equiv r \pmod{d}$. Hence

$$\frac{\nu}{d} \equiv \frac{r}{sd} - \frac{\bar{r}}{s} \pmod{1}$$

where $\bar{r}$ denotes the multiplicative inverse of $r$ modulo $s$, that is $\bar{r}r \equiv 1 \pmod{s}$. Here the fraction $\bar{r}/s$ has much smaller denominator than that of $\nu/d$ whereas the other term is small, namely

$$\frac{|r|}{sd} < \frac{1}{2s^2} ,$$

except in the case $r = s = 1$ where equality holds. Therefore the points $\nu/d$ behave as if they repel each other and are distanced considerably further apart than would appear at first glance. Precisely, if $\nu_1/d_1$ and $\nu_2/d_2$ are distinct with $r_1$ and $r_2$ having the same sign and $\frac{2}{3} \leqslant \frac{s_1}{s_2} \leqslant \frac{3}{2}$ then

$$\left\| \frac{\nu_1}{d_1} - \frac{\nu_2}{d_2} \right\| > \frac{1}{s_1 s_2} - \max\left(\frac{1}{2s_1^2}, \frac{1}{2s_2^2}\right) \geqslant \frac{1}{4 s_1 s_2} > \frac{1}{4\sqrt{d_1 d_2}} .$$

Thus if the moduli are confined to an interval $\frac{8}{9}D < d \leqslant D$ then the points $\nu/d$ are spaced by $1/4D$ rather than $1/D^2$. Applying the large sieve inequality of Davenport-Halberstam [DH] for these points we derive

LEMMA 3.2.  *For any complex numbers $\alpha_n$ we have*

$$\sum_{D < d \leqslant 2D} \sum_{\nu^2 + 1 \equiv 0 (\bmod d)} \left| \sum_{n \leqslant N} \alpha_n e\left(\frac{\nu n}{d}\right) \right|^2 \ll (D + N) \|\alpha\|^2$$

*where $\|\alpha\|$ denotes the $\ell_2$-norm of $\alpha = (\alpha_n)$ and the implied constant is absolute.*

By Cauchy's inequality Lemma 3.2 yields

(3.6) $$\sum_{d \leqslant D} \sum_{\nu^2 + 1 \equiv 0 (d)} \left| \sum_{n \leqslant N} \alpha_n e\left(\frac{\nu n}{d}\right) \right| \ll D^{\frac{1}{2}} (D + N)^{\frac{1}{2}} \|\alpha\| .$$

From this we shall derive a bound for general linear forms in the arithmetic functions

(3.7) $$\rho(k, \ell; d) = \sum_{\nu^2 + \ell^2 \equiv 0 (\bmod d)} e(\nu k/d) .$$

LEMMA 3.3.  *For any complex numbers $\xi(k, \ell)$ we have*

$$\sum_{d \leqslant D} \left| \sum_{\substack{0 < k \leqslant K \\ 0 < \ell \leqslant L}} \xi(k, \ell) \rho(k, \ell; d) \right| \ll \left(D + \sqrt{DKL}\right)(DKL)^\varepsilon \|\xi\|$$

*where $\|\xi\|$ denotes the $\ell_2$-norm of $\xi = (\xi(k, \ell))$; that is*

$$\|\xi\|^2 = \sum_{\substack{0 < k \leqslant K \\ 0 < \ell \leqslant L}} |\xi(k, \ell)|^2 ,$$

*and the implied constant depends only on $\varepsilon$.*

The functions $\rho(k,\ell;d)$ serve as "Weyl harmonics" for the equidistribution of roots of the congruence

$$(3.8) \qquad \nu^2 + \ell^2 \equiv 0 \pmod{d} .$$

Note that $\rho(0,\ell;d) = \rho(\ell;d)$ is the multiplicative function which appears in the expected main term $M_d(x)$ and this is expressed simply in terms of $\rho(d)$ by (3.2). If $k \neq 0$ then $\rho(k,\ell;d)$ is more involved but one can at least reduce this to the case $\ell = 1$. Specifically, letting $(d,\ell^2) = \gamma\delta^2$ with $\gamma$ squarefree so $d = \gamma\delta^2 d'$, $\ell = \gamma\delta\ell'$, one shows that

$$(3.9) \qquad \rho(k,\ell;d) = \delta\rho(k'\ell',1;d')$$

provided that $k = \delta k'$ is a multiple of $\delta$, while $\rho(k,\ell;d)$ vanishes if $k$ is not divisible by $\delta$. By this we obtain

$$\sum_{d \leqslant D} \Big| \sum_{\substack{0 < k \leqslant K \\ 0 < \ell \leqslant L}} \xi(k,\ell)\rho(k,\ell;d) \Big|$$
$$\leqslant \sum_{\gamma\delta^2 d \leqslant D} \sum \sum \delta \Big| \sum_{\substack{0 < k \leqslant K/\delta \\ 0 < \ell \leqslant L/\gamma\delta \\ (\ell,d)=1}} \xi(\delta k, \gamma\delta\ell)\rho(k\ell,1;d) \Big| .$$

Ignoring the condition $(\ell,d) = 1$ we would get the bound of Lemma 3.3 by applying (3.6) directly. However this co-primality condition can be inserted at no extra cost by Möbius inversion and this completes the proof of Lemma 3.3.

Now we are ready to prove Lemma 3.1. We begin by smoothing the sum $A_d(x)$ with a function $f(u)$ supported on $[0,x]$ such that

$$\begin{aligned} f(u) &= 1 & &\text{if } 0 < u \leqslant x - y , \\ f^{(j)}(u) &\ll y^{-j} & &\text{if } x - y < u < x , \end{aligned}$$

where $y$ will be chosen later subject to $x^{\frac{1}{2}} < y < x$ and the implied constant depends only on $j$. Our intention is to apply Fourier analysis to the sum

$$A_d(f) = \sum_{n \equiv 0 \pmod{d}} a_n f(n)$$

rather than directly to $A_d(x)$. By a trivial estimation the difference is

$$(3.10) \qquad \sum_{d \leqslant D} \big| A_d(x) - A_d(f) \big| \ll y x^{-\frac{1}{4}+\varepsilon} .$$

In $A_d(f)$ we split the summation over $a$ into classes modulo $d$ getting

$$A_d(f) = \sum_b \mathfrak{z}(b) \sum_{\alpha^2 + b^2 \equiv 0(d)} \sum_{a \equiv \alpha(d)} f(a^2 + b^2) .$$

It is convenient to first remove the contribution coming from terms with $b = 0$, since these are not covered by Lemma 3.3. This contribution is

$$\mathfrak{z}(0) \sum_{a^2 \equiv 0(d)} f(a^2) = \mathfrak{z}(0) \sum_a f((ad_1 d_2)^2) \ll \frac{\sqrt{x}}{d_1 d_2}.$$

For the nonzero values of $b$ we expand the above inner sum into Fourier series by Poisson's formula getting

$$\sum_{a \equiv \alpha(d)} f(a^2 + b^2) = \frac{1}{d} \sum_k e\left(\frac{\alpha k}{d}\right) \int_{-\infty}^{\infty} f(t^2 + b^2) e\left(\frac{tk}{d}\right) dt.$$

Hence the smooth sum $A_d(f)$ has the expansion

(3.11) $$A_d(f) = \frac{2}{d} \sum_{b \neq 0} \mathfrak{z}(b) \sum_k \rho(k, b; d) I(k, b; d) + O\left(\frac{\sqrt{x}}{d_1 d_2}\right)$$

where $I(k, b; d)$ is the Fourier integral

$$I(k, b; d) = \int_0^{\infty} f(t^2 + b^2) \cos(2\pi tk/d) dt.$$

The main term comes from $k = 0$ which gives

$$M_d(f) = \frac{2}{d} \sum_b \mathfrak{z}(b) \rho(b; d) I(0, b; d).$$

Since in this case the integral approximates to the sum, precisely

$$2I(0, b; d) = \sum_{a^2 + b^2 \leqslant x} 1 + O\left(y(x + y - b^2)^{-\frac{1}{2}}\right),$$

the difference between the expected main terms satisfies

$$M_d(f) - M_d(x) \ll \frac{y}{d} \sum_{c^4 \leqslant x} \rho(c^2; d)(x + y - c^4)^{-\frac{1}{2}}.$$

Summing over moduli we first derive by the same arguments which led us to (3.5) that

$$\sum_{d \leqslant D} d^{-1} \rho(c^2; d) \ll (\log 2D)^2,$$

and then summing over $c$ we arrive at

(3.12) $$\sum_{d \leqslant D} | M_d(f) - M_d(x) | \ll yx^{-\frac{1}{4}} (\log x)^2.$$

For positive frequencies $k$ we shall estimate $I(k, b; d) = I(-k, b; d)$ by repeated partial integration. We have

$$\frac{\partial^j}{\partial t^j} f(t^2 + b^2) = \sum_{0 \leqslant 2i \leqslant j} c_{ij} t^{j-2i} f^{(j-i)}(t^2 + b^2) \ll \left(\frac{\sqrt{x}}{y}\right)^j,$$

with some positive constants $c_{ij}$, whence

$$I(k, b; d) \ll \sqrt{x} \left(\frac{d\sqrt{x}}{ky}\right)^j$$

for any $j \geqslant 0$. This shows that $I(k, b; d)$ is very small if $k \geqslant K = Dy^{-1}x^{\frac{1}{2}+\varepsilon}$ by choosing $j = j(\varepsilon)$ sufficiently large. Estimating the tail of the Fourier series (3.11) trivially we are left with

$$A_d(f) = M_d(f) + \frac{4}{d} \sum_{b \neq 0} \mathfrak{z}(b) \sum_{0 < k \leqslant K} \rho(k, b; d) I(k, b; d) + O\left(\frac{\sqrt{x}}{d_1 d_2}\right).$$

To separate the modulus $d$ from $k, b$ in the Fourier integral we write

$$I(k, b; d) = \sqrt{x} k^{-1} \int_0^\infty f(xt^2 k^{-2} + b^2) \cos(2\pi t \sqrt{x}/d) dt$$

by changing the variable $t$ into $t\sqrt{x}/k$. Note that the new variable lies in the range $0 < t < k$. Hence $\big| A_d(f) - M_d(f) \big|$ is bounded by

$$\frac{4\sqrt{x}}{d} \int_0^K \bigg| \sum\sum_{\substack{0 < b \leqslant \sqrt{x} \\ t < k \leqslant K}} \frac{\mathfrak{z}(b)}{k} f\left(\frac{xt^2}{k^2} + b^2\right) \rho(k, b; d) \bigg| dt + O\left(\frac{\sqrt{x}}{d_1 d_2}\right).$$

Recall that $\mathfrak{z}(b)$ is supported on squares; $b = c^2$ with $|c| \leqslant C = x^{\frac{1}{4}}$. Applying Lemma 3.3 to the relevant triple sum and then integrating over $0 < t < K$ we obtain

$$\sum_{d \leqslant D} d \big| A_d(f) - M_d(f) \big| \ll \sqrt{x} \left(D + C\sqrt{DK}\right)(CK)^{\frac{1}{2}+\varepsilon}$$

$$\ll D^{\frac{3}{2}} y^{-1} x^{\frac{11}{8}+\varepsilon}.$$

Hence the smooth remainder satisfies

(3.13) $$\sum_{d \leqslant D} \big| A_d(f) - M_d(f) \big| \ll D^{\frac{1}{2}} y^{-1} x^{\frac{11}{8}+\varepsilon}.$$

Finally, on combining (3.10), (3.12) and (3.13) we obtain (3.3) by the choice $y = D^{\frac{1}{4}} x^{\frac{13}{16}}$.

From now on $\mathfrak{z}(b)$ is equal to 2 if $b = c^2 > 0$, $\mathfrak{z}(0) = 1$, and $\mathfrak{z}(b) = 0$ otherwise. In other words

(3.14) $$\mathfrak{z}(b) = \sum_{c^2 = b} 1$$

where $c$ is any integer. Note that $\mathfrak{z}(b)$ is the Fourier coefficient of the classical theta function. For this choice of $\mathfrak{z}$ we shall evaluate the main term $M_d(x)$ more precisely.

LEMMA 3.4.  *For $d$ cubefree we have*

(3.15) $$M_d(x) = g(d)\kappa x^{\frac{3}{4}} + O\left(h(d)x^{\frac{1}{2}}\right)$$

*where $\kappa$ is the constant given by the elliptic integral* (1.2) *and $g(d)$, $h(d)$ are the multiplicative functions given by*

(3.16) $\quad g(p)p = 1 + \chi_4(p)\left(1 - \frac{1}{p}\right), \quad g(p^2)p^2 = 1 + \rho(p)\left(1 - \frac{1}{p}\right),$

$\quad h(p)p = 1 + 2\rho(p), \quad\quad\quad h(p^2)p^2 = p + 2\rho(p),$

*except that $g(4) = \frac{1}{4}$.*

*Proof.* We have

$$M_d(x) = \frac{2}{d} \sum_{|c| \leqslant x^{\frac{1}{4}}} \rho(c^2; d) \left\{(x - c^4)^{\frac{1}{2}} + O(1)\right\}.$$

Since $d$ is cubefree we can write $d = d_1 d_2^2$ with $d_1 d_2$ squarefree, so that we have $\rho(\ell^2; d) = (\ell, d_2)\rho(d_1 d_2/(\ell, d_1 d_2))$ except for $d_2$ even and $\ell$ odd, in which case $\rho(\ell^2; d) = 0$. Hence, for $d$ not divisible by 4 we have

$$M_d(x) = \frac{2}{d} \sum_{\substack{\nu_1 | d_1 \\ \nu_2 | d_2}} \nu_2 \rho\left(\frac{d_1 d_2}{\nu_1 \nu_2}\right) \sum_{\substack{|c| \leqslant x^{\frac{1}{4}} \\ (c, d_1 d_2) = \nu_1 \nu_2}} \left\{(x - c^4)^{\frac{1}{2}} + O(1)\right\}$$

$$= \frac{2}{d} \sum_{\substack{\nu_1 | d_1 \\ \nu_2 | d_2}} \nu_2 \rho\left(\frac{d_1 d_2}{\nu_1 \nu_2}\right) \left\{\varphi\left(\frac{d_1 d_2}{\nu_1 \nu_2}\right)\frac{2\kappa x^{\frac{3}{4}}}{d_1 d_2} + O\left(\tau\left(\frac{d_1 d_2}{\nu_1 \nu_2}\right)x^{\frac{1}{2}}\right)\right\}.$$

This formula gives (3.15) with

$$g(d)d = \left(\sum_{\nu_1 | d_1} \rho\left(\frac{d_1}{\nu_1}\right)\varphi\left(\frac{d_1}{\nu_1}\right)d_1^{-1}\right)\left(\sum_{\nu_2 | d_2} \rho\left(\frac{d_2}{\nu_2}\right)\varphi\left(\frac{d_2}{\nu_2}\right)\frac{\nu_2}{d_2}\right),$$

$$h(d)d = \left(\sum_{\nu_1 | d_1} \rho\left(\frac{d_1}{\nu_1}\right)\tau\left(\frac{d_1}{\nu_1}\right)\right)\left(\sum_{\nu_2 | d_2} \rho\left(\frac{d_2}{\nu_2}\right)\tau\left(\frac{d_2}{\nu_2}\right)\nu_2\right),$$

which completes the proof of Lemma 3.4 in this case. For $d$ cubefree and divisible by 4 the above argument goes through except that, as noted, $\rho(\ell^2, d) = 0$ for $\ell$ odd. This implies that, in the summation, $c$ and hence $\nu_2$ must be restricted to even numbers. This makes the value of $g(4)$ exceptional. □

We define the error term

(3.17) $$r_d(x) = A_d(x) - g(d)A(x).$$

By Lemma 3.4 for $d = 1$ we get

(3.18)   $$A(x) = 4\kappa x^{\frac{3}{4}} + O\left(x^{\frac{1}{2}}\right) ;$$

thus for $d$ cubefree the error term satisfies

$$r_d(x) = A_d(x) - M_d(x) + O\left(h(d)x^{\frac{1}{2}}\right) .$$

Note that

$$\sum_{d \leqslant x}{}^3 h(d) \leqslant \prod_{p \leqslant x} (1 + h(p))\left(1 + h(p^2)\right) \ll (\log x)^4,$$

where the superscript 3 restricts to cubefree numbers. This together with Lemma 3.1 implies

PROPOSITION 3.5.  *We have for all $t \leqslant x$,*

(3.19)   $$\sum_{d \leqslant D}{}^3 |r_d(t)| \ll D^{\frac{1}{4}} x^{\frac{9}{16} + \varepsilon} .$$

The restriction to cubefree moduli in (3.19) is not necessary but it is sufficient for our needs. The fact that we are able to make this restriction will be technically convenient in a number of places specifically because cubefree numbers $d$ possess the property that they can be decomposed as $d = d_1 d_2^2$ with $d_1$, $d_2$ squarefree and $(d_1, d_2) = 1$.

Proposition 3.5 verifies one of the two major hypotheses of the ASP (Asymptotic Sieve for Primes), namely (2.9) with $D = x^{\frac{3}{4} - 5\varepsilon}$ by a comfortable margin and indeed is, apart from the $\varepsilon$, the best that one can hope for. The hypotheses (2.4), (2.5), and (2.6) are easily verified by an examination of (3.16). The asymptotic formula (2.7) is derived from the Prime Number Theorem for the primes in residue classes modulo four. Next, the crude bounds (2.1), (2.2) and (2.8) are obvious in our case. More precisely, one can derive by elementary arguments that $A_d(x) \ll d^{-1}\tau(d)A(x)$ uniformly for $d \leqslant x^{\frac{1}{2} - \varepsilon}$ in place of (2.8). Therefore we are left with the problem of establishing the second major hypothesis of the ASP, namely the bilinear form bound (2.11).

## 4. The bilinear form in the sieve: Refinements

Throughout $a_n$ denotes the number of integral solutions $a$, $c$ to

(4.1)   $$a^2 + c^4 = n .$$

Recall from the previous section that (see (3.18))

(4.2)   $$A(x) = \sum_{n \leqslant x} a_n = 4\kappa x^{\frac{3}{4}} + O\left(x^{\frac{1}{2}}\right) .$$

In this section we give a preliminary analysis of the bilinear forms

$$(4.3) \qquad \mathcal{B}(x; N) = \sum_m \Big| \sum_{\substack{N < n \leqslant 2N \\ mn \leqslant x \\ (n, m\Pi) = 1}} \beta(n) a_{mn} \Big|$$

with coefficients $\beta(n)$ given by (2.12) and $\Pi$ the product of primes $p \leqslant P$ with $P$ in the range

$$(4.4) \qquad (\log \log x)^2 \leqslant \log P \leqslant (\log x)(\log \log x)^{-2}.$$

Although the sieve does not require any lower bound for $P$, that is $\Pi = 1$ is permissible, we introduce this as a technical device which greatly simplifies a large number of computations. With slightly more work we could relax the lower bound for $P$ to a suitably large power of $\log x$ and still obtain the same results.

Note the bound

$$\mathcal{B}(x; N) \ll A(x)(\log x)^4$$

uniformly in $N \leqslant x^{\frac{1}{2}}$. This follows from (3.1) by a trivial estimation, but we need the stronger bound (2.11). We shall establish the following improvement:

PROPOSITION 4.1.    *Let $\eta > 0$ and $A > 0$. Then we have*

$$(4.5) \qquad \mathcal{B}(x; N) \ll A(x)(\log x)^{4-A}$$

*for every $N$ with*

$$(4.6) \qquad x^{\frac{1}{4} + \eta} < N < x^{\frac{1}{2}} (\log x)^{-B}$$

*and the coefficients $\beta(n, C)$ given by (2.12) with $1 \leqslant C \leqslant N^{1-\eta}$. Here $B$ and the implied constant in (4.5) need to be taken sufficiently large in terms of $\eta$ and $A$.*

By virtue of the results presented in the previous sections Proposition 4.1 is more than sufficient to infer the formula

$$(4.7) \qquad \sum_{p \leqslant x} a_p \log p = HA(x) \left\{ 1 + O\left( \frac{\log \log x}{\log x} \right) \right\}$$

(it suffices to have (4.5) with $A = 2^{26} + 4$ and $x^{3/8 - \eta} < N < x^{1/2}(\log x)^{-B}$ for some $\eta > 0$ and $B > 0$). In this formula $H$ is given by (2.17) with $g(p)p = 1 + \chi_4(p)(1 - \frac{1}{p})$, whence

$$(4.8) \qquad H = \prod_p \left( 1 - \chi_4(p) p^{-1} \right) = L(1, \chi_4)^{-1} = \frac{4}{\pi}.$$

Therefore (4.7), (4.8) and (4.2) yield the asymptotic formula (1.1) of our main theorem. Note that in the formulation of Theorem 1 we restricted to representations by positive integers thus obtaining a constant equal to one fourth of that in (4.7).

It remains to prove Proposition 4.1, and this is the heart of the problem. In this section we make a few technical refinements of the bilinear form $\mathcal{B}(x; N)$ which will be useful in the sequel.

First of all the coefficients $\beta(n)$ can be quite large which causes a problem in Section 9. More precisely we have $|\beta(n)| \leqslant \tau(n)$ so the problem occurs for a few $n$ for which $\tau(n)$ is exceptionally large. We remove these terms now because it will be more difficult to control them later. Let $\mathcal{B}'(x; N)$ denote the partial sum of $\mathcal{B}(x; N)$ restricted by

$$\tau(n) \leqslant \tau \tag{4.9}$$

where $\tau$ will be chosen as a large power of $\log x$. The complementary sum is estimated trivially by

$$\sum\sum_{\substack{mn \leqslant x \\ \tau(n) > \tau}} \mu^2(mn)\tau(n)a_{mn} \leqslant \tau^{-1} \sum\sum_{mn \leqslant x} \mu^2(mn)\tau(n)^2 a_{mn} = \tau^{-1} \sum_{n \leqslant x}^{\flat} \tau_5(n) a_n \ .$$

By Lemma 2.2 we have $\tau_5(n) \leqslant \tau(n)^{\log 5/\log 2} \leqslant (2\tau(d))^7$ for some $d \mid n$ with $d \leqslant n^{1/3}$. Hence the above sum is bounded by

$$\sum_{d \leqslant x^{\frac{1}{3}}} (2\tau(d))^7 A_d(x) \ll A(x) \sum_{d \leqslant x^{\frac{1}{3}}} \tau(d)^7 g(d) \ll A(x)(\log x)^{2^7} \ ,$$

which gives

$$\mathcal{B}(x; N) = \mathcal{B}'(x; N) + O\left(\tau^{-1} A(x)(\log x)^{128}\right) \ . \tag{4.10}$$

To make this bound admissible for (4.5) we assume that

$$\tau \geqslant (\log x)^{A+124} \ . \tag{4.11}$$

While the restriction (4.9) will help us to estimate the error term in the lattice point problem it is not desired for the main term because the property $\tau(n) \leqslant \tau$ is not multiplicatively stable (to the contrary of $(n, \Pi) = 1$). In the resulting main term in Section 10 we shall remove the restriction (4.9) by the same method which allowed us to install it here.

In numerous transformations of $\mathcal{B}(x; N)$ we shall be faced with technical problems such as separation of variables or handling abnormal structures. When resolving these problems we wish to preserve the nature of the coefficients $\beta(n)$ (think of $\beta(n)$ as being the Möbius function). Thus we should avoid any technique which uses long integration because it corrupts $\beta(n)$.

To get hold of the forthcoming problems we reduce the range of the inner sum of $\mathcal{B}'(x; N)$ to short segments of the type

$$(4.12) \qquad N' < n \leqslant (1+\theta)N'$$

where $\theta^{-1}$ will be a large power of $\log N$, and we replace the restriction $mn \leqslant x$ by $mN \leqslant x$. This reduction can be accomplished by splitting into at most $\theta^{-1}$ such sums and estimating the residual contribution trivially. In fact we get a better splitting by means of a smooth partition of unity. This amounts to changing $\beta(n)$ into

$$(4.13) \qquad \beta(n) = p(n)\mu(n) \sum_{c|n, c \leqslant C} \mu(c)$$

where $p$ is a smooth function supported on the segment (4.12) for some $N'$ which satisfies $N < N' < 2N$. It will be sufficient that $p$ be twice differentiable with

$$(4.14) \qquad p^{(j)} \ll (\theta N)^{-j}, \quad j = 0, 1, 2 .$$

One needs at most $2\theta^{-1}$ such partition functions to cover the whole interval $N < n \leqslant 2N$ with multiplicity one except for the points $n$ with $|mn - x| < \theta x$, $|n - N| < \theta N$ or $|n - 2N| < \theta N$. However, these boundary points contribute at most $O\left(\theta A(x)(\log x)^4\right)$ by a straightforward estimation so we have

$$\mathcal{B}'(x; N) = \sum_p \mathcal{B}'_p(x; N) + O\left(\theta A(x)(\log x)^4\right)$$

where $p$ ranges over the relevant partition functions and $\mathcal{B}'_p(x; N)$ is the corresponding smoothed form of $\mathcal{B}'(x; N)$. To make the above bound for the residual contribution admissible for (4.5) we assume that

$$(4.15) \qquad \theta = (\log x)^{-A'}$$

with $A' \geqslant A$. We do not specialize $A'$ for the time being, in fact not until Section 18, but it will be much larger than $A$. In other words $\theta$ is quite a bit smaller than the factor

$$(4.16) \qquad \vartheta = (\log x)^{-A} ,$$

which we aim to save in the bound (4.5). Since the number of smoothed forms does not exceed $2\theta^{-1}$ it suffices to show that each of these satisfies

$$(4.17) \qquad \mathcal{B}'_p(x; N) \ll \vartheta\theta A(x)(\log x)^4 .$$

Next we split the outer summation into dyadic segments

$$(4.18) \qquad M < m \leqslant 2M ;$$

to save space we shall sometimes write this as $m \sim M$. Since the contribution of terms with $m \leqslant \vartheta x N^{-1}$ is absorbed by the bound (4.17) we are left with

$$(4.19) \qquad \vartheta x < MN < x \, .$$

Note that (4.6) and (4.19) imply $M \geqslant N$ since we may require $B \geqslant 2A$.

After the splitting we obtain bilinear forms of type

$$(4.20) \qquad \mathcal{B}^*(M, N) = \sum\sum_{(m,n)=1} \alpha(m)\beta(n) a_{mn} \, ,$$

where we allow $\alpha(m)$ to be any complex numbers supported on (4.18) with $|\alpha(m)| \leqslant 1$, while $\beta(n)$ are given by (4.13). Here and hereafter, in order to save frequent writing of the summation conditions

$$(4.21) \qquad (n, \Pi) = 1$$

$$(4.22) \qquad \tau(n) \leqslant \tau$$

we rather regard these as restrictions on the support of $\beta(n)$. Occasionally it will be appropriate to remind the reader of this convention. It now suffices to show that for every $N$ satisfying (4.6) and $MN$ satisfying (4.19) we have

$$(4.23) \qquad \mathcal{B}^*(M, N) \ll \vartheta \theta (MN)^{\frac{3}{4}} (\log MN)^4.$$

## 5. The bilinear form in the sieve: Transformations

Typically for general bilinear forms one applies Cauchy's inequality in order to smooth and then to execute the outer summation. However in the case of our special form $\mathcal{B}^*(M, N)$ the application of Cauchy's inequality at this stage would be premature. This is due to the multiplicity of representations

$$(5.1) \qquad a_{mn} = \sum_{a^2+c^4=mn} 1$$

where $a$, $c$ run over all integers. This multiplicity is locked into the inner sum and we do not wish to amplify it by squaring since that would have fatal effects on the harmonic analysis when the time comes to count the lattice points. Therefore, we release that part of the multiplicity which is accommodated in the outer variable and we smooth out this part rather than amplify by applying Cauchy's inequality. In order to be able to extract this hidden part of the multiplicity we shall write the solutions $a^2 + b^2 = mn$ in $a, b \in \mathbb{Z}$ in terms of Gaussian integers $w, z \in \mathbb{Z}[i]$. Not only is our entry to the Gaussian domain necessary but it will also clarify the arguments. The arithmetic we are going to apply lies truly in the field $\mathbb{Q}(i)$. On the other hand, some of our arguments such as the quadratic reciprocity law, seem more familiar when performed in $\mathbb{Q}$. Thus we make no appeal to properties of the biquadratic residue symbol, instead we work with the Dirichlet symbol (see (19.11)) which is an extension of that of Jacobi to $\mathbb{Z}[i]$.

Since $(m,n) = 1$ we have by the unique factorization in $\mathbb{Z}[i]$

$$(5.2) \qquad a_{mn} = \frac{1}{4} \sum_{|w|^2 = m} \sum_{|z|^2 = n} \mathfrak{z}(\operatorname{Re} \overline{w} z)$$

where $\frac{1}{4}$ accounts for the four units $1, i, i^2, i^3$ in $\mathbb{Z}[i]$. Since $n$ is odd so is $z$. Multiplying $w$ and $z$ by a unit one can rotate $z$ to a number satisfying

$$(5.3) \qquad z \equiv 1 \pmod{2(1+i)}.$$

Such a number is called primary; it is determined uniquely by its ideal. In terms of coordinates $z = r + is$ the congruence (5.3) means

$$(5.4) \qquad r \equiv 1 \pmod{2}$$

and

$$(5.5) \qquad s \equiv r - 1 \pmod{4}$$

so that $r$ is odd and $s$ is even.

By (5.2) we can express the bilinear form (4.20) as

$$(5.6) \qquad \mathcal{B}^*(M, N) = \sum\sum_{(w\overline{w}, z\overline{z}) = 1} \alpha_w \beta_z \mathfrak{z}(\operatorname{Re} \overline{w} z)$$

where $\alpha_w = \alpha(|w|^2)$ and $\beta_z = \beta(|z|^2)$. Here we assume that $z$ runs over primary numbers so the multiplicity four does not occur in (5.6). In the sequel we regard $\beta_z$ as a function supported on numbers having a fixed residue class modulo eight, say

$$(5.7) \qquad z \equiv z_0 \pmod{8}$$

where $z_0$ is primary. This can be accomplished by splitting $\mathcal{B}^*(M, N)$ into eight such classes. Recall that we also have the restrictions for the support of $\beta_z$ coming from (4.21), (4.22). These read as

$$(5.8) \qquad (z, \Pi) = 1$$

$$(5.9) \qquad \tau(|z|^2) \leqslant \tau.$$

To obtain the factorization (5.2) it was very convenient to have the condition $(m, n) = 1$. This condition, which meanwhile has become $(w\overline{w}, z\overline{z}) = 1$ in (5.6), is now a hindrance since we want $w$ to run freely. We shall remove the condition $(w\overline{w}, z\overline{z}) = 1$ by estimating trivially the complementary form

$$\widetilde{\mathcal{B}}(M, N) = \sum\sum_{(w\overline{w}, z\overline{z}) \neq 1} \alpha_w \beta_z \mathfrak{z}(\operatorname{Re} \overline{w} z).$$

To this end we take advantage of (5.8) getting

$$\widetilde{\mathcal{B}}(M,N) \ll \sum_{p>P} \sum\sum_{\substack{M<u^2+v^2\leqslant 2M \\ u^2+v^2\equiv 0(p)}} \sum\sum_{\substack{N<r^2+s^2\leqslant 2N \\ r^2+s^2\equiv 0(p) \\ ur+vs=\square}} |\beta(r^2+s^2)|.$$

Since $\beta(n)$ is supported on odd squarefree numbers we have $(r,s) = 1$. Note that $p \nmid rs$. Put $ur+vs = c^2$ with $|c| \leqslant 2(MN)^{\frac{1}{4}}$. Given $c,r,s$ the residue class of $u \pmod{ps/(c,p)}$ is fixed and then $v$ is determined. Therefore the number of points $w = u+iv$ is bounded by $O(1+\sqrt{M}(c,p)/ps)$. By symmetry we can replace this by $O(1+\sqrt{M}(c,p)/pr)$ and hence by $O(1+\sqrt{M}(c,p)/p\sqrt{N})$. Summing over $c$ and noting that $P \ll (MN)^{\frac{1}{4}}$, we find that our complementary form satisfies

$$\begin{aligned}\widetilde{\mathcal{B}}(M,N) &\ll (MN)^{\frac{1}{4}}\left(1+\sqrt{M}/P\sqrt{N}\right) \sum_p \sum_{\substack{N<n\leqslant 2N \\ n\equiv 0(\bmod p)}} \tau(n)|\beta(n)| \\ &\ll (MN)^{\frac{1}{4}}\left(1+\sqrt{M}/P\sqrt{N}\right) N(\log N)^2 .\end{aligned}$$

Hence, adding $\widetilde{\mathcal{B}}(M,N)$ to $\mathcal{B}^*(M,N)$ we conclude that

(5.10) $$\mathcal{B}^*(M,N) = \mathcal{B}(M,N) + O\left((M^{\frac{1}{4}}N^{\frac{5}{4}} + P^{-1}M^{\frac{3}{4}}N^{\frac{3}{4}})(\log N)^3\right)$$

where $\mathcal{B}(M,N)$ is the free bilinear form

$$\mathcal{B}(M,N) = \sum_w \sum_z \alpha_w \beta_z \Im(\operatorname{Re}\overline{w}z) .$$

Note that the error term in (5.10) is admissible for (4.23) if $N \leqslant \vartheta\theta\sqrt{MN}$ and $\vartheta\theta P \geqslant 1$. The first condition is satisfied if

(5.11) $$B \geqslant \tfrac{3}{2}A + A',$$

by virtue of (4.19) and (4.6). The second condition requires $P \geqslant (\log x)^{A+A'}$. Actually in (5.24) below we shall impose a stronger condition for $P$.

Next we are going to assume that $\beta_z$ is supported in a narrow sector

(5.12) $$\varphi < \arg z \leqslant \varphi + 2\pi\theta$$

for some $-\pi < \varphi < \pi$ where $\theta$ is the same as in (4.12). This can be accomplished by splitting according to a smooth partition of unity (without any residual contribution because there is no boundary). We need only a $\mathcal{C}^2$-class partition. In other words we attach to $\beta_z$ a periodic function $q(\alpha)$ of period $2\pi$ supported on $\varphi < \alpha \leqslant \varphi + 2\pi\theta$ such that

$$q^{(j)} \ll \theta^{-j} , \quad j = 0,1,2 .$$

Thus from now on

(5.13) $$\beta_z = q(\alpha)p(n)\mu(n) \sum_{c|n, c \leqslant C} \mu(c)$$

where $\alpha = \arg z$ and $n = |z|^2$. Note that $\beta_z$ is supported on numbers $z = r + is$ with $|z|^2 = r^2 + s^2$ squarefree, which implies $(r,s) = 1$ so $z$ is primitive. Since $z$ is also odd we have $(z, \bar{z}) = 1$ and this property will prove to be convenient in several places.

The intersection of the annulus (4.12) with the sector (5.12) is a polar box

(5.14) $$\mathfrak{B} = \{z : N' < |z|^2 \leqslant (1+\theta)N', \; \varphi < \arg z \leqslant \varphi + 2\pi\theta\}$$

of volume
$$\operatorname{vol} \mathfrak{B} = \pi\theta^2 N' \sim \pi\theta^2 N \; .$$

Since the number of polar boxes is $O(\theta^{-2})$ we now need to prove that for the bilinear form $\mathcal{B}(M, N)$ restricted smoothly to a box we have

(5.15) $$\mathcal{B}(M, N) \ll \vartheta\theta^2(MN)^{\frac{3}{4}}(\log MN)^4$$

whereas the trivial bound is

(5.16) $$\mathcal{B}(M, N) \ll \theta^2(MN)^{\frac{3}{4}}(\log N)^2 \; .$$

Indeed, arguing along the same lines as for (5.10) we have

$$\begin{aligned} \mathcal{B}(M, N) &\ll \sum\sum_{M < u^2 + v^2 \leqslant 2M} \sum_{\substack{r + is \in \mathfrak{B} \\ ur + vs = \square}} |\beta(r^2 + s^2)| \\ &\ll M^{\frac{3}{4}} N^{-\frac{1}{4}} \sum_{r + is \in \mathfrak{B}} |\beta(r^2 + s^2)| \; . \end{aligned}$$

By Lemma 2.2
$$|\beta(r^2 + s^2)| \leqslant \tau(r^2 + s^2) \leqslant 9 \sum_{\substack{d|(r^2 + s^2) \\ d \leqslant N^{\frac{1}{3}}}} \tau(d) \; .$$

Given such a $d$, we have
$$\#\{r + is \in \mathfrak{B}; r^2 + s^2 \equiv 0(d)\} \ll \theta^2 N \rho(d) d^{-1} \; .$$

Moreover we have
$$\sum_{d \leqslant N^{\frac{1}{3}}} \rho(d)\tau(d)d^{-1} \ll (\log N)^2 \; .$$

These estimates yield the bound (5.16). With more work one could replace $(\log N)^2$ by $\log N$ but we seek the saving of a factor $\vartheta^{-1}$ which is an arbitrary power of $\log N$.

We can assume that

$$\left| \varphi \left(\mathrm{mod}\, \tfrac{\pi}{2}\right) \right| > \pi\vartheta$$

because the other sectors, altogether of angle $\leqslant 2\pi\vartheta$, contribute no more than the bound (4.23) by the estimate (5.16). For $z = r + is$ in any remaining polar box we are not near either axis, and hence

$$\vartheta\sqrt{N} < |r|,\ |s| < 2\sqrt{N}\ .$$

Be aware also that $r, s$ have fixed signs depending only on $\varphi$.

Recall that $\alpha_w$ are bounded numbers with $|w|^2$ in the dyadic segment (4.18). By Cauchy's inequality

$$(5.17) \qquad \mathcal{B}^2(M, N) \ll M\mathcal{D}(M, N)$$

where

$$(5.18) \qquad \mathcal{D}(M, N) = \sum_w f(w) \left| \sum_z \beta_z \mathfrak{Z}(\mathrm{Re}\,\overline{w}z) \right|^2\ .$$

Here we have introduced a smooth majorant $f$ to simplify the forthcoming harmonic analysis. We choose an $f$ that is supported in the annulus

$$(5.19) \qquad \tfrac{1}{2}\sqrt{M} \leqslant |w| \leqslant 2\sqrt{M}\ .$$

Also, it is convenient to take $f$ to be radial, that is $f(w) = f(|w|)$. Now we need to prove that

$$(5.20) \qquad \mathcal{D}(M, N) \ll \vartheta^2 \theta^4 M^{\frac{1}{2}} N^{\frac{3}{2}} (\log MN)^8\ .$$

Squaring out we get

$$(5.21) \qquad \mathcal{D}(M, N) = \sum_w f(w) \sum_{z_1} \sum_{z_2} \beta_{z_1} \overline{\beta}_{z_2} \mathfrak{Z}(\mathrm{Re}\,\overline{w}z_1) \mathfrak{Z}(\mathrm{Re}\,\overline{w}z_2)\ .$$

Here we want to insert the condition $(z_1, z_2) = 1$ which will give us a sum that is easier to work with. Since $z_1 z_2$ is coprime with $\Pi$ it will turn out, as we next show, that we can do this at a small cost.

First we require a trivial bound for $\mathcal{D}(M, N)$. To begin note that we have $|\mathcal{D}(M, N)| \leqslant \tau^2 D(M, N)$ where

$$D(M, N) = \sum_{|w|^2 \sim M}\ \sideset{}{^*}\sum_{|z_1|^2, |z_2|^2 \sim N} \sideset{}{^*}\sum \mathfrak{Z}(\mathrm{Re}\,\overline{w}z_1)\mathfrak{Z}(\mathrm{Re}\,\overline{w}z_2).$$

Here the $*$ indicates summation over primitive $z$.

LEMMA 5.1. *For every $M \geqslant N \geqslant 2$,*

$$D(M, N) \ll \left( M^{\frac{3}{4}} N^{\frac{3}{4}} + M^{\frac{1}{2}} N^{\frac{3}{2}} \right) (\log MN)^{514}\ .$$

*Proof.* The contribution of the diagonal $|z_1| = |z_2|$ is

$$D_=(M,N) \ll \sum_w \sum_z{}^* \mathfrak{Z}(\operatorname{Re}\overline{w}z) \ll (MN)^{\frac{3}{4}} \log MN$$

by the argument that gave (5.16). The remaining terms of $D(M,N)$, those off the diagonal, have

$$\Delta(z_1, z_2) = \tfrac{1}{2i}(\bar{z}_1 z_2 - z_1 \bar{z}_2) \neq 0$$

because $(z_1, \bar{z}_1) = (z_2, \bar{z}_2) = 1$ and $|z_1| \sim |z_2|$. These contribute

$$D_{\neq}(M,N) \ll \sum\sum_{|c_1|,|c_2| \leqslant 2(MN)^{\frac{1}{4}}} \sum^*\sum^*_{\substack{|z_1| \neq |z_2| \\ \Delta | (c_1^2 z_2 - c_2^2 z_1) \neq 0}} 1 \ .$$

That $c_1^2 z_2 - c_2^2 z_1 \neq 0$ follows from (6.3) and (6.4) below.

Using the rectangular coordinates $z_1 = r_1 + is_1$, $z_2 = r_2 + is_2$ we have $\Delta = r_1 s_2 - r_2 s_1$, $c_1^2 r_2 \equiv c_2^2 r_1 \pmod{|\Delta|}$ and $c_1^2 s_2 \equiv c_2^2 s_1 \pmod{|\Delta|}$. By symmetry we can assume $c_1^2 s_2 - c_2^2 s_1 \neq 0$. For given $c_1, c_2, s_1, s_2, \Delta \neq 0$, the number $r_1$ is fixed mod $s_1/(s_1, s_2)$ and then $r_2$ is determined. The number of pairs $r_1, r_2$ is bounded by $\sqrt{N}(s_1, s_2)/s_1$. Hence, letting $\delta = (s_1, s_2)$, $s_1 = \delta s_1^*$, $s_2 = \delta s_2^*$, we get

$$D_{\neq}(M,N) \ll \sqrt{N} \sum_\delta \tau(\delta) \sum_{s_1^*} \sum_{s_2^*} \frac{1}{s_1^*} \sum\sum_{c_1^2 s_2^* \neq c_2^2 s_1^*} \tau(c_1^2 s_2^* - c_2^2 s_1^*) \ .$$

By Lemma 2.2 with $k=4$ there exists $d \leqslant (8\sqrt{MN})^{\frac{1}{4}}$ such that we have $c_1^2 s_2^* \equiv c_2^2 s_1^* \pmod{d}$ and $\tau(c_1^2 s_2^* - c_2^2 s_1^*) \ll \tau(d)^8$. Hence

$$\sum\sum_{c_1^2 s_2^* \neq c_2^2 s_1^*} \tau(c_1^2 s_2^* - c_2^2 s_1^*) \ll \sum_{d < (MN)^{\frac{1}{4}}} \tau(d)^8 \sum\sum_{c_1^2 s_2^* \equiv c_2^2 s_1^* \pmod{d}} 1$$

$$\ll (MN)^{\frac{1}{2}} \sum_{d < (MN)^{\frac{1}{4}}} \tau(d)^9 d^{-1} \ll (MN)^{\frac{1}{2}} (\log MN)^{2^9}.$$

Summing over $s_1^*, s_2^*$ then $\delta$ we conclude that

$$D_{\neq}(M,N) \ll M^{\frac{1}{2}} N^{\frac{3}{2}} (\log MN)^{2+2^9} \ .$$

Combining this with the estimate for the diagonal contribution $D_=$ we complete the proof of Lemma 5.1. □

Now we are ready to reduce the sum $\mathcal{D}(M,N)$ to the corresponding sum restricted to $(z_1, z_2) = 1$. We denote the latter by

$$\mathcal{D}^*(M,N) = \sum_w f(w) \sum\sum_{(z_1, z_2) = 1} \beta_{z_1} \bar{\beta}_{z_2} \mathfrak{Z}(\operatorname{Re}\overline{w}z_1) \mathfrak{Z}(\operatorname{Re}\overline{w}z_2) \ .$$

We shall prove that the difference $\widetilde{\mathcal{D}}(M,N)$ between these two sums satisfies

$$\widetilde{\mathcal{D}}(M,N) \ll \tau^2 \left( M^{\frac{3}{4}} N^{\frac{3}{4}} + P^{-1} M^{\frac{1}{2}} N^{\frac{3}{2}} \right) (\log MN)^{516}.$$

Indeed, denoting Gaussian primes by $\pi$ we find that $\widetilde{\mathcal{D}}(M,N)$ is bounded by

$$\sum_{P<|\pi|^2\leqslant N} \sum_{M<|w|^2\leqslant 2M} \sideset{}{^*}\sum_{\frac{N}{|\pi|^2}<|z_1|^2,\ |z_2|^2\leqslant \frac{2N}{|\pi|^2}} \sideset{}{^*}\sum |\beta_{\pi z_1}\beta_{\pi z_2}|\mathfrak{z}(\operatorname{Re}\overline{w}\pi z_1)\mathfrak{z}(\operatorname{Re}\overline{w}\pi z_2)$$

$$\ll \tau^2 D(MP_1, NP_1^{-1})(\log MN)^2$$

for some $P_1$ with $P < P_1 \leqslant N$. The result now follows from Lemma 5.1.

Subtracting the above bound for $\widetilde{\mathcal{D}}(M,N)$ from $\mathcal{D}(M,N)$ we conclude that

$$(5.22)\quad \mathcal{D}(M,N) = \mathcal{D}^*(M,N) + O\left(\tau^2\left(M^{\frac{3}{4}}N^{\frac{3}{4}} + P^{-1}M^{\frac{1}{2}}N^{\frac{3}{2}}\right)(\log MN)^{516}\right).$$

Observe that the first error term in (5.22) is admissible for (5.20) provided that

$$(5.23)\qquad\qquad \tau \leqslant x^{\frac{1}{3}\eta},$$

by virtue of (4.19) and (4.6), and the second is admissible if

$$(5.24)\qquad\qquad P \geqslant \tau^2 (\log x)^{2A+4A'+508}.$$

Under these conditions (5.23) and (5.24) it remains to prove that

$$(5.25)\qquad\qquad \mathcal{D}^*(M,N) \ll \vartheta^2 \theta^4 M^{\frac{1}{2}} N^{\frac{3}{2}} (\log MN)^8 .$$

Changing the order of summation we arrange our new sum as

$$(5.26)\qquad\qquad \mathcal{D}^*(M,N) = \sum\sum_{(z_1,z_2)=1} \beta_{z_1}\overline{\beta}_{z_2}\mathcal{C}(z_1,z_2)$$

where

$$(5.27)\qquad\qquad \mathcal{C}(z_1,z_2) = \sum_w f(w)\mathfrak{z}(\operatorname{Re}\overline{w}z_1)\mathfrak{z}(\operatorname{Re}\overline{w}z_2) .$$

The last is a free sum over Gaussian integers $w$. Note that the restrictions on the support of $\beta_z$ which we have imposed so far and the summation condition $(z_1, z_2) = 1$ in (5.26) imply that $(z_1, \overline{z}_1) = (z_2, \overline{z}_2) = 1$ and $z_1 \equiv z_2 \pmod 8$.

## 6. Counting points inside a biquadratic ellipse

In this section we evaluate approximately $\mathcal{C}(z_1, z_2)$ defined by (5.27). The problem reduces to counting lattice points inside the curve

$$t_1^4 - 2\gamma t_1^2 t_2^2 + t_2^4 = x$$

for fixed $\gamma$ with $0 < \gamma < 1$, namely $\gamma = \cos(\alpha_2 - \alpha_1)$ where $\alpha_j = \arg z_j$. This particular curve arises because of the particular choice (3.14) of $\mathfrak{z}(b)$ as the Fourier coefficient of the classical theta function. The counting requires quite subtle harmonic analysis because the points involved are also constrained by a congruence to a large modulus, so there are very few points relative to the area of the region.

The modulus to which we have referred above is the determinant

$$(6.1) \qquad \Delta = \Delta(z_1, z_2) = \operatorname{Im} \overline{z}_1 z_2 = \frac{1}{2i}(\overline{z}_1 z_2 - z_1 \overline{z}_2) = |z_1 z_2| \sin(\alpha_2 - \alpha_1) \;.$$

Note that $(\Delta, |z_1 z_2|^2) = 1$ because $(z_1, z_2) = 1$ and $(z_1, \overline{z}_1) = (z_2, \overline{z}_2) = 1$. In particular, this co-primality condition implies that $\Delta$ does not vanish. For $z_1, z_2$ in the same polar box (5.14) we have

$$(6.2) \qquad 1 \leqslant |\Delta| < 4\pi \theta N \;.$$

The variable of summation $w$ in (5.27) which runs over Gaussian integers can be parameterized by two squares of rational integers;

$$(6.3) \qquad \operatorname{Re} \overline{w} z_1 = c_1^2, \qquad \operatorname{Re} \overline{w} z_2 = c_2^2 \;,$$

say, with $c_1, c_2 \in \mathbb{Z}$. Indeed these values determine $w$ by

$$i \Delta w = c_1^2 z_2 - c_2^2 z_1 \;.$$

As $w$ ranges freely over $\mathbb{Z}[i]$ the above equation is equivalent to the congruence

$$(6.4) \qquad c_1^2 z_2 \equiv c_2^2 z_1 \pmod{|\Delta|} \;.$$

Using the rectangular coordinates $z_1 = r_1 + is_1$ and $z_2 = r_2 + is_2$ it appears that (6.4) reads as two rational congruences

$$c_1^2 r_1 \equiv c_2^2 r_2 \pmod{|\Delta|}$$
$$c_1^2 s_1 \equiv c_2^2 s_2 \pmod{|\Delta|}$$

with

$$(6.5) \qquad \Delta = r_1 s_2 - r_2 s_1 \;.$$

In fact (6.4) actually reduces to a single congruence because $z_2/z_1$ is congruent to a rational number modulo $|\Delta|$,

$$(6.6) \qquad \frac{z_2}{z_1} \equiv \operatorname{Re} \frac{z_2}{z_1} = \frac{r_1 r_2 + s_1 s_2}{r_1^2 + s_1^2} \pmod{|\Delta|} \;.$$

Applying the above transformations we write (5.27) as

$$(6.7) \qquad \mathcal{C}(z_1, z_2) = \sum\sum_{c_1^2 z_2 \equiv c_2^2 z_1 \,(\operatorname{mod} |\Delta|)} f\left((c_1^2 z_2 - c_2^2 z_1)/\Delta\right) \;.$$

Next we split the summation into residue classes modulo $|\Delta|$, say

$$\mathcal{C}(z_1, z_2) = \sum\sum_{\gamma_1^2 z_2 \equiv \gamma_2^2 z_1 (\bmod |\Delta|)} \mathcal{C}(z_1, z_2; \gamma_1, \gamma_2) ,$$

where

$$\mathcal{C}(z_1, z_2; \gamma_1, \gamma_2) = \sum\sum_{(c_1, c_2) \equiv (\gamma_1, \gamma_2)(\bmod |\Delta|)} f\left((c_1^2 z_2 - c_2^2 z_1)/\Delta\right) .$$

Then for each pair of classes we execute the summation by Poisson's formula obtaining for $\mathcal{C}(z_1, z_2; \gamma_1, \gamma_2)$ the Fourier series

$$|\Delta|^{-1} |z_1 z_2|^{-\frac{1}{2}} \sum_{h_1} \sum_{h_2} F\left(h_1 |\Delta z_2|^{-\frac{1}{2}}, h_2 |\Delta z_1|^{-\frac{1}{2}}\right) e\left((\gamma_1 h_1 + \gamma_2 h_2)|\Delta|^{-1}\right),$$

where

$$(6.8) \qquad F(u_1, u_2) = \iint f\left(\frac{z_2}{|z_2|} t_1^2 - \frac{z_1}{|z_1|} t_2^2\right) e(u_1 t_1 + u_2 t_2) dt_1\, dt_2 .$$

Hence we obtain

$$(6.9) \qquad \mathcal{C}(z_1, z_2) = |z_1 z_2|^{-\frac{1}{2}} \sum_{h_1} \sum_{h_2} F\left(h_1 |\Delta z_2|^{-\frac{1}{2}},\ h_2 |\Delta z_1|^{-\frac{1}{2}}\right) G(h_1, h_2)$$

where $G(h_1, h_2)$ is the sum over the residue classes $\gamma_1, \gamma_2$ to modulus $|\Delta|$,

$$(6.10) \qquad G(h_1, h_2) = \frac{1}{|\Delta|} \sum\sum_{\gamma_1^2 z_2 \equiv \gamma_2^2 z_1 (\bmod |\Delta|)} e\left((\gamma_1 h_1 + \gamma_2 h_2)|\Delta|^{-1}\right) .$$

Have in mind that $F(u_1, u_2)$ and $G(h_1, h_2)$ depend also on $z_1, z_2$. Naturally the main contribution comes from $h_1 = h_2 = 0$. In this case we display the dependence on $z_1, z_2$ by writing

$$(6.11) \qquad F_0(z_1, z_2) = \iint f\left(\frac{z_2}{|z_2|} t_1^2 - \frac{z_1}{|z_1|} t_2^2\right) dt_1\, dt_2$$

and

$$(6.12) \qquad G_0(z_1, z_2) = \frac{1}{|\Delta|} \#\left\{\gamma_1, \gamma_2;\ \gamma_1^2 z_2 \equiv \gamma_2^2 z_1 (\bmod |\Delta|)\right\}$$

which stand for $F(0,0)$ and $G(0,0)$ respectively. Thus the contribution of the zero frequencies to $\mathcal{C}(z_1, z_2)$ is

$$(6.13) \qquad \mathcal{C}_0(z_1, z_2) = |z_1 z_2|^{-\frac{1}{2}} F_0(z_1, z_2) G_0(z_1, z_2) .$$

In the next two sections we compute the Fourier integral $F(u_1, u_2)$ and the exponential sum $G(h_1, h_2)$.

## 7. The Fourier integral $F(u_1, u_2)$

Recall that $f$ is a radial function, say
$$f(w) = f(|w|) = \mathfrak{f}(|w|^2) ,$$
where $\mathfrak{f}$ is a smooth function supported on $[\frac{1}{4}M, 4M]$. We shall assume that $|\mathfrak{f}^{(j)}| \leqslant M^{-j}$ for $0 \leqslant j \leqslant 4$. Putting
$$g(t_1, t_2) = \left| \frac{z_2}{|z_2|} t_1^2 - \frac{z_1}{|z_1|} t_2^2 \right|^2$$
we have
$$(7.1) \qquad F(u_1, u_2) = \iint \mathfrak{f}(g(t_1, t_2)) \, e(u_1 t_1 + u_2 t_2) \, dt_1 \, dt_2 .$$
Here the domain of integration is restricted by the support of $\mathfrak{f}(g)$. We have
$$g(t_1, t_2) = t_1^4 - 2 \operatorname{Re} \frac{\overline{z_1} z_2}{|z_1 z_2|} t_1^2 t_2^2 + t_2^4$$
and
$$\frac{\overline{z_1} z_2}{|z_1 z_2|} = \cos(\alpha_2 - \alpha_1) + i \sin(\alpha_2 - \alpha_1) = \gamma + i\delta ,$$
say, where $\alpha_j = \arg z_j$. Note that $\gamma$ is close to 1 and $|\delta|$ is small because $z_1, z_2$ are in the same polar box (5.14). Precisely
$$(7.2) \qquad \delta = \Delta |z_1 z_2|^{-1} = \sin(\alpha_2 - \alpha_1) ,$$
so that
$$(7.3) \qquad (2N)^{-1} < |\delta| < 4\pi\theta .$$
Applying the above notation we write several useful expressions for the quartic form $g(t_1, t_2)$:
$$g(t_1, t_2) = t_1^4 - 2\gamma t_1^2 t_2^2 + t_2^4$$
$$= (t_1^2 - \gamma t_2^2)^2 + \delta^2 t_2^4 = (t_2^2 - \gamma t_1^2)^2 + \delta^2 t_1^4.$$
Hence
$$4g(t_1, t_2) \geqslant 2\delta^2(t_1^4 + t_2^4) \geqslant \delta^2(t_1^2 + t_2^2)^2 .$$
We have also
$$(t_1^2 - t_2^2)^2 \leqslant g(t_1, t_2) \leqslant (t_1^2 + t_2^2)^2 .$$
Since $\frac{1}{4}M < g < 4M$ by the support of $\mathfrak{f}$ these inequalities imply
$$|t_1^2 - t_2^2| < 2M^{\frac{1}{2}} ,$$
$$\tfrac{1}{2}M^{\frac{1}{2}} < t_1^2 + t_2^2 < 4M^{\frac{1}{2}} |\delta|^{-1} .$$

Hence the area of integration (7.1) is bounded by $O\left(M^{\frac{1}{2}} \log \delta^{-2}\right)$ so by trivial estimation

$$F(u_1, u_2) \ll M^{\frac{1}{2}} \log N. \tag{7.4}$$

If $u_1 \neq 0$ we can integrate by parts four times with respect to $t_1$ getting

$$F(u_1, u_2) = (2\pi u_1)^{-4} \iint \frac{\partial^4}{\partial t_1^4} \mathfrak{f}(g(t_1, t_2))\, e\,(u_1 t_1 + u_2 t_2)\, dt_1\, dt_2\ .$$

We have

$$\frac{\partial^4}{\partial t_1^4}\mathfrak{f}(g) = g''''\mathfrak{f}' + \left(4g'g''' + 3g''g''\right)\mathfrak{f}'' + 6g'g'g''\mathfrak{f}''' + g'g'g'g'\mathfrak{f}''''\ ,$$

and by the above inequalities this is bounded by

$$M^{-1} + (t_1^2 + t_2^2)^2 M^{-2} + t_1^2(t_1^2 + t_2^2)(t_1^2 - \gamma t_2^2)^2 M^{-3} + t_1^4(t_1^2 - \gamma t_2^2)^4 M^{-4}\ .$$

Since

$$(t_1^2 - \gamma t_2^2)^2 + (t_2^2 - \gamma t_1^2)^2 + \delta^2(t_1^4 + t_2^4) = 2g(t_1, t_2) < 8M$$

we deduce that

$$\frac{\partial^4}{\partial t_1^4}\mathfrak{f}(g) \ll \delta^{-2} M^{-1}\ .$$

This gives, by trivial estimation combined with (7.4),

$$F(u_1, u_2) \ll \left(1 + u_1^4 \delta^2 M\right)^{-1} M^{\frac{1}{2}} \log N\ .$$

By symmetry this bound also holds with $u_1$ replaced by $u_2$. Taking the geometric mean of these two bounds we arrive at

$$F(u_1, u_2) \ll \left(1 + u_1^2 |\delta| \sqrt{M}\right)^{-1} \left(1 + u_2^2 |\delta| \sqrt{M}\right)^{-1} M^{\frac{1}{2}} \log N\ .$$

Hence we get by (7.2):

LEMMA 7.1.  *For $u_1 = h_1 |\Delta z_2|^{-1/2}$ and $u_2 = h_2 |\Delta z_1|^{-1/2}$ the Fourier integral* (6.8) *satisfies*

$$F(u_1, u_2) \ll \left(1 + h_1^2 H^{-2}\right)^{-1} \left(1 + h_2^2 H^{-2}\right)^{-1} M^{\frac{1}{2}} \log N \tag{7.5}$$

*where $H = M^{-\frac{1}{4}} N^{\frac{3}{4}}$.*

In particular, for $h_1 = h_2 = 0$ the estimate (7.5) becomes

$$F_0(z_1, z_2) \ll M^{\frac{1}{2}} \log N\ , \tag{7.6}$$

but in this case we need a more precise formula. We compute as follows

$$F_0(z_1, z_2) = \int_{-\infty}^{\infty} \int_{-\infty}^{\infty} \mathfrak{f}(g(t_1, t_2))\, dt_1\, dt_2$$

$$= 4 \int_0^{\infty} \int_0^{\infty} f\left((t_1^4 - 2\gamma t_1^2 t_2^2 + t_2^4)^{\frac{1}{2}}\right) dt_1\, dt_2$$

$$= 4 \int_0^{\infty} \int_0^{\infty} f\left(u^2(t^4 - 2\gamma t^2 + 1)^{\frac{1}{2}}\right) u\, du\, dt$$

$$= \hat{f}(0) E(\gamma),$$

where $\hat{f}(0)$ is the integral mean-value of $f$,

(7.7) $$\hat{f}(0) = \int_0^{\infty} f(u)\, du \ll M^{\frac{1}{2}}$$

and $E(\gamma)$ denotes the elliptic integral

(7.8) $$E(\gamma) = \int_0^{\infty} (t^2 - 2\gamma t + 1)^{-\frac{1}{2}} t^{-\frac{1}{2}} dt\ .$$

Since $\gamma$ is close to 1 we have a satisfactory asymptotic expansion (cf. (3.138.7) and (8.113.3) of [GR])

(7.9) $$E(\gamma) = \log 4\delta^{-2} + O(\delta^2 \log \delta^{-2})\ .$$

Insertion of (7.2) in this gives

LEMMA 7.2.    For $z_1$, $z_2$ in the box (5.14) the integral (6.11) satisfies

(7.10) $$F_0(z_1, z_2) = \hat{f}(0) 2\log 2|z_1 z_2/\Delta| + O\left(\Delta^2 M^{\frac{1}{2}} N^{-2} \log N\right)\ .$$

## 8. The arithmetic sum $G(h_1, h_2)$

Recall that $G(h_1, h_2)$ is given by (6.10). This is a kind of Weyl sum for the equidistribution of roots $\gamma_1, \gamma_2$ of the quadratic form $\gamma_1^2 z_2 - \gamma_2^2 z_1$ modulo $|\Delta|$. We write (uniquely)

(8.1) $$\Delta = \Delta_1 \Delta_2^2$$

where $\Delta_1$ is squarefree and $\Delta_2 \geqslant 1$. The solutions to

(8.2) $$\gamma_1^2 z_2 \equiv \gamma_2^2 z_1 \pmod{|\Delta|}$$

satisfy $(\gamma_1^2, \Delta) = (\gamma_2^2, \Delta) = d_1 d_2^2$, say, where $d_1, d_2 \geqslant 1$ and $d_1$ is squarefree. This implies $d_1 | \Delta_1$, $d_2 | \Delta_2$ and $(d_1, \Delta_2/d_2) = 1$. Moreover $\gamma_1 = d_1 d_2 \eta_1$, $\gamma_2 = d_1 d_2 \eta_2$ where $\eta_1, \eta_2$ run over the residue classes to modulus $|\Delta|/d_1 d_2$ and coprime with $|\Delta|/d_1 d_2^2$. Accordingly $G(h_1, h_2)$ splits into

$$G(h_1, h_2) = |\Delta|^{-1} \sum_{\substack{b_1 d_1 = |\Delta_1| \\ b_2 d_2 = \Delta_2 \\ (d_1, b_2) = 1}} \sum_{\substack{\eta_1, \eta_2 \pmod{b_1 b_2^2 d_2} \\ (\eta_1 \eta_2, b_1 b_2) = 1 \\ \eta_1^2 z_2 \equiv \eta_2^2 z_1 \pmod{b_1 b_2^2}}} e\left((\eta_1 h_1 + \eta_2 h_2)/b_1 b_2^2 d_2\right) .$$

The innermost sum vanishes unless $h_1 \equiv h_2 \equiv 0 \pmod{d_2}$ so the full sum is equal to

$$|\Delta|^{-1} \sum_{\substack{b_1 d_1 = |\Delta_1| \\ (d_1, b_2) = 1}} \sum_{\substack{b_2 d_2 = \Delta_2 \\ d_2 | (h_1, h_2)}} d_2^2 \sideset{}{^*}\sum_{\substack{\eta_1, \eta_2 \pmod{b_1 b_2^2} \\ \eta_1^2 z_2 \equiv \eta_2^2 z_1 \pmod{b_1 b_2^2}}} \sideset{}{^*}\sum e\left((\eta_1 h_1 + \eta_2 h_2)/b_1 b_2^2 d_2\right) .$$

Changing $\eta_2$ into $\omega \eta_1$ we conclude that $G(h_1, h_2)$ is given by

(8.3) $\quad |\Delta|^{-1} \sum_{\substack{b_1 d_1 = |\Delta| \\ (d_1, b_2) = 1}} \sum_{\substack{b_2 d_2 = \Delta_2 \\ d_2 | (h_1, h_2)}} d_2^2 \sum_{\omega^2 \equiv z_2/z_1 \pmod{b_1 b_2^2}} R\left((h_1 + \omega h_2) d_2^{-1}; b_1 b_2^2\right)$

where $R(h; b)$ denotes the Ramanujan sum

$$R(h; b) = \sideset{}{^*}\sum_{\eta \pmod b} e\left(\frac{\eta h}{b}\right) .$$

Using the well-known bound $|R(h; b)| \leqslant (h, b)$, we obtain

$$\begin{aligned}
|R\left((h_1 + \omega h_2) d_2^{-1}; b_1 b_2^2\right)| &\leqslant \left((h_1 + \omega h_2) d_2^{-1}, b_1 b_2^2\right) \\
&\leqslant \left((h_1^2 - \omega^2 h_2^2) d_2^{-2}, b_1 b_2^2\right) \\
&\leqslant \left(z_1 h_1^2 - z_2 h_2^2, \Delta\right) d_2^{-2}.
\end{aligned}$$

Denote by $n(z; b)$ the number of solutions in rational integer classes $\omega$ modulo $b$ of the quadratic congruence

(8.4) $\qquad\qquad\qquad\qquad \omega^2 \equiv z \pmod{b}.$

Incidentally notice that for $z \equiv -a^2 \pmod{b}$ we get the arithmetic function $n(-a^2; b) = \rho(a; b)$ which was considered in Section 3. Of course, $n(z; b)$ vanishes if $z$ is not congruent to a rational integer, however in our case $z = z_2/z_1$ is rational modulo $|\Delta|$ and prime to $\Delta$; see (6.6). Applying the trivial bound $n(z_2/z_1; b_1 b_2^2) \leqslant 4\tau(b_1 b_2)$ together with the above bound for Ramanujan's sum we deduce by (8.3) that

LEMMA 8.1.  *For any $h_1, h_2$ the exponential sum (6.10) satisfies*

(8.5) $\qquad |G(h_1, h_2)| \leqslant 4\tau_3(\Delta) |\Delta|^{-1} (z_1 h_1^2 - z_2 h_2^2, \Delta) .$

In particular for $h_1 = h_2 = 0$ the estimate (8.5) becomes

(8.6) $\qquad\qquad\qquad\qquad G_0(z_1, z_2) \leqslant 4\tau_3(\Delta)$

which is almost best possible but in this case we need an exact formula. Since $R(0; b) = \varphi(b)$ we obtain by (8.3)

$$(8.7) \qquad G_0(z_1, z_2) = \sum_{b_1 d_1 = |\Delta_1|} d_1^{-1} \sum_{\substack{b_2 d_2 = \Delta_2 \\ (b_2, d_1) = 1}} (b_1 b_2^2)^{-1} \varphi(b_1 b_2^2) n(z_2/z_1; b_1 b_2^2) \ .$$

On the other hand by (6.12) we have $|\Delta| G_0(z_1, z_2) = N(z_2/z_1; |\Delta|)$ where $N(a; q)$ denotes the number of solutions to

$$(8.8) \qquad a\gamma_1^2 \equiv \gamma_2^2 \pmod{q} \ .$$

Since $N(a; q)$ is multiplicative we have

$$(8.9) \qquad G_0(z_1, z_2) = \prod_{p^\nu \| \Delta} p^{-\nu} N(z_2/z_1; p^\nu) \ .$$

This expression reduces the problem to local computations.

LEMMA 8.2. *For $p \neq 2$ and $(a, p) = 1$ we have*

$$(8.10) \qquad p^{-\nu} N(a; p^\nu) = 1 + \left(1 - \frac{1}{p}\right) \left(\left[\frac{\nu}{2}\right] + \left[\frac{\nu+1}{2}\right] \left(\frac{a}{p}\right)\right) \ .$$

*For $p = 2$, $\nu \geq 1$, $a \equiv 1 \pmod 8$ we have*

$$(8.11) \qquad 2^{-\nu} N(a; 2^\nu) = \nu \ .$$

*Proof.* One could proceed by counting the solutions to (8.8) directly, however we use the formula (8.7). This gives

$$p^{-\nu} N(a; p^\nu) = \frac{1}{2p}(1 - (-1)^\nu) + \sum_{\substack{0 \leq \alpha \leq \nu \\ \alpha \equiv \nu \pmod 2}} p^{-\alpha} \varphi(p^\alpha) n(a; p^\alpha)$$

where the first term is present only if $2 \nmid \nu$ and it comes from $d_1 = p$ in (8.7). For $p \neq 2$ and $\alpha \geq 1$ we have

$$(8.12) \qquad n(a; p^\alpha) = 1 + \left(\frac{a}{p}\right)$$

which leads to (8.10). For $p = 2$ and $a \equiv 1 \pmod 8$ we have $n(a; p^\alpha) = (4, 2^{\alpha-1})$ yielding (8.11). □

In the formula (8.10) we write

$$\left(\left[\frac{\nu}{2}\right] + \left[\frac{\nu+1}{2}\right] \left(\frac{a}{p}\right)\right) = \left(\frac{a}{p}\right) + \left(\frac{a}{p^2}\right) + \cdots + \left(\frac{a}{p^\nu}\right)$$

which leads to the following global expression:

COROLLARY 8.3. *For $q$ odd and $(a, q) = 1$ we have*

$$(8.13) \qquad N(a; q) = q \sum_{d|q} \frac{\varphi(d)}{d} \left(\frac{a}{d}\right) .$$

By (8.9), (8.11), and (8.13) we infer that

$$(8.14) \qquad G_0(z_1, z_2) = \nu \sum_{\substack{d|\Delta \\ d \text{ odd}}} \frac{\varphi(d)}{d} \left(\frac{z_2/z_1}{d}\right)$$

where $\nu$ is the order of 2 in $\Delta$, that is $\Delta = 2^\nu \Delta'$ with $\Delta'$ odd. Note that $z_2/z_1 \equiv 1 \pmod 8$ due to (5.7), so $\nu \geqslant 3$.

To accommodate the factor $\nu$ we extend the Jacobi symbol to even moduli by setting

$$(8.15) \qquad \left(\frac{a}{d}\right) = \left(\frac{a}{d'}\right) \quad \text{if } 2 \nmid a ,$$

where $d'$ denotes the odd part of $d$. Now we conclude by (8.14)

LEMMA 8.4. *Suppose that $(z_1, z_2) = (z_1, \overline{z}_1) = (z_2, \overline{z}_2) = 1$ and also $z_1 \equiv z_2 \pmod 8$. Then the number of congruence pairs of solutions in (6.12) with given $\Delta = \Delta(z_1, z_2) = \operatorname{Im} \overline{z}_1 z_2$ is expressed by*

$$(8.16) \qquad G_0(z_1, z_2) = 2 \sum_{4d|\Delta} \frac{\varphi(d)}{d} \left(\frac{z_2/z_1}{d}\right) .$$

## 9. Bounding the error term in the lattice point problem

In this section we combine the results of the previous two sections completing the estimation of the error term in the lattice point problem of Section 6. The estimation of the main term will take the rest of the paper.

By (6.9), (7.5) and (8.5) we obtain

$$(9.1) \qquad \mathcal{C}(z_1, z_2) = \mathcal{C}_0(z_1, z_2) + O\left(\tau_3(\Delta)|\Delta|^{-1} \mathcal{H}(z_1, z_2) M^{\frac{1}{2}} N^{-\frac{1}{2}} \log N\right)$$

where $\mathcal{C}_0(z_1, z_2)$ is given by (6.13) and

$$(9.2) \quad \mathcal{H}(z_1, z_2) = \sum\sum_{(h_1, h_2) \neq (0,0)} (z_1 h_1^2 - z_2 h_2^2, \Delta)(1 + h_1^2 H^{-2})^{-1}(1 + h_2^2 H^{-2})^{-1} .$$

Recall that $H = M^{-\frac{1}{4}} N^{\frac{3}{4}}$. For aesthetic reasons only we would like to estimate $\mathcal{H}(z_1, z_2)$ for individual $z_1, z_2$, however the effective range of summation is too short to do so. Note that for $N$ in (4.6) and $MN$ satisfying (4.19) we have

$$(9.3) \qquad x^\eta < H < N^{\frac{1}{2}} .$$

Thus $h_1, h_2$ are small indeed, nevertheless $(z_1h_1^2 - z_2h_2^2, \Delta)$ can be quite large for some points $z_1, z_2$. For this reason we take advantage of the additional summation over $z_1, z_2$ which is present in our main problem.

Given $h_1, h_2$ not both zero we put

$$\Lambda = \Lambda(z_1, z_2) = z_1h_1^2 - z_2h_2^2 . \tag{9.4}$$

We begin by estimating the sum

$$\mathcal{Z}(h_1, h_2) = \sum\sum_{(z_1,z_2)=1} |\beta_{z_1}\beta_{z_2}|\tau_3(\Delta)(\Lambda, \Delta)|\Delta|^{-1} \tag{9.5}$$

where $\Delta = \Delta(z_1, z_2)$ is the determinant defined in (6.1). We have $1 \leqslant |\Delta| < N$ and, because of (5.9),

$$|\beta_z| \leqslant \tau(|z|^2) \leqslant \tau. \tag{9.6}$$

Hence

$$\mathcal{Z}(h_1, h_2) \leqslant 2\tau^2(\log N)D^{-1} \sum\sum_{\substack{z_1,z_2 \in \mathfrak{B} \\ D \leqslant |\Delta| < 2D}} \tau_3(\Delta)(\Lambda, \Delta)$$

for some $1 \leqslant D \leqslant N$ where $\mathfrak{B}$ is the polar box (5.14). Next we group terms according to the value of $(\Lambda, \Delta) = d$ say, getting

$$\mathcal{Z}(h_1, h_2) \leqslant 2\tau^2(\log N)D^{-1} \sum_{d<2D} d \sum\sum_{\substack{z_1,z_2 \in \mathfrak{B} \\ D \leqslant |\Delta| < 2D \\ \Delta \equiv \Lambda \equiv 0 (\bmod\, d)}} \tau_3(\Delta) .$$

For further computations we use the rectangular coordinates $z_1 = r_1 + is_1$ and $z_2 = r_2 + is_2$ with $r_1, r_2, s_1, s_2$ satisfying (5.4), (5.5). Observe the relations

$$\Lambda(r_1, r_2) \equiv \Lambda(s_1, s_2) \equiv 0 (\bmod\, d),$$

$$\Lambda(r_1, r_2)s_2 - \Lambda(s_1, s_2)r_2 = \Delta(z_1, z_2)h_1^2,$$

$$\Lambda(r_1, r_2)s_1 - \Lambda(s_1, s_2)r_1 = \Delta(z_1, z_2)h_2^2.$$

Since $\Delta(z_1, z_2) \neq 0$ and $h_1^2 + h_2^2 \neq 0$ we have either $\Lambda(r_1, r_2)$ or $\Lambda(s_1, s_2)$ different from zero, thus we may assume that

$$0 \neq \Lambda(r_1, r_2) \equiv 0 (\bmod\, d).$$

Given $r_1, r_2$ and any value of $\Delta$ the number of ordered pairs $s_1, s_2$ which give $z_1, z_2$ in $\mathfrak{B}$ satisfying $r_1s_2 - r_2s_1 = \Delta$ is $O(\vartheta^{-1}(r_1, r_2))$. Moreover we have

$$\sum_{\substack{D \leqslant |\Delta|<2D \\ \Delta \equiv 0(\bmod\, d)}} \tau_3(\Delta) \ll \tau_3(d)d^{-1}D(\log 2D)^2 .$$

Therefore
$$\mathcal{Z}(h_1, h_2) \ll \vartheta^{-1}\tau^2(\log N)^3 \sum_{\substack{d \\ 0 \neq \Lambda(r_1,r_2) \equiv 0(d)}} \sum_{r_1} \sum_{r_2} \tau_3(d)(r_1, r_2)$$
$$= \vartheta^{-1}\tau^2(\log N)^3 \sum\sum_{\Lambda(r_1,r_2) \neq 0} (r_1, r_2)\tau_4(\Lambda(r_1, r_2)).$$

Once again we have a difficulty to complete the summation, but this time in the variables $r_1, r_2$ because of the divisor function $\tau_4(\Lambda(r_1, r_2))$, especially when $\rho = (r_1, r_2)$ is large. For this reason we return to the summation in $h_1, h_2$. We need to estimate
$$L(r_1, r_2) = \sum\sum_{r_1 h_1^2 \neq r_2 h_2^2} \tau_4(r_1 h_1^2 - r_2 h_2^2)(1 + h_1^2 H^{-2})^{-1}(1 + h_2^2 H^{-2})^{-1}$$

where $H$ satisfies (9.3). Note that we can restrict this series to $h_1, h_2 \ll H^3$ because the contribution of the other terms is absorbed by the term with $h_1 = 1$, $h_2 = 0$. In this truncated series we estimate $\tau_4(r_1 h_1^2 - r_2 h_2^2)$ by $\tau(q)^c$ for some $q \leq H$ with $r_1 h_1^2 \equiv r_2 h_2^2 \pmod{q}$, where $c$ is a constant depending only on $\eta$. Specifically we may use Lemma 2.2 for $n = |r_1 h_1^2 - r_2 h_2^2|$ and note that $n < 4\sqrt{N}H^6 < H^{1/2\eta}$ by (9.3) so it suffices to take $c = \eta^{-1}\log\eta^{-1}$. Therefore
$$L(r_1, r_2) \ll \sum_{q \leq H} \tau(q)^c \sum\sum_{r_1 h_1^2 \equiv r_2 h_2^2 \pmod{q}} (1 + h_1^2 H^{-2})^{-1}(1 + h_2^2 H^{-2})^{-1}$$

where the restriction to $h_1, h_2 \ll H^3$ is no longer required. Splitting into residue classes to modulus $q$ we get

(9.7) $$L(r_1, r_2) \ll H^2 \sum_{q \leq H} \tau(q)^c q^{-2} N(r_1, r_2; q)$$

where $N(a, b; q)$ denotes the number of solutions $\gamma_1, \gamma_2$ to

(9.8) $$a\gamma_1^2 \equiv b\gamma_2^2 \pmod{q}.$$

If $(b, q) = 1$ then $N(a, b; q)$ is equal to $N(ab; q)$ and the latter was evaluated in Corollary 8.3 in the case $(2ab, q) = 1$. Now we give a general estimate.

LEMMA 9.1. *For any integers $a, b, q$ with $q \geq 1$ we have*

(9.9) $$N(a, b; q) \leq ([a, b], q)q\tau(q).$$

*Proof.* By multiplicativity we can assume that $q$ is a prime power, say $q = p^\nu$. For $q$ prime the bound is trivial if $p | ab$, while if not it reduces for each $\gamma_2$ to the congruence $\gamma_1^2 \equiv r \pmod{q}$ which has at most two solutions. From now on suppose $\nu \geq 2$. If $p \mid a$ and $p \mid b$ then we reduce to the case $q = p^{\nu-1}$ by dividing through by $p$. If $p \mid a$ and $p \nmid b$ then $p^2 \mid a$ so we can divide by $p^2$ and reduce to the case $q = p^{\nu-2}$ getting (by induction) a bound $p^{-4}(ab,q)q\tau(qp^{-2})$ for the solutions modulo $qp^{-2}$. Multiplying this bound by $p^4$ we get the result.

There remains the case $p \nmid ab$. If $p \mid \gamma_2$ then also $p \mid \gamma_1$ and these give by induction a contribution no more than

$$qp^{-2}\tau(qp^{-2})p^2 = q\tau(q) - 2q .$$

The other solutions satisfy $p \nmid \gamma_1\gamma_2$ and for each of the $\varphi(q)$ values of $\gamma_2$ there are at most two values of $\gamma_1$ except for $q$ even in which case there are at most four. Therefore the primitive solutions contribute at most $2\varphi(q) < 2q$ for $q$ odd and $4\varphi(q) = 2q$ for $q$ even. Adding these contributions we obtain (9.9). □

Inserting (9.9) into (9.7) we get

$$L(r_1, r_2) \ll H^2 \sum_{q \leqslant H} q^{-1}\tau(q)^{c+1}([r_1, r_2], q) \ll \tau([r_1, r_2])^{c+2} H^2 (\log H)^{2^{c+1}} .$$

Moreover we have

$$\sum_{r_1}\sum_{r_2}(r_1, r_2)\tau([r_1, r_2])^{c+2} \leqslant \sum_\rho \rho\tau(\rho)^{c+2}\left(\sum_{\rho r < 2\sqrt{N}}\tau(r)^{c+2}\right)^2$$
$$\ll N(\log N)^{2^{c+4}}.$$

From these estimates we conclude by (5.26) and (9.1) that

(9.10) $$\mathcal{D}^*(M, N) = \mathcal{D}_0(M, N) + O\left(\vartheta^{-1}\tau^2 N^2 (\log N)^b\right)$$

for some $b$ depending on $\eta$; precisely $b = 2^{c+4} + 2^{c+1} + 4 \ll \eta^{-1/\eta}$ where $\eta$ is fixed in Proposition 4.1. Here the main term is

(9.11) $$\mathcal{D}_0(M, N) = \sum\sum_{(z_1, z_2)=1} \beta_{z_1}\overline{\beta}_{z_2} \mathcal{C}_0(z_1, z_2),$$

and the error term is, by (4.6), (4.15), (4.16), and (4.19), admissible for (5.25) provided that

(9.12) $$\tau \leqslant (\log x)^{\frac{1}{2}B - \frac{7}{4}A - 2A' - \frac{1}{2}b}.$$

## 10. Breaking up the main term

It remains to estimate the main term (9.11). Recall that $\mathcal{C}_0(z_1, z_2)$ is given by (6.13) in terms of the Fourier integral (6.11) and the arithmetic sum (6.12). Although we refer to the sum $\mathcal{D}_0(M, N)$ as the "main term" since it originated from the leading term in the lattice point problem it is in fact, in contrast to what is usually called by that name, smaller than it would appear at first glance. This however is for reasons more subtle than those responsible for the size of what we have called the error term.

The lacunarity of the original sequence $\mathcal{A} = (a_n)$ has featured in the estimation for the error term in (9.10), and it is no longer an issue in the main term $\mathcal{D}_0(M, N)$. Thus it is quite easy to derive by trivial estimation the bound

$$(10.1) \qquad \mathcal{D}_0(M, N) \ll \theta^4 M^{\frac{1}{2}} N^{\frac{3}{2}} (\log N)^8$$

which of course barely misses what one wants, namely (5.20). The required improvement by a factor $\vartheta^2$ (actually a saving of arbitrary power of $\log N$) will result from the cancellation of terms due to the sign change of $\beta_{z_1} \bar{\beta}_{z_2}$ which involve the Möbius function. These are twisted by the arithmetic kernel $G_0(z_1, z_2)$ which is rather sophisticated; it involves the Jacobi symbol

$$(10.2) \qquad \chi_d(z_2/z_1) = \left( \frac{z_2/z_1}{d} \right)$$

which in turn originated (in Section 8) from the roots to the congruence

$$(10.3) \qquad \omega^2 \equiv z_2/z_1 \pmod{|\Delta|} .$$

Had the twist of $\beta_{z_1} \bar{\beta}_{z_2}$ been the smooth function $F_0(z_1, z_2)$ alone, or for that matter a separable arithmetic function of suitable character, we would be able to receive the cancellation quickly from an excursion into the zero-free region for Hecke $L$-functions with Grossencharacters. However the presence of the symbol $\chi_d(z_2/z_1)$ to very large moduli $d$ obscures the situation and we need modern technology to resolve the problems. Among our arguments one can find some traces of automorphic theory but we do not dwell on these here. Retrospectively, the treatment of the main term $\mathcal{D}_0(M, N)$ should be regarded as the core of our proof of the main Theorem 1, though it did not seem a central matter when we got the very first ideas.

In this section we relax some factors in the main term by using familiar approximations and we break up what is left into three parts according to the size of the moduli (small, medium, large) to be treated separately by three distinct methods in the forthcoming sections.

Inserting (6.13) into (9.11) we get by the approximations (7.10) and (8.6)

$$\mathcal{D}_0(M, N) = 2\hat{f}(0) \sum\sum_{(z_1, z_2)=1} \beta_{z_1} \bar{\beta}_{z_2} G_0(z_1, z_2) |z_1 z_2|^{-\frac{1}{2}} \log 2|z_1 z_2/\Delta|$$
$$+ O\left( M^{\frac{1}{2}} N^{-\frac{5}{2}} (\log N) \sum\sum_{(z_1, z_2)=1} |\beta_{z_1} \bar{\beta}_{z_2}| \tau_3(\Delta) \Delta^2 \right).$$

Recall that $\beta_z$ are supported on primary numbers in a polar box (5.14) and are restricted by (5.8), (5.9). We no longer need, nor wish to have, this last condition $\tau(|z|^2) \leqslant \tau$. We remove this in the same way as we have installed it between (4.9) – (4.11). We also simplify slightly by inserting
$$|z_1 z_2|^{-\frac{1}{2}} = (1 + O(\theta)) N^{-\frac{1}{2}}$$
(see (5.14)), and we use (6.2), (7.7), (8.6) to arrive at

(10.4) $\quad \mathcal{D}_0(M,N) = 2\hat{f}(0)N^{-\frac{1}{2}}T(\beta) + O\left((\tau^{-1} + \theta)Y(\beta)M^{\frac{1}{2}}N^{-\frac{1}{2}}\log N\right)$ .

Here

(10.5) $$T(\beta) = \sum\sum_{(z_1,z_2)=1} \beta_{z_1}\overline{\beta}_{z_2} G_0(z_1, z_2) \log 2|z_1 z_2/\Delta|$$

and

(10.6) $$Y(\beta) = \sum\sum_{(z_1,z_2)=1} |\beta_{z_1}\overline{\beta}_{z_2}|\tau(|z_1|^2)\tau(|z_2|^2)\tau_3(\Delta) .$$

From now on the condition $\tau(|z|^2) \leqslant \tau$ no longer exists, however we are left with the parameter $\tau$ in (9.10) and (10.4) which may be chosen at will.

LEMMA 10.1. *For any $\beta_z$ supported in the polar box (5.14) and satisfying $|\beta_z| \leqslant \tau(|z|^2)$ we have*

(10.7) $$Y(\beta) \ll \theta^4 N^2 (\log N)^{2^{19}} .$$

*Proof.* By Lemma 2.2 there exist $d, d_1, d_2 < N^{\frac{1}{4}}$ such that $d, d_1, d_2$ are mutually co-prime,

(10.8) $\quad d \mid \Delta(z_1, z_2), \quad d_1 \mid |z_1|^2, \quad d_2 \mid |z_2|^2$

(10.9) $\quad \tau_3(\Delta) \ll \tau(d)^{16} , \quad \tau(|z_1|^2) \ll \tau(d_1)^8 , \quad \tau(|z_2|^2) \ll \tau(d_2)^8 .$

Since the number of points $z_1, z_2 \in \mathfrak{B}$ satisfying (10.8) and (10.9) is bounded by $O\left(\theta^4 (N \log N)^2 \tau(dd_1 d_2)(dd_1 d_2)^{-1}\right)$ we obtain
$$Y(\beta) \ll \theta^4 (N \log N)^2 \sum_d \sum_{d_1} \sum_{d_2} \tau(dd_1 d_2)^{17}(dd_1 d_2)^{-1},$$
which yields (10.7). □

Inserting (10.7) into (10.4) we get

(10.10) $\quad \mathcal{D}_0(M,N) = 2\hat{f}(0)N^{-\frac{1}{2}}T(\beta) + O\left((\tau^{-1} + \theta)\theta^4 M^{\frac{1}{2}} N^{\frac{3}{2}} (\log N)^{2^{20}}\right)$ .

Here the error term satisfies the bound (5.25) provided

(10.11) $\quad \tau \geqslant \vartheta^{-2}(\log x)^{2^{20}} = (\log x)^{2A+2^{20}}$

and $\theta \leqslant \vartheta^2 (\log N)^{-2^{20}}$. The latter condition is assured since $A'$ can be chosen much larger than $A$, for example

$$(10.12) \qquad A' \geqslant 2A + 2^{20}$$

would suffice. The lower bound (10.11) is not in contradiction with the upper bound (9.12) provided that $B = B(\eta, A)$ in the statement of Proposition 4.1 is sufficiently large.

Now the last step is to estimate $T(\beta)$. Our target is

PROPOSITION 10.2.  *For the $\beta_z$ given by* (5.13) *and restricted by* (5.7), (5.8) *we have*

$$T(\beta) \ll N^2 (\log N)^{-\sigma} + P^{-1} N^2 \log N$$

*for any $\sigma > 0$, the implied constant depending on $\sigma$.*

*Remarks.* For the proof of Proposition 10.2 the lower bound for $P$ given in (4.4) is never utilized. Of course, since this assumption is now implicit in (5.8), the second term on the right side in this proposition is actually superfluous.

Inserting (8.16) in (10.5) and changing the order of summation we get

$$(10.13) \qquad T(\beta) = 2 \sum_d \frac{\varphi(d)}{d} \sum\sum_{\substack{(z_1, z_2) = 1 \\ \Delta(z_1, z_2) \equiv 0 (4d)}} \beta_{z_1} \overline{\beta}_{z_2} \left( \frac{z_2/z_1}{d} \right) \log 2 \left| \frac{z_1 z_2}{\Delta} \right| .$$

Note that $1 \leqslant d \leqslant N$ because $1 \leqslant |\Delta| \leqslant N$. We split this sum into

$$(10.14) \qquad T(\beta) = U(\beta) + V(\beta) + W(\beta)$$

where

$$(10.15)$$
$$U(\beta) = 2 \sum_{d \leqslant X} \frac{\varphi(d)}{d} \sum\sum_{\substack{(z_1, z_2) = 1 \\ \Delta(z_1, z_2) \equiv 0 (4d)}} f\left(\frac{|\Delta|}{d}\right) \beta_{z_1} \overline{\beta}_{z_2} \left( \frac{z_2/z_1}{d} \right) \log 2 \left| \frac{z_1 z_2}{\Delta} \right|$$

$$(10.16)$$
$$V(\beta) = 2 \sum_{d > X} \frac{\varphi(d)}{d} \sum\sum_{\substack{(z_1, z_2) = 1 \\ \Delta(z_1, z_2) \equiv 0 (4d)}} f\left(\frac{|\Delta|}{d}\right) \beta_{z_1} \overline{\beta}_{z_2} \left( \frac{z_2/z_1}{d} \right) \log 2 \left| \frac{z_1 z_2}{\Delta} \right|$$

$$(10.17)$$
$$W(\beta) = 2 \sum_d \frac{\varphi(d)}{d} \sum\sum_{\substack{(z_1, z_2) = 1 \\ \Delta(z_1, z_2) \equiv 0 (4d)}} f^*\left(\frac{|\Delta|}{d}\right) \beta_{z_1} \overline{\beta}_{z_2} \left( \frac{z_2/z_1}{d} \right) \log 2 \left| \frac{z_1 z_2}{\Delta} \right|$$

where $X$ will be chosen as a sufficiently large power of $\log N$ and $f, f^*$ are smooth functions whose graphs are

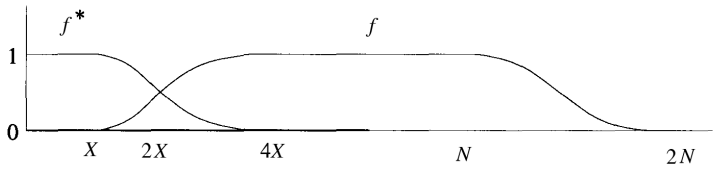

Note that $V(\beta)$ ranges over $d$ with $X < d < |\Delta|X^{-1}$ (this can be void) while $W(\beta)$ ranges over $d$ with $d^* = |\Delta|/4d < X$. Separating these cases with the smooth partition $f + f^*$ of unity will later help simplify some technical details. Each of these three parts will be estimated separately with considerable effort by very different methods.

First we shall deal with $V(\beta)$ because for the medium size moduli our arguments are quite general. No special properties of $\beta_z$ are needed. The source of cancellation is the sign change of the symbol $\chi_d(z_2/z_1)$, not the sign change of the Möbius function which is a component of $\beta_z$. Therefore we can afford to contaminate $\beta_z$ considerably when processing technical matters such as separation of variables. Thus, here we cut the range at $|\Delta|/4X$ which depends on the points $z_1, z_2$ to save technical work later in the range of large moduli. It will take us the next four sections to establish a general result (Proposition 14.1) which is adequate for application to $V(\beta)$ (see Proposition 15.1).

## 11. Jacobi-twisted sums over arithmetic progressions

Given a sequence $\mathcal{A} = (a_n)$ of complex numbers, it is natural to study its distribution in various residue classes $a \bmod d$. The goal is to establish a good approximation to

$$(11.1) \qquad A(x; d, a) = \sum_{\substack{n \leqslant x \\ n \equiv a (\bmod d)}} a_n$$

by a simple function of $d$ and with a relatively small error term uniformly for $d$ in a large range. For example, in sieve theory one considers the zero residue class in which case $A(x; d, 0)$ is well approximated by $g(d)A(x; 1, 0)$ with $g$ a nice multiplicative function. When $a$ is not the zero class $\bmod d$ then the expected main term for $A(x; d, a)$ is slightly different. Let us focus on the primitive classes $a$, that is with $(a, d) = 1$. Among these classes a reasonable sequence $\mathcal{A} = (a_n)$ is expected to be uniformly distributed, which means that

$$(11.2) \qquad r(x; d, a) = \sum_{\substack{n \leqslant x \\ n \equiv a (\bmod d)}} a_n - \frac{1}{\varphi(d)} \sum_{\substack{n \leqslant x \\ (n, d) = 1}} a_n$$

is quite small. Indeed, the theory of $L$-functions, if applicable, usually leads to the bound

$$(11.3) \qquad r(x;d,a) \ll \|\mathcal{A}\| x^{\frac{1}{2}} (\log x)^{-9A}$$

for any $A \geqslant 1$ with the implied constant depending only on $A$. In the case $a_n = \Lambda(n)$ this result is known as the Siegel-Walfisz theorem.

Although (11.3) holds uniformly for all primitive residue classes it is useful only for relatively small moduli; the bound (11.3) is trivial for $d > (\log x)^{9A}$. However, by the large sieve inequality one can extend (11.3) on average over such classes to moduli as large as $D = x(\log x)^{-A}$. Specifically, one shows that

$$(11.4) \qquad \sum_{d \leqslant D} \sideset{}{^*}\sum_{a \,(\mathrm{mod}\, d)} |r(x;d,a)|^2 \ll \|\mathcal{A}\|^2 x (\log x)^{-A}$$

with the same $A$ as in (11.3). In the case $a_n = \Lambda(n)$ this result, in a slightly less explicit form, is known as the Barban-Davenport-Halberstam theorem. The above generalization is, apart from the explicit power of $\log x$, given in [BFI]. Since we do not use (11.4) in this paper there is no need to provide a proof.

In this section we begin to investigate twisted sums in arithmetic progressions of type

$$(11.5) \qquad \sum\sum_{\bar{r}s \equiv a \,(\mathrm{mod}\, d)} \alpha_{rs} \left(\frac{r}{d'}\right)$$

and we shall continue the study of these sums in the following three sections. Here $\bar{r}$ denotes the multiplicative inverse modulo $d$. Note this implies $(r,d) = 1$. Also $d'$ denotes the odd part of $d$ and $\left(\frac{r}{d'}\right)$ is the Jacobi symbol. We pulled out $d'$ because it would be confusing to use here the convention (8.15). Keep in mind that we consider (11.5) with any complex numbers $\alpha_{rs}$ supported in the dyadic box

$$(11.6) \qquad R < r \leqslant 2R \quad \text{and} \quad S < s \leqslant 2S \, .$$

We define the local variance

$$(11.7) \qquad E(d) = \sum_{a \,(\mathrm{mod}\, d)} \Big| \sum\sum_{\bar{r}s \equiv a \,(\mathrm{mod}\, d)} \alpha_{rs} \left(\frac{r}{d'}\right) \Big|^2 \, .$$

Note that $E(d)$ is not restricted to primitive classes. We have

$$(11.8) \qquad E(d) = \sum\sum\sum\sum_{\bar{r}_1 s_1 \equiv \bar{r}_2 s_2 \,(\mathrm{mod}\, d)} \alpha_{r_1 s_1} \bar{\alpha}_{r_2 s_2} \left(\frac{r_1 r_2}{d'}\right) \, .$$

Applying additive characters we can also write

$$(11.9) \qquad E(d) = \frac{1}{d} \sum_{a \,(\mathrm{mod}\, d)} \Big| \sum_r \sum_s \alpha_{rs} \left(\frac{r}{d'}\right) e\left(\frac{a\bar{r}s}{d}\right) \Big|^2 \, .$$

For the variance $E^*(d)$ restricted to the primitive residue classes one can apply multiplicative characters getting

$$(11.10) \qquad E^*(d) = \frac{1}{\varphi(d)} \sum_{\chi (\bmod d)} \Big| \sum_r \sum_s \alpha_{rs} \Big(\frac{r}{d'}\Big) \chi(r) \bar{\chi}(s) \Big|^2 \ .$$

Our aim is to estimate the global variance

$$(11.11) \qquad V(D) = \sum_{D < d \leqslant 2D} \sum_{a (\bmod d)} \Big| \sum \sum_{\bar{r}s \equiv a (\bmod d)} \alpha_{rs} \Big(\frac{r}{d'}\Big) \Big|^2 \ .$$

If the variables $r, s$ are separable in the sense that the coefficients $\alpha_{rs}$ factor as

$$(11.12) \qquad \alpha_{rs} = \beta_r \gamma_s \ ,$$

or they are a linear combination of such things, then an application of the large sieve for character sums could be contemplated (cf. [BFI]). But nothing like (11.12) holds in our case! The absence of such a factorization and the presence of the Jacobi symbol necessitate new ideas to pursue the goal. Further comments on this point are made at the end of Section 16.

In this section we establish a basic estimate for $V(D)$ which will be extended in the next two sections and then summarized in Section 14. These four sections can be more or less considered as an independent unit. For this reason we shall feel free to use $V$ and later $W$ with meanings differing from those in Section 10, to which we shall later return.

We express our estimates in terms of the $\ell_2$-norm of the vector $\alpha = (\alpha_{rs})$;

$$(11.13) \qquad \|\alpha\|^2 = \sum_r \sum_s |\alpha_{rs}|^2 \ .$$

PROPOSITION 11.1. *Let $D, R, S \geqslant 1$. For any complex numbers $\alpha_{rs}$ supported in the box* (11.6) *we have*

$$(11.14) \qquad V(D) \ll \Big\{ D^{-\frac{1}{2}} RS + \Big(D\sqrt{RS} + RS^{\frac{3}{4}} + SR^{\frac{3}{4}}\Big)(RS)^\varepsilon \Big\} \|\alpha\|^2$$

*with any $\varepsilon > 0$ and the implied constant depending only on $\varepsilon$.*

*Remarks.* We have the trivial bound (use (11.8) and Lemma 11.2 below)

$$E(d) \ll \Big(d^{-1} RS + \sqrt{RS}\Big) \|\alpha\|^2 \ .$$

Hence
$$V(D) \ll (RS + D\sqrt{RS}) \|\alpha\|^2 \ .$$

Therefore the improvement in (11.14) appears in the first term. This first term $D^{-\frac{1}{2}}RS$ is not weakened by $(RS)^\varepsilon$, so it gives a nontrivial bound for all but very small moduli. One cannot do better because the moduli $d \sim D$ which are squares contribute to $V(D)$ at least $D^{-\frac{1}{2}}RS$. However if we restricted $d$ to squarefree numbers then our argument would give $D^{-1}RS$ in place of $D^{-\frac{1}{2}}RS$. The second term $D\sqrt{RS}$ is fine for any $D$ slightly smaller than $\sqrt{RS}$. The last two terms in (11.14) can probably be improved by refining our treatment but there is no reason to do so. These two terms contribute less than the trivial bound as long as $R$ and $S$ have the same order of magnitude in the logarithmic scale, that is if $R \gg S^\varepsilon$ and $S \gg R^\varepsilon$.

We precede the proof of Proposition 11.1 with three easy lemmas.

LEMMA 11.2. *For $(a,b,d) = 1$ the number $\mathcal{N}_d(R,S)$ of solutions to $ar \equiv bs \pmod{d}$ in positive integers $r \leqslant R$ and $s \leqslant S$ satisfies*

$$(11.15) \qquad \mathcal{N}_d(R,S) \leqslant d^{-1}RS + \sqrt{RS} \;.$$

*Proof.* Dividing the congruence by $\alpha = (a,d)$ and $\beta = (b,d)$, then counting the solutions in two ways we infer the bound

$$\min\left\{\frac{R}{\beta}\left(\frac{S\beta}{d}+1\right), \frac{S}{\alpha}\left(\frac{R\alpha}{d}+1\right)\right\} \leqslant d^{-1}RS + \min(R,S)$$

which yields (11.15). □

LEMMA 11.3. *For $e \geqslant 1$ and $Q, R \geqslant 2$ we have*

$$(11.16) \qquad \sum_{\substack{q \leqslant Q \\ q' \neq \square}} \Big| \sum_{\substack{r \leqslant R \\ (r,e)=1}} \left(\frac{r}{q'}\right) \Big|^2 \ll \min\left\{\tau(e)^2 Q^2, QR + R^3\right\} (\log QR)^2 \;.$$

Here, and hereafter, $q' \neq \square$ means that $q'$ is not the square of an integer.

*Proof.* By the Polyá-Vinogradov bound the inner sum is $O\left(\tau(e)\sqrt{q}\log q\right)$ giving the first estimate of (11.16). To get the second estimate we ignore the condition $q' \neq \square$, then we square out and change the order of summation obtaining

$$\sum_{r_1}\sum_{r_2} \Big| \sum_q \left(\frac{r_1 r_2}{q'}\right) \Big| \;.$$

The terms with $r_1 r_2 = \square$ contribute $O(QR \log R)$ by trivial estimation. The remaining terms contribute $O\left(R^3(\log Q)(\log R)\right)$ by the Polyá-Vinogradov estimate. Combining these contributions we derive (11.16). □

LEMMA 11.4. *We have*

$$(11.17) \qquad \sum_{\substack{d_1,d_2 \leqslant D \\ d'_1 d'_2 \neq \square}} \frac{1}{[d_1,d_2]} \Big| \sum_{\substack{r \leqslant R \\ (r,d_1 d_2)=1}} \left(\frac{r}{d'_1 d'_2}\right) \Big|^2 \ll R^{\frac{3}{2}}(DR)^\varepsilon$$

*for any $\varepsilon > 0$, the implied constant depending on $\varepsilon$.*

*Proof.* The condition $(r, d_1 d_2) = 1$ is redundant unless $d_1 d_2$ is even in which case it simply means that $r$ is odd. Thus we can ignore this condition by inclusion-exclusion. Hence the sum is bounded by

$$\sum_{e \leqslant D} \sum_{\substack{b_1, b_2 \leqslant D \\ b'_1 b'_2 \neq \square}} (eb_1 b_2)^{-1} | \sum_{\substack{r \leqslant R \\ (r,e)=1}} \left( \frac{r}{b'_1 b'_2} \right) |^2$$

$$\ll D^\varepsilon \sum_{e \leqslant D} \frac{1}{eQ} \sum_{\substack{q \leqslant Q \\ q' \neq \square}} | \sum_{\substack{r \leqslant R \\ (r,e)=1}} \left( \frac{r}{q'} \right) |^2$$

for some $Q \leqslant D^2$. By Lemma 11.3 this is $\ll \min\{Q, R + Q^{-1} R^3\}(DR)^\varepsilon$. This bound yields (11.17), the worst $Q$ being $R^{\frac{3}{2}}$. □

Our treatment of $V(D)$ goes via the dual sum

(11.18) $$W(D) = \sum_r \sum_s | \sum_d \sum_{a \equiv \bar{r} s \pmod{d}} \gamma_{ad} \left( \frac{r}{d'} \right) |^2 \ .$$

Here $r, s, d, a$ run over the same ranges as in $V(D)$ and $\gamma_{ad}$ are any complex numbers. By the duality principle familiar from the theory of the large sieve (see, for example, page 32 of [Bo2]) the estimate (11.14) is equivalent to

(11.19) $$W(D) \ll \left\{ D^{-\frac{1}{2}} RS + \left( D\sqrt{RS} + RS^{\frac{3}{4}} + SR^{\frac{3}{4}} \right)(RS)^\varepsilon \right\} \|\gamma\|^2 \ .$$

Now we are going to prove (11.19). First we enlarge $W(D)$ by attaching a smooth majorant $f(s)$ such that

$$f(s) \geqslant 0 \quad , \quad f(s) \geqslant 1 \text{ if } S < s \leqslant 2S,$$
$$\hat{f}(0) = 2S \quad , \quad \hat{f}(t) \ll S(1 + |t|S)^{-2} \ .$$

Then we square out and change the order of summation getting

$$W \leqslant \sum_{d_1} \sum_{d_2} \sum_{a \pmod{q}} \gamma_{ad_1} \bar{\gamma}_{ad_2} \sum \sum_{\bar{r} s \equiv a \pmod{q}} f(s) \left( \frac{r}{d'_1 d'_2} \right)$$

where $q = [d_1, d_2]$ is the least common multiple of $d_1$ and $d_2$. Here and from now on, often without writing it, we assume that $R < r \leqslant 2R$.

The terms with $d'_1 d'_2 = \square$ contribute

$$W^\square = \sum \sum_{d'_1 d'_2 = \square} \sum \sum_{\substack{a_1 \pmod{d_1} \\ a_2 \pmod{d_2}}} \gamma_{a_1 d_1} \bar{\gamma}_{a_2 d_2} \sum_{(r, d'_1 d'_2)=1} \sum_{\substack{\bar{r} s \equiv a_1 \pmod{d_1} \\ \bar{r} s \equiv a_2 \pmod{d_2}}} f(s)$$

$$\ll \sum \sum_{d'_1 d'_2 = \square} \sum_{a_1 \pmod{d_1}} |\gamma_{a_1 d_1}|^2 \sum \sum_{\bar{r} s \equiv a_1 \pmod{d_1}} f(s) \ .$$

Hence by Lemma 11.2,

$$W^\square \ll \left(D^{-1}RS + \sqrt{RS}\right) \sum_{d_1} \nu(d_1) \sum_{a_1 \pmod{d_1}} |\gamma_{a_1 d_1}|^2$$

where $\nu(d_1)$ is the number of $d_2 \sim D$ such that $d_1' d_2' = \square$. We have $\nu(d_1) \ll \sqrt{D}$ (or even better $\nu(d_1) \ll 1$ if $d_1, d_2$ were squarefree) giving

(11.20) $$W^\square \ll \left(D^{-\frac{1}{2}}RS + \sqrt{DRS}\right) \|\gamma\|^2$$

where

$$\|\gamma\|^2 = \sum_d \sum_{a \pmod d} |\gamma_{ad}|^2 .$$

Next we estimate the contribution $W^\diamond$ of terms with $d_1' d_2' \neq \square$. By Cauchy's inequality

(11.21) $$|W^\diamond|^2 \leqslant \sigma(\gamma) \sum\sum_{d_1' d_2' \neq \square} W(d_1, d_2)$$

where

$$\sigma(\gamma) = \sum_{d_1} \sum_{d_2} \sum_{a \pmod q} |\gamma_{ad_1} \gamma_{ad_2}|^2 \leqslant \|\gamma\|^4 ,$$

and $W(d_1, d_2)$ is a local variance to modulus $q = [d_1, d_2]$, namely

$$W(d_1, d_2) = \sum_{a \pmod q} \left| \sum\sum_{\bar{r}s \equiv a \pmod q} f(s) \left(\frac{r}{d_1' d_2'}\right) \right|^2 .$$

By Poisson summation for the sum over $s$ we have the Fourier expansion

$$\sum_{s \equiv ar \pmod q} f(s) = \frac{1}{q} \sum_h \hat{f}\left(\frac{h}{q}\right) e\left(\frac{ahr}{q}\right) .$$

Hence, by grouping terms according to the product $hr$ we get

$$W(d_1, d_2) = \frac{1}{q^2} \sum_{a \pmod q} \left| \sum_k c_k\, e\left(\frac{ak}{q}\right) \right|^2$$

where

$$c_k = \sum\sum_{\substack{hr=k \\ (r,q)=1}} \hat{f}\left(\frac{h}{q}\right) \left(\frac{r}{d_1' d_2'}\right) .$$

Next, by the popular inequality $|x+y|^2 \leqslant 2|x|^2 + 2|y|^2$ and the orthogonality of additive characters

$$W(d_1, d_2) \leqslant \frac{2}{q}|c_0|^2 + \frac{2}{q} \sum\sum_{\substack{k_1 \equiv k_2 \pmod q \\ k_1 k_2 \neq 0}} c_{k_1} \bar{c}_{k_2} .$$

The zero coefficient is the real character sum

$$c_0 = \hat{f}(0) \sum_{(r,q)=1} \left(\frac{r}{d'_1 d'_2}\right)$$

with $\hat{f}(0) = 2S$ and the other coefficients satisfy $c_k \ll \tau(k) S (1 + |k|S/qR)^{-2}$, by a trivial estimation. Hence we infer that

$$W(d_1, d_2) \ll q^{-1} S^2 \, \Big| \sum_{(r,q)=1} \left(\frac{r}{d'_1 d'_2}\right) \Big|^2 + (R+S) R (qR)^\varepsilon \, .$$

Summing over $d_1$, $d_2$ we obtain by Lemma 11.4 that

(11.22) $$W^\diamond \ll \left( SR^{\frac{3}{4}} + DR + D\sqrt{RS} \right) (DR)^\varepsilon \|\gamma\|^2 \, .$$

Adding (11.20) and (11.22) we get

(11.23) $$W(D) \ll \left\{ D^{-\frac{1}{2}} RS + \left( SR^{\frac{3}{4}} + RS^{\frac{3}{4}} + D\sqrt{RS} + DR \right) (RS)^\varepsilon \right\} \|\gamma\|^2 \, .$$

Here we have added the extra term $RS^{3/4}$ to gain some symmetry, and we replaced $(DR)^\varepsilon$ by $(RS)^\varepsilon$ because if $D > RS$ the estimate (11.23) is trivial.

To complete the proof of (11.19) it remains to remove the term $DR$ in (11.23). First we look at the sum $W^*(D)$ reduced by the condition $(a, d') = 1$, in other words $W^*(D)$ is the sum $W(D)$ for vectors $\gamma_{ad}$ with $(a, d') = 1$. Clearly, if we switch $r$ with $s$ and also change $\gamma_{ad}$ to $\gamma_{ad}\left(\frac{a}{d'}\right)$ then $W^*(D)$ is not altered. Therefore due to this symmetry we may assume that $R \leqslant S$. Applying (11.23) we get

(11.24) $$W^*(D) \ll \left\{ D^{-\frac{1}{2}} RS + \left( SR^{\frac{3}{4}} + RS^{\frac{3}{4}} + D\sqrt{RS} \right) (RS)^\varepsilon \right\} \|\gamma\|^2 \, .$$

Now we deduce the same bound for $W(D)$. To this end, we transform $W(D)$ as follows

$$W(D) = \sum_r \sum_s \Big| \sum_{e|s'} \sum_{\substack{(d',a)=1 \\ a \equiv \bar{r}s/e \pmod d}} \gamma_{ea,ed} \left(\frac{r}{ed'}\right) \Big|^2$$

$$\leqslant \sum_r \sum_s \sum_{e|s'} \sum_{f|s'} e^\alpha f^{-\alpha} \Big| \sum_{\substack{(d',a)=1 \\ a \equiv \bar{r}s/e \pmod d}} \gamma_{ea,ed} \left(\frac{r}{ed'}\right) \Big|^2$$

$$\leqslant \sum_e \sum_f \tau(e) e^\alpha f^{-\alpha} \sum_{r \sim R} \sum_{s \sim S/ef} \Big| \sum_{\substack{(d',a)=1 \\ a \equiv \bar{r}sf \pmod d}} \gamma_{ea,ed} \left(\frac{r}{d'}\right) \Big|^2 \, .$$

Here we have written $1 = e^\alpha e^{-\alpha}$ with $0 < \alpha < \frac{1}{2}$ which yields a factor $e^\alpha f^{-\alpha}$ by applying Cauchy's inequality. This factor is needed for technical reasons, namely to make one term in the following bound free of the divisor function while losing only slightly in the remaining terms. Of course, such a refinement is not crucial but it will simplify some forthcoming work. Now applying (11.24) we get

$$W(D) \ll \sum_e \sum_f \tau(e) e^\alpha f^{-\alpha} (ef)^{-\frac{1}{2}} \sum\sum_{(d,a)=1} |\gamma_{ea,ed}|^2$$
$$\left\{ (Df)^{-\frac{1}{2}} RS + \left( D\sqrt{RS} + RS^{\frac{3}{4}} + SR^{\frac{3}{4}} \right) (RS/ef)^\varepsilon \right\}.$$

We choose $\alpha = \frac{1-\varepsilon}{2}$ so the series for $f$ converges and the function $\tau(e) e^{\alpha - \frac{1}{2}}$ is bounded. By the above estimates, and since

$$\sum_e \sum\sum_{(d,a)=1} |\gamma_{ea,ed}|^2 = \|\gamma\|^2 ,$$

we obtain (11.19) for $W(D)$ as claimed. Finally, by the duality principle, (11.19) implies (11.14) for $V(D)$. □

## 12. Flipping moduli

The bound (11.14) is nontrivial in the range

$$(\log 2RS)^A < D < (RS)^{\frac{1}{2}-\varepsilon} .$$

In this section we leapfrog by reflection over the middle to cover the range

$$(RS)^{\frac{1}{2}+\varepsilon} < D < RS(\log 2RS)^{-A} .$$

For technical convenience we assume that our vectors $\alpha = (\alpha_{rs})$ satisfy

(12.1) $$(r, 2s) = 1 .$$

Furthermore by splitting into four residue classes we can assume without loss of generality that $r$ is fixed modulo eight. We have by the reciprocity law

$$\left( \frac{r}{d'} \right) = \pm \left( \frac{d}{r} \right)$$

where $\pm$ depends on $d$ and on $r \pmod 8$ but not on $r$ in any other fashion. Therefore $V(D)$ for our vectors can be written as

(12.2) $$V(D) = \sum_{D<d\leqslant 2D} \sum_{a (\mod d)} \Big| \sum\sum_{\bar r s \equiv a (\mod d)} \alpha_{rs} \left( \frac{d}{r} \right) \Big|^2 .$$

PROPOSITION 12.1.  *Let $D, R, S \geqslant 1$. For any complex numbers $\alpha_{rs}$ with $(r, 2s) = 1$ supported in the box (11.6) we have*

(12.3) $$V(D) \leqslant \mathcal{L}(D, R, S) \sum_r \sum_s \tau(r) |\alpha_{rs}|^2$$

where $\mathcal{L}(D, R, S)$ satisfies the bound

$$\text{(12.4)} \quad \mathcal{L}(D, R, S) \ll D + D^{\frac{1}{3}}(RS)^{\frac{2}{3}}(\log 2RS)^4 \\ + \left[D^{-1}(RS)^{\frac{3}{2}} + RS^{\frac{3}{4}} + SR^{\frac{3}{4}}\right](RS)^{\varepsilon},$$

with any $\varepsilon > 0$ and the implied constant depending only on $\varepsilon$.

Proposition 12.1 will be derived from Proposition 11.1 by way of interpreting the congruence $r_1 s_2 \equiv r_2 s_1 \pmod{d}$ as $r_1 s_2 \equiv r_2 s_1 \pmod{k}$ where $k$ is the complementary divisor. The idea is reminiscent to that used by C. Hooley [Ho] in a similar context (the Barban-Davenport-Halberstam theorem), however the presence of Jacobi symbol makes our case quite subtle since the reciprocity law is essentially employed in the process of flipping moduli. There will be some complications related to various co-primality conditions. In order to distribute the burden of complications evenly throughout the proof we are going to prepare the following variation of Proposition 11.1.

PROPOSITION 11.1*. *Let $m, D, R, S \geqslant 1$. For any complex numbers $\alpha_{rs}$ and $\beta_{rs}$ supported in the box (11.6) we have*

$$\text{(12.5)} \quad \sum_{D < d \leqslant 2D} \left| \sum\sum\sum\sum_{\substack{r_1 s_2 \equiv r_2 s_1 \pmod{dm} \\ (r_1, r_2) = 1}} \alpha_{r_1 s_1} \overline{\beta}_{r_2 s_2} \left(\frac{d}{r_1 r_2}\right) \right| \leqslant \mathcal{M}(mD, R, S) \|\alpha\| \|\beta\|$$

where $\mathcal{M}(D, R, S)$ satisfies the bound

$$\mathcal{M}(D, R, S) \ll D^{-\frac{1}{2}} RS \log 2R + \left(D\sqrt{RS} + RS^{\frac{3}{4}} + SR^{\frac{3}{4}}\right)(RS)^{\varepsilon}$$

with any $\varepsilon > 0$ and the implied constant depending only on $\varepsilon$.

*Proof.* As compared to the local variance (11.8) the bilinear form in (12.5) has a few extra features. To handle the extra modulus $m$ we write $\ell = dm$ and note that $(r_1 r_2, m) = 1$, so that we can separate $\left(\frac{d}{r_1 r_2}\right) = \left(\frac{m}{r_1 r_2}\right)\left(\frac{\ell}{r_1 r_2}\right)$, and then transfer $\left(\frac{m}{r_1 r_2}\right)$ to the coefficients $\alpha_{r_1 s_1} \beta_{r_2 s_2}$. The new modulus $\ell$ runs over the segment $mD < \ell \leqslant 2mD$ and satisfies $\ell \equiv 0 \pmod{m}$ but we shall ignore the latter property by positivity. Then we remove the condition $(r_1, r_2) = 1$ by Möbius inversion. After these transformations we see that the sum (12.5) is bounded by

$$\sum_\rho \sum_{\ell \sim mD} \sum_{a \pmod{\ell}} \left| \sum\sum_{\overline{r}s \equiv a \pmod{\ell}} \alpha_{\rho rs}\left(\frac{m\ell}{r}\right) \right| \left| \sum\sum_{\overline{r}s \equiv a \pmod{\ell}} \beta_{\rho rs}\left(\frac{m\ell}{r}\right) \right|.$$

If we denote, for a given $\rho$, the above inner multiple sum by $L(\alpha, \beta)$ then, on applying Cauchy's inequality, $L^2(\alpha, \beta) \leqslant L(\alpha, \alpha) L(\beta, \beta)$. Now, Proposition 11.1 is applicable to $L(\alpha, \alpha)$ giving

$$L(\alpha,\alpha) \ll \left\{ \frac{RS}{\rho\sqrt{mD}} + \left[mD\sqrt{RS} + RS^{\frac{3}{4}} + SR^{\frac{3}{4}}\right](RS)^{\varepsilon} \right\} \sum_r \sum_s |\alpha_{\rho rs}|^2.$$

Similarly $L(\beta, \beta)$ satisfies the same estimate with $\alpha_{\rho rs}$ replaced by $\beta_{\rho rs}$. From these two estimates one arrives at Proposition 11.1* by summing over $\rho$.   □

Now we proceed to the proof of Proposition 12.1. We begin by attaching to $V(D)$ a smooth majorant $f(y)$ supported on $\frac{1}{2}D \leqslant y \leqslant 3D$, then we square out getting

$$(12.6) \qquad V(D) \leqslant \sum_d f(d) \sum\sum\sum\sum_{r_1 s_2 \equiv r_2 s_1 (\mathrm{mod}\, d)} \alpha_{r_1 s_1} \overline{\alpha}_{r_2 s_2} \left(\frac{d}{r_1 r_2}\right).$$

Next we remove the terms near the diagonal. To do this smoothly we use the function $g(x)$ whose graph is with $H$ to be chosen later subject to $1 \leqslant H \leqslant RSD^{-1}$. Notice that for $r_1, s_1$ and $r_2, s_2$ in the box (11.6) we have

$$\left| \frac{s_1}{r_1} - \frac{s_2}{r_2} \right| < \frac{2S}{R},$$

so the factor $g\left(\left|\frac{s_1}{r_1} - \frac{s_2}{r_2}\right|\right)$ may be inserted into (12.6) without alteration, except for the points $r_1, s_1$ and $r_2, s_2$ with

$$\left| \frac{s_1}{r_1} - \frac{s_2}{r_2} \right| < \frac{2S}{HR}.$$

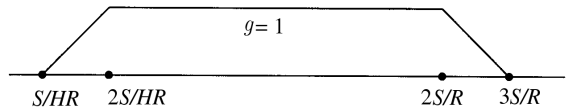

The contribution of these exceptional points is estimated trivially by

$$V_0(f, g) \ll D\|\alpha\|^2 + \sum_r \sum_s \nu(r, s)|\alpha_{rs}|^2$$

where the first term comes from the points exactly on the diagonal $r_1 s_2 = r_2 s_1$ (note that this equation implies $r_1 = r_2$ and $s_1 = s_2$ by virtue of (12.1)), and the second from the rest. Here $\nu(r, s)$ is the number of integers $k$ with $1 \leqslant k < 8RS(DH)^{-1}$ and rationals $s_1/r_1$ such that

$$\left| \frac{s_1}{r_1} - \frac{s}{r} \right| < \frac{2S}{HR} \quad \text{and} \quad r_1 s \equiv rs_1 \pmod{k}.$$

Each $k$ is obtained by flipping $d$ to the complementary divisor of $|r_1 s - r s_1|$. Given $k$ one shows, as in the proof of Lemma 11.2, that the number of rationals $s_1/r_1$ satisfying the above conditions is $\nu_k(r,s) \ll RS(Hk)^{-1} + \sqrt{RS}$. Summing over $k$ we obtain

$$\nu(r,s) \ll \frac{RS}{H}\left(\log \frac{2RS}{DH}\right) + \frac{RS}{DH}\sqrt{RS} \ .$$

This yields

(12.7) $\qquad V_0(f,g) \ll \left(D + H^{-1} RS \log 2RS + H^{-1} D^{-1} R^{\frac{3}{2}} S^{\frac{3}{2}}\right) \|\alpha\|^2 \ .$

The remaining points contribute

$$V(f,g) = \sum_d f(d) \sum\sum\sum\sum_{r_1 s_2 \equiv r_2 s_1 (\bmod d)} \alpha_{r_1 s_1} \overline{\alpha}_{r_2 s_2}\left(\frac{d}{r_1 r_2}\right) g\left(\Big|\frac{s_1}{r_1} - \frac{s_2}{r_2}\Big|\right) \ .$$

We reduce the variables $r_1, r_2$ by the common divisor $c = (r_1, r_2)$ and remove the resulting condition $(c,d) = 1$ by Möbius inversion getting

(12.8) $\qquad V(f,g) = \sum_c \sum_{m|c} \mu(m) V_{cm}(f,g)$

where $V_{cm}(f,g)$ is the sum

(12.9) $\qquad \sum_d f(dm) \sum\sum\sum\sum_{\substack{r_1 s_2 \equiv r_2 s_1 (dm) \\ (r_1, r_2)=1}} \alpha_{cr_1 s_1} \overline{\alpha}_{cr_2 s_2}\left(\frac{dm}{r_1 r_2}\right) g\left(\frac{1}{c}\Big|\frac{s_1}{r_1} - \frac{s_2}{r_2}\Big|\right) .$

For every $c$ and $m|c$ we estimate $V_{cm}(f,g)$ separately. Here we flip $d$ to the complementary divisor of $|r_1 s_2 - r_2 s_1| m^{-1}$. We write $|r_1 s_2 - r_2 s_1| = dmq$ and interpret each and every term involving $d$ in terms of $q$. We have

$$f(dm) = f\left(\Big|\frac{s_1}{r_1} - \frac{s_2}{r_2}\Big|\frac{r_1 r_2}{q}\right)$$

and using the reciprocity law we get (recall that $r_1, r_2$ are co-prime, odd, and congruent modulo eight)

$$\left(\frac{dm}{r_1 r_2}\right) = \left(\frac{s_1}{r_1}\right)\left(\frac{s_2}{r_2}\right)\left(\frac{q}{r_1 r_2}\right) \ .$$

This transforms $V_{cm}(f,g)$ into the sum

(12.10) $\qquad \sum_q \sum\sum\sum\sum_{\substack{r_1 s_2 \equiv r_2 s_1 (mq) \\ (r_1, r_2)=1}} \beta_{r_1 s_1} \overline{\beta}_{r_2 s_2}\left(\frac{q}{r_1 r_2}\right) B\left(\frac{1}{c}\left(\frac{s_1}{r_1} - \frac{s_2}{r_2}\right), \frac{cr_1 r_2}{q}\right)$

where

(12.11) $\qquad B(x,y) = f(|x|y) g(|x|)$

and
$$\beta_{rs} = \left(\frac{s}{r}\right)\alpha_{crs} .$$

Next we separate the variables $q, r_1, s_1, r_2, s_2$ involved in $B(x,y)$ for the $x$, $y$ relevant to (12.10). Before doing this we observe that $q$ runs over the segment

(12.12) $$\frac{RS}{cDH} < q < \frac{24RS}{cD}$$

which fact follows by examining the support of $g$. Having recorded this information we represent $B(x,y)$ as the Fourier-Mellin transform in $x$ and $y$ respectively,

(12.13) $$f(|x|y)g(|x|) = \iint h(u,t)e(ux)y^{it} du\, dt .$$

Here we have by the inverse transform
$$h(u,t) = \frac{1}{2\pi}\iint f(|x|y)g(|x|)e(-ux)y^{-1-it} dx\, dy$$
$$= \frac{1}{\pi}\int_0^\infty x^{it}g(x)\cos(2\pi ux)dx \int_0^\infty f(y)y^{-1-it}dy .$$

Integrating the Mellin transform by parts four times we get
$$\int_0^\infty f(y)y^{-1-it}dy \ll (t^2+1)^{-2} ,$$

and integrating the Fourier transform by parts up to two times we get
$$\int_0^\infty x^{it}g(x)\cos(2\pi ux)dx \ll (t^2+1)\min\left\{\frac{S}{R}, \frac{1}{|u|}, \frac{HR}{u^2 S}\right\}\log 2H .$$

From these estimates it follows that the $L_1$-norm of $h(u,t)$ satisfies

(12.14) $$\iint |h(u,t)| du\, dt \ll (\log 2H)^2 .$$

Applying (12.13) to (12.10) we infer that

(12.15) $$V_{cm}(f,g) \ll (\log 2H)^2 \sum_q \left|\sum\sum\sum\sum_{\substack{r_1 s_2 \equiv r_2 s_1 (mq) \\ (r_1,r_2)=1}} \gamma_{r_1 s_1}\tilde{\gamma}_{r_2 s_2}\left(\frac{q}{r_1 r_2}\right)\right|$$

where

(12.16) $$\gamma_{rs} = r^{it}e\left(\frac{us}{cr}\right)\left(\frac{s}{r}\right)\alpha_{crs}$$

for some real $u$ and $t$ (the $\sim$ denotes complex conjugation except for $r^{it}$ which remains unchanged).

At last we are ready to apply Proposition 11.1*. This gives, by (12.12), (12.15) and (12.16),

$$V_{cm}(f,g) \ll \left\{ (cm)^{-\frac{1}{2}}(DHRS)^{\frac{1}{2}}(\log 2RS)^3 \right.$$
$$\left. + \left[ c^{-\frac{3}{2}}mD^{-1}(RS)^{\frac{3}{2}} + RS^{\frac{3}{4}} + SR^{\frac{3}{4}} \right](RS)^\varepsilon \right\} \sum_r \sum_s |\alpha_{crs}|^2 \ .$$

Summing over $m$ and $c$ as in (12.8) we obtain

(12.17) $$V(f,g) \leqslant \mathcal{H}(D,H,R,S) \sum_r \sum_s \tau(r)|\alpha_{rs}|^2$$

where $\mathcal{H}(D,H,R,S)$ satisfies the bound

$$\mathcal{H}(D,H,R,S) \ll (DHRS)^{\frac{1}{2}}(\log 2RS)^4 + \left[ D^{-1}(RS)^{\frac{3}{2}} + RS^{\frac{3}{4}} + SR^{\frac{3}{4}} \right](RS)^\varepsilon.$$

Adding (12.17) to (12.7) and choosing $H = (D^{-1}RS)^{1/3}$ we complete the proof of (12.3). $\square$

## 13. Enlarging moduli

After the results of the last two sections we need a nontrivial bound for $V(D)$ in the middle range

$$(RS)^{\frac{1}{2}-\varepsilon} < D < (RS)^{\frac{1}{2}+\varepsilon} \ .$$

We deal with these moduli by establishing an inequality between $V(D)$ and $V(D^+)$ for $D^+ > D$, thereby allowing us to appeal to the result for larger moduli given in the previous section. We accomplish this by considering, for each given $d$ the special moduli $dp^2$ where $p$ runs through a set of primes. Since the Jacobi symbol for the enlarged modulus $dp^2$ is essentially the same as that for $d$ this gives us a multiple counting of the original sum. Without this multiplicity we would gain nothing and for this reason the method fails for general characters. Different arguments but of similar spirit appear already in the paper [BD] of Bombieri and Davenport.

Let $D, R, S \geqslant 1$ and $\alpha_{rs}$ be any complex numbers for $(r, 2s) = 1$ supported in the box (11.6). For any prime $p$ we write

$$\sum\sum_{\bar{r}s \equiv a(\bmod d)} \alpha_{rs}\left(\frac{d}{r}\right) = \sum\sum_{\substack{\bar{r}s \equiv a(\bmod d) \\ p \nmid r}} \alpha_{rs}\left(\frac{d}{r}\right) + \sum\sum_{\substack{\bar{r}s \equiv a(\bmod d) \\ p \mid r}} \alpha_{rs}\left(\frac{d}{r}\right) \ .$$

Hence we infer

(13.1) $$V(D) \leqslant 2 \sum_{D < d \leqslant 2D} \left( E_p(d) + E'_p(d) \right)$$

where
$$E_p(d) = \sum_{a(\bmod d)} | \sum\sum_{\substack{\bar{r}s \equiv a(\bmod d) \\ p \nmid r}} \alpha_{rs} \left(\frac{d}{r}\right) |^2$$

and $E'_p(d)$ is given by the same formula but with the condition $p|r$ in place of $p \nmid r$. We estimate the latter trivially as follows

$$\sum_d E'_p(d) \leqslant \sum_d \sum\sum\sum\sum_{r_1 s_2 \equiv r_2 s_1 (\bmod d)} |\alpha_{pr_1 s_1} \alpha_{pr_2 s_2}|$$
$$\ll \{D + p^{-1}(RS)^{1+\varepsilon}\} \sum_r \sum_s |\alpha_{prs}|^2 .$$

Multiplying (13.1) through by $p^{-1} \log p$ and summing over $P < p \leqslant 2P$ we obtain

(13.2) $$V(D) \ll \sum_{p \sim P} \frac{\log p}{p} \sum_{d \sim D} E_p(d) + (DP^{-1} + P^{-2}RS)(RS)^\varepsilon \|\alpha\|^2 .$$

In general we have the following inequality (monotonicity of the local variance)

$$\begin{aligned} A(d) &= \sum_{a(\bmod d)} | \sum_{n \equiv a(\bmod d)} \alpha_n |^2 \\ &= \frac{1}{d} \sum_{b(\bmod d)} | \sum_n \alpha_n e\left(\frac{bn}{d}\right) |^2 \\ &\leqslant \frac{1}{d} \sum_{b(\bmod dm)} | \sum_n \alpha_n e\left(\frac{bn}{dm}\right) |^2 = mA(dm) . \end{aligned}$$

We apply this with $m = p^2$ to $E_p(d)$ getting

(13.3) $$E_p(d) \leqslant p^2 E(dp^2)$$

since $\left(\frac{dp^2}{r}\right) = \left(\frac{d}{r}\right)$ for $p \nmid r$. Inserting (13.3) into (13.2) we obtain the following

PRINCIPLE OF ENLARGING MODULI. *For every $D, P, R, S \geqslant 1$ we have*

(13.4) $$V(D) \ll V(D^+) P \log DP + (DP^{-1} + P^{-2}RS)(RS)^\varepsilon \|\alpha\|^2$$

*for some $D^+$ with $DP^2 \leqslant D^+ \leqslant 4DP^2$ and where the implied constant depends only on $\varepsilon$.*

Actually $D^+$ may be taken to be $2^j DP^2$ for one of $j = 0, 1, 2$. Note that in the recurrent term on the right side we have enlarged the range of moduli by a factor $P^2$ while losing only $P \log DP$ in the bound. Combining Proposition 12.1 with the above principle we infer

PROPOSITION 13.1. *Let $D, R, S \geqslant 1$. For any complex numbers $\alpha_{rs}$ with $(r, 2s) = 1$ supported in the box (11.6) we have*

$$(13.5) \quad V(D) \ll \left\{ D(R+S)^{\frac{1}{4}}(RS)^{\frac{1}{4}} + D^{-\frac{1}{2}}(R+S)^{\frac{1}{8}}(RS)^{\frac{9}{8}} \right\} (RS)^{\varepsilon} \|\alpha\|^2 \ .$$

*Proof.* We can assume that

$$(13.6) \quad (R+S)^{\frac{1}{4}}(RS)^{\frac{1}{4}} \leqslant D \leqslant (R+S)^{-\frac{1}{4}}(RS)^{\frac{3}{4}}$$

or else the bound (13.5) is trivial. Applying (12.3) to $V(D^+)$ in (13.4) we obtain

$$V(D) \ll \left\{ D^{\frac{1}{3}} P^{\frac{5}{3}}(RS)^{\frac{2}{3}} + (DP)^{-1}(RS)^{\frac{3}{2}} \right.$$
$$\left. + P(R+S)^{\frac{1}{4}}(RS)^{\frac{3}{4}} + DP^3 + P^{-2}RS \right\} (PRS)^{\varepsilon} \|\alpha\|^2 \ .$$

We choose $P = D^{-\frac{1}{2}}(R+S)^{-\frac{1}{8}}(RS)^{\frac{3}{8}}$ getting (13.5). □

## 14. Jacobi-twisted sums: Conclusion

We combine the results of Sections 11, 12, 13 to formulate a bound for $V(D)$ which is nontrivial throughout the range

$$(\log 2RS)^A < D < RS(\log 2RS)^{-A}$$

just as in the Barban-Davenport-Halberstam Theorem.

PROPOSITION 14.1. *Let $D, R, S \geqslant 1$. For any complex numbers $\alpha_{rs}$ with $(r, 2s) = 1$ supported in the box $R < r < 2R$, $S < s < 2S$, we have*

$$\sum_{D < d \leqslant 2D} \sum_{a(\bmod d)} \Big| \sum\sum_{\bar{r}s \equiv a(\bmod d)} \alpha_{rs} \left(\frac{r}{d}\right) \Big|^2 \leqslant \mathcal{N}(D, R, S) \sum_r \sum_s \tau(r) |\alpha_{rs}|^2$$

*where $\mathcal{N}(D, R, S)$ satisfies the bound*

$$\mathcal{N}(D, R, S) \ll D + D^{-\frac{1}{2}} RS + D^{\frac{1}{3}} (RS)^{\frac{2}{3}} (\log 2RS)^4 + (R+S)^{\frac{1}{12}} (RS)^{\frac{11}{12} + \varepsilon}$$

*for any $\varepsilon > 0$, the implied constant depending only on $\varepsilon$. Here $\left(\frac{r}{d}\right)$ is the Jacobi symbol if $d$ is odd and is extended for $d$ even by (8.15).*

*Proof.* This follows by application of Proposition 11.1 if $D \leqslant D_1$, Proposition 13.1 if $D_1 < D < D_2$, and Proposition 12.1 if $D \geqslant D_2$, where

$$D_1 = (R+S)^{\frac{1}{12}}(RS)^{\frac{5}{12}} \quad \text{and} \quad D_2 = (R+S)^{-\frac{1}{8}}(RS)^{\frac{5}{8}} \ .$$

## 15. Estimation of $V(\beta)$

Now we are ready to estimate $V(\beta)$ which is the partial sum of (10.13) over the moduli $d$ in the middle range as defined in (10.16). We write $z_1 = r_1 + is_1$ and $z_2 = r_2 + is_2$ so $\Delta(z_1, z_2) = r_1 s_2 - r_2 s_1$ and the congruence condition $\Delta(z_1, z_2) \equiv 0 \pmod{4d}$ reads as

(15.1) $$r_1 s_2 \equiv r_2 s_1 \pmod{4d}.$$

For any $p | \Delta$ we have
$$z_2/z_1 \equiv r_2/r_1 \quad \text{if } p \nmid r_1,$$
$$z_2/z_1 \equiv s_2/s_1 \quad \text{if } p \nmid s_1.$$

Since $(r_1, s_1) = 1$ we always have one or both cases. Moreover we have $(d, r_1) = (d, r_2) = e$, say. Dividing (15.1) by $e$ and changing notation to remove the factor $e$ from $d, r_1, r_2$ we infer that $V(\beta)$ is given by

(15.2)
$$2 \sum\sum_{(z_1,z_2)=1} \beta_{z_1} \overline{\beta}_{z_2} \sum\sum_{\substack{ed>X \\ r_1 s_2 \equiv r_2 s_1 (4d)}} f\left(\frac{|\Delta|}{ed}\right) \frac{\varphi(ed)}{ed} \left(\frac{s_1 s_2}{e}\right) \left(\frac{r_1 r_2}{d}\right) \log 2 \left|\frac{z_1 z_2}{\Delta}\right|$$

where now $z_1 = er_1 + is_1$, $z_2 = er_2 + is_2$.

Our next goal is to separate the variables $z_1, z_2$. The condition $(z_1, z_2) = 1$ can be removed precisely by the Möbius formula, however we take advantage of the restriction (5.8) which allows us to ignore this condition at the admissible cost $O(P^{-1}N^2)$. To get it at this cost use Lemma 2.2 with $k = 4$ to reduce $d$ before estimating trivially. Here one naturally loses a few logarithms coming from the structure of the arithmetic functions involved but these are compensated by the saving of a large number of logarithms due to the size of the box so the bound is actually smaller than stated. This remark, which we shall not repeat, will later apply to a few 'trivial' estimates of a similar nature.

To separate the variables $z_1, z_2$ constrained by $f(|\Delta|/ed)$ we use the same technique as in Section 12. We employ the function $g$ whose graph is

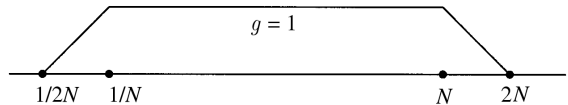

Put $B(x, y) = f(|x|y)g(|x|)$. The cutting factor $f(|\Delta|/ed)$ can be replaced in (15.2) by $B(x, y)$ at the special points

(15.3) $$x = \frac{s_1}{r_1} - \frac{s_2}{r_2} \quad \text{and} \quad y = \frac{|r_1 r_2|}{d}.$$

Now, using the Fourier-Mellin transform of $B(x,y)$ (see (12.13)), we separate the variables at the cost of the $L_1$-norm of the transform which is bounded by $O((\log N)^2)$; see (12.14).

To separate the variables in the logarithmic factor

(15.4) $$\log 2|z_1 z_2/\Delta| = -\log \tfrac{1}{2}|\sin(\alpha_2 - \alpha_1)|,$$

see (7.2), we use the Fourier series expansion

(15.5) $$-h(\alpha) \log \tfrac{1}{2}|\sin \alpha| = \sum_\ell c_\ell e^{i\ell\alpha} .$$

Here we have attached a "mollifier" function $h(\alpha)$ for the purpose of accelerating the convergence. We assume that $h(\alpha)$ is even, periodic of period $2\pi$ and vanishing at $\alpha = 0$. There are many good choices. For our purpose it suffices to take

(15.6) $$h(\alpha) = \min\left\{\|\tfrac{\alpha}{\pi}\|N, 1\right\}.$$

We have

(15.7) $$\sum_\ell |c_\ell| \ll (\log N)^2.$$

Note that for $\alpha = \alpha_2 - \alpha_1$ the mollifier $h(\alpha)$ does not alter the value (15.4) because $\pi N^{-1} < |\alpha| < \pi$. Inserting the Fourier series (15.5) in (15.2) we achieve the separation of variables at the cost of $O((\log N)^2)$.

After the separation of variables we get by (15.2)

$$V(\beta) \ll (\log N)^4 \sum\sum_{X<ed<Y} \sum_{a(\bmod d)} \left| \sum\sum_{ar \equiv s(4d)} \xi_{ers} \beta_{er+is}\left(\frac{r}{d}\right) \right|^2 + P^{-1}N^2$$

where $Y = NX^{-1}$ (this limit was redundant before the separation of variables) and $\xi_{ers}$ are complex numbers with $|\xi_{ers}| = 1$ (these are character values coming out of the separation process). Given $e$ we apply Proposition 14.1 to the sums over $d, a, r, s$ with $D, R, S$ such that $X < eD < Y$, $eR < \sqrt{2N}$ and $S < \sqrt{2N}$, showing that $V(\beta)$ is bounded by

$$(\log N)^8 \left\{X^{-\frac{1}{2}}N + Y^{\frac{1}{3}}N^{\frac{2}{3}} + N^{\frac{23}{24}+\varepsilon}\right\} \sum_r \sum_s \tau_3(r)|\beta_{r+is}|^2 + P^{-1}N^2 .$$

Here we estimate $|\beta_{r+is}|$ by $\tau(r^2 + s^2)$ but retain the range of summation in the box (5.14). Then we relax $\tau(r^2 + s^2)$ by applying Lemma 2.2 so that the summation over $s$ can be executed with sufficient accuracy. Finally summing over $r$ we arrive at

$$\sum_r \sum_s \tau_3(r)|\beta_{r+is}|^2 \ll \theta^2 N(\log N)^{1996} .$$

Hence we conclude

PROPOSITION 15.1. *For $1 \leqslant X \leqslant N^{\frac{1}{9}}$ we have*

$$V(\beta) \ll \left(P^{-1} + X^{-\frac{1}{3}}\right) N^2 . \tag{15.8}$$

## 16. Estimation of $U(\beta)$

In this section we estimate $U(\beta)$ using classical methods of prime number theory in the Gaussian domain. We write

$$U(\beta) = 2 \sum_{d \leqslant X} d^{-1} \varphi(d) U_d(\beta) \tag{16.1}$$

where

$$U_d(\beta) = \sum\sum_{\substack{(z_1, z_2) = 1 \\ \Delta(z_1, z_2) \equiv 0 \,(\mathrm{mod}\, 4d)}} f\left(\frac{|\Delta|}{d}\right) \beta_{z_1} \bar\beta_{z_2} \left(\frac{z_2/z_1}{d}\right) \log 2 \left| \frac{z_1 z_2}{\Delta} \right| . \tag{16.2}$$

Here the moduli are small and the source of cancellation will be the sign change of the Möbius function which is part of $\beta_z$. Therefore we can afford to estimate $U_d(\beta)$ individually for every $d \leqslant X$. In view of the multiplicative structure of $\beta_z$ it is natural to employ Hecke characters.

First note that the condition $\Delta(z_1, z_2) \equiv 0 \pmod{4d}$ is equivalent to

$$z_1 \equiv \omega z_2 \pmod{4d} \tag{16.3}$$

for some rational residue class $\omega(\mathrm{mod}\, 4d)$. Thus we can write

$$U_d(\beta) = \sum_{\omega(\mathrm{mod}\, 4d)} \left(\frac{\omega}{d}\right) \sum\sum_{\substack{(z_1, z_2) = 1 \\ z_1 \equiv \omega z_2 (\mathrm{mod}\, 4d)}} f\left(\frac{|\Delta|}{d}\right) \beta_{z_1} \bar\beta_{z_2} \log 2 \left|\frac{z_1 z_2}{\Delta}\right| . \tag{16.4}$$

Now we remove the condition $(z_1, z_2) = 1$ which costs us $O(d^{-1} P^{-1} N^2)$ by a trivial estimation using (5.8). Moreover by a trivial estimation we can remove the terms of (16.4) near the diagonal, say those with $|\alpha_2 - \alpha_1| < 2\pi H^{-1}$, at the cost of $O(d^{-1} H^{-1} N^2)$. This can be done smoothly by means of a function $h(\alpha)$ essentially as in (15.6); here specifically

$$h(\alpha) = \min\left\{\left\|\frac{\alpha}{\pi}\right\| H, 1\right\} . \tag{16.5}$$

We obtain

$$U_d(\beta) = \sum_{\omega(\mathrm{mod}\, 4d)} \left(\frac{\omega}{d}\right) \sum\sum_{z_1 \equiv \omega z_2 (\mathrm{mod}\, 4d)} \beta_{z_1} \bar\beta_{z_2} u(\alpha_1 - \alpha_2) \tag{16.6}$$
$$+ O\left(d^{-1}(P^{-1} + H^{-1}) N^2\right)$$

where (see (7.2))

(16.7) $$u(\alpha) = -h(\alpha)\log\tfrac{1}{2}|\sin\alpha|.$$

Note that the factor $f(|\Delta|/d)$ is not needed in (16.6). Next we expand $u(\alpha)$ into its Fourier series
$$u(\alpha) = \sum_k \hat{u}(k)e^{ik\alpha}$$
which converges quite rapidly. Indeed, we derive by (16.5) that

(16.8) $$\hat{u}(k) \ll (1 + k^2 H^{-2})^{-1}\log H.$$

Hence we can truncate the Fourier series with a small error term, namely

(16.9) $$u(\alpha) = \sum_{|k|\leqslant K}\hat{u}(k)e^{ik\alpha} + O(K^{-1}H^2\log H).$$

Inserting this in (16.4) we get

(16.10) $$U_d(\beta) = \sum_{|k|\leqslant K}\hat{u}(k)\sum_{\omega(\bmod 4d)}\left(\frac{\omega}{d}\right)\sum\sum_{z_1\equiv\omega z_2(\bmod 4d)}\beta_{z_1}\bar\beta_{z_2}\Big(\frac{z_1\bar z_2}{|z_1 z_2|}\Big)^k$$
$$+ O\left(d^{-1}(P^{-1} + H^{-1} + K^{-1}H^2)N^2\right).$$

Now we detect the congruence (16.3) by multiplicative characters of the group $(\mathbb{Z}[i]/4d\mathbb{Z}[i])^*$ getting

(16.11) $$\sum_\omega\left(\frac{\omega}{d}\right)\sum\sum_{z_1\equiv\omega z_2}\beta_{z_1}\bar\beta_{z_2}\Big(\frac{z_1\bar z_2}{|z_1 z_2|}\Big)^k = \frac{1}{\Phi(d)}\sum_\chi \mathcal{J}(\chi)|S_\chi^k(\beta)|^2$$

where $\Phi(d)$ is the order of the group which is given by

(16.12) $$\Phi(d) = 4d\varphi(4d)\prod_{p|d}\left(1 - \frac{\chi_4(p)}{p}\right),$$

and where

(16.13) $$\mathcal{J}(\chi) = \sum_{\omega\ (\bmod 4d)}\chi(\omega)\left(\frac{\omega}{d}\right).$$

Note that the sum $\mathcal{J}(\chi)$ is incomplete because $\omega$ runs over the subgroup of rational classes modulo $4d$ whereas $\chi$ is a character on the group of all classes modulo $4d$ in $\mathbb{Z}[i]$. For every such $\chi$ in (16.11) we have the character sum

(16.14) $$S_\chi^k(\beta) = \sum_z \beta_z\chi(z)(z/|z|)^k.$$

Inserting (16.11) in (16.10) and the latter in (16.1) we arrive at the formula

(16.15) $$U(\beta) = \sum_{|k|\leqslant K}\hat{u}(k)\sum_{d\leqslant X}\frac{2\varphi(d)}{d\Phi(d)}\sum_{\chi(\bmod 4d)}\mathcal{J}(\chi)|S_\chi^k(\beta)|^2$$
$$+ O\left((P^{-1} + H^{-1} + K^{-1}H^2)N^2\log N\right).$$

It remains to estimate the character sums $S_\chi^k(\beta)$. To this end we shall employ the theory of Hecke $L$-functions. The foundation of this theory was described by Hecke in the setting of a general number field but we shall take advantage of simplifications that come with the Gaussian field. Note that for $m \in \mathbb{Z}$ and $\chi$ a multiplicative character on the residue classes modulo $4d$ in $\mathbb{Z}[i]$ the function

(16.16) $$\psi(z) = \chi(z)(z/|z|)^m$$

defines a character on odd ideals of $\mathbb{Z}[i]$ by setting

$$\psi(\mathfrak{a}) = \psi(z)$$

where $z$ is the unique generator of $\mathfrak{a}$ which is primary.

We say that the Gaussian integer $z$ is odd if it is co-prime to $1+i$. Recall from (5.3) that we say $z$ is primary if $z \equiv 1 \pmod{2(1+i)}$ (see also (5.4) and (5.5)). Every odd $z$ is conjugate to exactly one primary integer; we denote this by $\hat{z}$. The only primary unit is $z = 1$. The product of primary numbers is primary and every primary number factors uniquely (up to permutation) as a product of primary numbers which are Gaussian primes.

To the character (16.16) we attach the $L$-function

(16.17) $$L(s, \psi) = \sum_{\mathfrak{a}} \psi(\mathfrak{a})(N\mathfrak{a})^{-s} .$$

Since $\psi$ is completely multiplicative $L$ has the Euler product

$$L(s, \psi) = \prod_{\mathfrak{p}} \left(1 - \psi(\mathfrak{p})(N\mathfrak{p})^{-s}\right)^{-1} .$$

This has a meromorphic continuation to $\mathbb{C}$, is entire apart from a simple pole at $s = 1$ in case of trivial $\psi$, which happens only if $m = 0$ and $\chi$ is trivial. There is also a functional equation for $\chi$ primitive proved by Hecke:

$$\left(\tfrac{4d}{2\pi}\right)^s \Gamma\left(s + \tfrac{1}{2}|m|\right) L(s, \psi) = W \left(\tfrac{4d}{2\pi}\right)^{1-s} \Gamma\left(1 - s + \tfrac{1}{2}|m|\right) L(1 - s, \bar{\psi})$$

where $W$ is a Gauss sum normalized so that $|W| = 1$. By the functional equation one derives, by applying the Phragmén-Lindelöf principle, crude but sufficient upper bounds for $L(s, \psi)$ and $L'(s, \psi)$. Namely we have

(16.18) $$L(s, \psi) \ll (d(|s| + |m|))^{1-\sigma} (\log(4d(|s| + |m|)))^2,$$

(16.19) $$L'(s, \psi) \ll (d(|s| + |m|))^{1-\sigma} (\log(4d(|s| + |m|)))^3,$$

in the strip $\tfrac{1}{2} \leqslant \sigma \leqslant 1$ where the implied constant is absolute.

We also need a zero-free region for $L(s, \psi)$. It follows from the above bounds by classical arguments that there are no zeros for $s = \sigma + it$ with

(16.20) $$\sigma > 1 - c/\log(4d + |m| + |t|)$$

where $c$ is a positive absolute constant, apart from a possible exceptional simple real zero in the case that $\psi$ is real.

If $\psi$ is real then $m = 0$ so $\psi$ is a real character on residue classes on $\mathbb{Z}[i]$ to modulus $4d$. One such character is given by

$$\psi(\mathfrak{a}) = \chi_q(N\mathfrak{a})$$

where $\chi_q$ is a real Dirichlet character to modulus $q = 4d$. In this case $L(s, \psi)$ is the product of two Dirichlet $L$-functions

$$L(s, \psi) = L(s, \chi_q)L(s, \chi_4\chi_q) \ .$$

This product is the automorphic $L$-function associated to a certain Eisenstein series for the group $\Gamma_0(q)$. There are a few other real characters on the group $\mathbb{Z}[i]/q\mathbb{Z}[i]$, which will be studied in Section 19. These give $L$-functions of cusp forms of weight one, level $q$ and the central character $\chi_4\chi_q$. Applying Siegel's method one shows

LEMMA 16.1.  *For any $0 < \varepsilon \leqslant \frac{1}{4}$ there exists a positive constant $c(\varepsilon)$ such that for every real character $\psi(\mathrm{mod}\, 4d)$ on $\mathbb{Z}[i]$ the Hecke $L$-function $L(s, \psi)$ does not vanish in the segment*

(16.21) $$s > 1 - c(\varepsilon)d^{-\varepsilon}.$$

*Proof.* We follow the simplified version of D. Goldfeld [Go] which applies in an extremely general context for $L$-functions attached to real characters. The constant $c(\varepsilon)$ is not computable from this proof nor from any other proof known to date.

Let $\psi_1(\mathrm{mod}\, 4d_1)$, $\psi_2(\mathrm{mod}\, 4d_2)$ be real characters on $\mathbb{Z}[i]$ such that $\psi_1$, $\psi_2$, and $\psi_1\psi_2$ are nontrivial. We consider the product of $L$-functions

$$L(s) = \zeta_K(s)L(s, \psi_1)L(s, \psi_2)L(s, \psi_1\psi_2)$$

where $\zeta_K(s) = \zeta(s)L(s, \chi_4)$ is the zeta-function of $\mathbb{Z}[i]$. By the Euler product we see that the Dirichlet series for $L(s)$ has nonnegative coefficients, specifically

$$L(s) = \exp \sum_{\mathfrak{p}} \sum_k \left(1 + \psi_1(\mathfrak{p}^k)\right)\left(1 + \psi_2(\mathfrak{p}^k)\right)(N\mathfrak{p})^{-ks} = \sum_1^\infty a_n n^{-s}$$

where $a_1 = 1$ and $a_n \geqslant 0$ for all $n$. At $s = 1$ there is a simple pole with

$$\mathrm{res}_{s=1}L(s) = \frac{\pi}{4}L(1, \psi_1)L(1, \psi_2)L(1, \psi_1\psi_2).$$

On the critical line we have by (16.18)

$$L(s) \ll (d_1 d_2 |s|)^2 \ .$$

Hence we derive by contour integration that

$$\sum_{1}^{\infty} a_n n^{-\beta} e^{-n/y} = \frac{1}{2\pi i} \int_{(2)} L(s+\beta)\Gamma(s) y^s ds$$

$$= \frac{\pi}{4} L(1,\psi_1) L(1,\psi_2) L(1,\psi_1\psi_2) \Gamma(1-\beta) y^{1-\beta}$$

$$+ L(\beta) + O\left(d_1^2 d_2^2 y^{\frac{1}{2}-\beta}\right)$$

where $\frac{1}{2} < \beta < 1$ and the error term bounds the integral on the line $\sigma = \frac{1}{2} - \beta$. On the left side the series is bounded below by its first term $e^{-1/y} > 1 - y^{-1}$. Therefore we have

$$\tfrac{\pi}{4} L(1,\psi_1) L(1,\psi_2) L(1,\psi_1\psi_2) \Gamma(1-\beta) y^{1-\beta} + L(\beta) > 1 + O\left(d_1^2 d_2^2 y^{\frac{1}{2}-\beta}\right).$$

Take $\beta = \beta_1$ any real zero of $L(s,\psi_1)$ in $\frac{3}{4} \leqslant \beta_1 < 1$ if it exists. Then $L(\beta_1) = 0$ and, choosing $y = c(d_1 d_2)^8$ where $c$ is a large absolute constant, we get

$$L(1,\psi_1) L(1,\psi_2) L(1,\psi_1\psi_2) \gg (1-\beta_1)(d_1 d_2)^{-8(1-\beta_1)}.$$

Since $L(1,\psi_1) \ll (\log 4d_1)^2$ and $L(1,\psi_1\psi_2) \ll (\log 4d_1 d_2)^2$ by (16.18) we conclude that

(16.22)  $$L(1,\psi_2) \gg (1-\beta_1)(d_1 d_2)^{-8(1-\beta_1)} (\log 4d_1 d_2)^{-4}$$

where the implied constant is absolute and can be effectively computed.

Now we are ready to complete the proof of the lemma. We argue as follows. Fix $0 < \varepsilon \leqslant \frac{1}{4}$. If, for every real character $\psi$, the $L(s,\psi)$ does not vanish in the segment $s > 1 - \varepsilon$ then we take $c(\varepsilon) = \varepsilon$ getting the result. Suppose then there exists a real character $\psi$ for which $L(s,\psi)$ vanishes at some point $s > 1 - \varepsilon$. Fix one such, say $\psi = \psi_1$ and let $\beta_1$ denote the largest such zero. Let $\psi_2$ be any nontrivial real character on $\mathbb{Z}[i]$. If $\psi_1\psi_2$ is trivial then $L(s,\psi_2)$ has the same zeros as $L(s,\psi_1)$ in $s > \frac{1}{2}$. In this case we take $c(\varepsilon) = 1 - \beta_1$ getting the result. If $\psi_1\psi_2$ is nontrivial then (16.22) yields the lower bound

$$L(1,\psi_2) \gg d_2^{-8\varepsilon} (\log 4d_2)^{-4}$$

where the implied constant is effectively computable in terms of $d_1$. On the other hand, if $\beta_2$ is a real zero of $L(s,\psi_2)$ in $s > 1 - \varepsilon$ then by the mean-value theorem and (16.19) we get the upper bound

$$L(1,\psi_2) = (1-\beta_2) L'(s,\psi_2) \ll (1-\beta_2) d_2^\varepsilon (\log 4d_2)^3.$$

Combining the upper and lower bounds we infer that $1-\beta_2 \gg d_2^{-10\varepsilon}$ completing the proof of Lemma 16.1 (just change $10\varepsilon$ into $\varepsilon$). □

The information accumulated in (16.18)–(16.21) suffices by classical arguments to establish the upper bound

(16.23)  $$L(s,\psi)^{-1} \ll d(\log(|m|+|t|+3))^2$$

in the region

(16.24) $$\sigma \geqslant 1 - c(\varepsilon)/d^\varepsilon \log(|m| + |t| + 3)$$

for any $\varepsilon > 0$ where $c(\varepsilon) > 0$ and the implied constant depends only on $\varepsilon$. To obtain such a bound for $1/L$ we may begin with the corresponding bound for $L'/L$, integrate it to the right to obtain a bound for $|\log L|$ and then exponentiate the latter.

Having the bound (16.23) we can now estimate each character sum

(16.25) $$S_\chi^k(\beta) = \sum_z \beta_z \psi(z)$$

by employing $L$-functions with Grossencharacters. Recall that $\beta_z$ is given by (5.13) subject to (5.7) and (5.8) (the restriction for the number of divisors (5.9) was already waived after (10.4) had been established). Here the congruence (5.7) can be detected by characters to modulus eight so we can relax it by changing $\chi$ in $\psi$, however we retain the property of $z$ being primary so as to keep the unique correspondence with ideals. The remaining restriction (5.8) that $z$ be free of small prime factors will be taken care of easily after (16.25) has been expressed in terms of $L$-series.

To this end we expand the "angle mollifier" $q(\alpha)$ into its Fourier series

$$q(\alpha) = \sum_\ell \hat{q}(\ell) e^{i\ell\alpha} .$$

By the properties of $q(\alpha)$ recorded just before (5.13) it follows that the Fourier coefficients are bounded by

(16.26) $$\hat{q}(\ell) \ll \theta(1 + \theta^2 \ell^2)^{-1} .$$

Thus we can truncate the Fourier series with a small error term

(16.27) $$q(\alpha) = \sum_{|\ell| \leqslant L} \hat{q}(\ell) e^{i\ell\alpha} + O(\theta^{-1} L^{-1}) .$$

This error term contributes to $S_\chi^k(\beta)$ at most $O(L^{-1} N)$ by a trivial estimation.

Next we write the "radius mollifier" $p(n)$ as the Mellin transform

$$p(n) = \frac{1}{2\pi i} \int_{(\sigma)} \hat{p}(s) n^{-s} ds$$

where $\frac{1}{2} \leqslant \sigma \leqslant 2$. Since $p(n)$ is supported on the segment (4.12) and its derivatives satisfy (4.14) it follows that $\hat{p}(s)$ is entire and bounded by

(16.28) $$\hat{p}(s) \ll \theta(1 + \theta^2 t^2)^{-1} N^\sigma .$$

Applying the above formulas to the character sum (16.25) we get

(16.29) $$S_\chi^k(\beta) = \sum_{|\ell| \leqslant L} \hat{q}(\ell) \frac{1}{2\pi i} \int_{(\sigma)} \hat{p}(s) Z_\chi^{k+\ell}(s) ds + O(L^{-1} N)$$

where
$$Z_\chi^m(s) = \sum_{(z,\Pi)=1} \chi(z)\bigl(z/|z|\bigr)^m \mu(n)\gamma(n) n^{-s},$$

$n=|z|^2$ and (see (2.12))
$$\gamma(n) = \sum_{c|n,\ c\leqslant C} \mu(c) .$$

We write the series for $Z_\chi^m(s)$ in terms of ideals (and, for notational simplicity, delete the scripts $\chi$, $m$ ) getting
$$Z(s) = \sum_{(\mathfrak{a},\Pi)=1} \psi(\mathfrak{a})\mu(N\mathfrak{a})\gamma(N\mathfrak{a})(N\mathfrak{a})^{-s} .$$

In the analytic theory of zeta-functions it is more convenient to work with Dirichlet series over rational integers. Therefore we group the ideals in accordance with the norm and associate with the Grossencharacter $\psi$ the Hecke eigenvalue

(16.30) $$\lambda(n) = \sum_{N\mathfrak{a}=n} \psi(\mathfrak{a}) = \sum_{z\bar{z}=n} \chi(z)\bigl(z/|z|\bigr)^m .$$

From the Euler product it follows that
$$\sum_\nu \lambda(p^\nu) T^\nu = \bigl(1 - \lambda(p)T + \varepsilon(p)T^2\bigr)^{-1}$$

where $\varepsilon(p) = 1$ if $p \equiv 1 \pmod 4$ and $\varepsilon(p) = 0$, otherwise. Now
$$\begin{aligned}
Z(s) &= \sum_{(n,\Pi)=1} \lambda(n)\mu(n)\gamma(n) n^{-s} \\
&= \sideset{}{^\flat}\sum_{\substack{c\leqslant C \\ (c,\Pi)=1}} \lambda(c) c^{-s} \sum_{(n,c\Pi)=1} \lambda(n)\mu(n) n^{-s} \\
&= \prod_{p>P} \bigl(1 - \lambda(p)p^{-s}\bigr) \sideset{}{^\flat}\sum_{\substack{c\leqslant C \\ (c,\Pi)=1}} \lambda(c) \prod_{p|c} \bigl(1 - \lambda(p)p^{-s}\bigr)^{-1} c^{-s} \\
&= M(s) P(s) Q(s) ,
\end{aligned}$$

say, where
$$M(s) = \prod_p \bigl(1 - \lambda(p)p^{-s}\bigr), \qquad P(s) = \prod_{p\leqslant P} \bigl(1 - \lambda(p)p^{-s}\bigr)^{-1},$$

and $Q(s)$ is the remaining sum over $c \leqslant C$. Since $M(s)L(s,\psi)$ is given by an Euler product which converges absolutely in $\operatorname{Re} s > \tfrac{1}{2}$ we get by (16.23)

(16.31) $$M(s) \ll d\bigl(\log(|m| + |t| + 3)\bigr)^2$$

in the range (16.24). For the estimation of $P(s)$ we assume in addition to (16.24) that

(16.32) $$1 - 1/\log P \leqslant \sigma \leqslant 1 .$$

In this region $|p^{-s}| = p^{-\sigma} < 3p^{-1}$ for all $p \leqslant P$ so we have the trivial bound

(16.33) $$P(s) \ll (\log P)^3 .$$

We also estimate $Q(s)$ trivially by

(16.34) $$Q(s) \ll C^{1-\sigma} \log C .$$

Multiplying the bounds (16.31), (16.33) and (16.34) we get

(16.35) $$Z(s) \ll dC^{1-\sigma} \left(\log(|m| + |t| + N)\right)^6$$

uniformly in the intersection of the regions (16.24) and (16.32).

Now we can estimate the integral

(16.36) $$I = \frac{1}{2\pi i} \int_{(\sigma)} \hat{p}(s) Z(s) ds .$$

We choose $T \geqslant |m| + 3$ and put

(16.37) $$\delta = \min \{c(\varepsilon)/d^{\varepsilon} \log T , \ 1/\log P\} .$$

We then move the contour of integration in (16.36) to the vertical segments

$$\begin{aligned} s &= 1 + it & \text{with } |t| \geqslant T, \\ s &= 1 - \delta + it & \text{with } |t| \leqslant T, \end{aligned}$$

and to the two connecting horizontal segments

$$s = \sigma \pm iT \quad \text{with } 1 - \delta \leqslant \sigma \leqslant 1 .$$

Integrating trivially along the above segments we infer by (16.28) and (16.35) that

$$\begin{aligned} I &\ll d \left(T^{-1} + (C/N)^{\delta}\right) N \left(\log(N+T)\right)^6 \\ &\leqslant d \left(T^{-1} + N^{-\eta\delta}\right) N \left(\log(N+T)\right)^6 \end{aligned}$$

since $C \leqslant N^{1-\eta}$ where $\eta$ is a positive constant from Proposition 4.1. Choosing $T = 3\exp(\sqrt{\log N})$ we get

(16.38) $$I \ll \left\{ \exp\left(-\eta c(\varepsilon) d^{-\varepsilon} \sqrt{\log N}\right) + \exp\left(-\eta \frac{\log N}{\log P}\right) \right\} dN (\log N)^6$$

uniformly for $|m| \leqslant 2\exp(\sqrt{\log N})$. By (16.38) and (16.29) we conclude, after summing over $\ell$ with $|\ell| \leqslant L = \exp(\sqrt{\log N})$ and using (16.26), that the same bound (16.38) holds true for the character sum (16.25). We state this as

LEMMA 16.2. *Let $\chi$ be a character to modulus $4d$ on $\mathbb{Z}[i]$. Then*

$$S_\chi^k(\beta) \ll \left\{ \exp\left(-\eta c(\varepsilon) d^{-\varepsilon} \sqrt{\log N}\right) + \exp\left(-\eta \frac{\log N}{\log P}\right) \right\} dN (\log N)^6$$

*uniformly for $|k| \leq \exp(\sqrt{\log N})$, the implied constant depending only on $\varepsilon$.*

Finally, by Lemma 16.2 and (16.15) we derive that $U(\beta)$ is bounded by

$$\left\{ \exp\left(-2\eta c(\varepsilon) X^{-\varepsilon} \sqrt{\log N}\right) + \exp\left(-2\eta \frac{\log N}{\log P}\right) \right\} HX^4 N^2 (\log HN)^{13}$$
$$+ \left(P^{-1} + H^{-1} + K^{-1}H^2\right) N^2 \log N$$

provided $K \leq \exp(\sqrt{\log N})$. We choose $H = X$, $K = X^3$ and restrict $X$ by

(16.39) $$\log X \ll \log \log N.$$

Recall that $P$ was already restricted by

(16.40) $$\log P \ll (\log N)(\log \log N)^{-2}.$$

By these choices and restrictions the above bound for $U(\beta)$ becomes:

PROPOSITION 16.3. *For $1 \leq X \leq (\log N)^\nu$ with any $\nu$ we have*

(16.41) $$U(\beta) \ll (P^{-1} + X^{-1}) N^2 \log N ,$$

*where the implied constant depends on $\nu$.*

*Remarks.* The restriction (16.40) can be weakened slightly, replaced by $\log P \ll (\log N)(\log \log N)^{-1}$ by an elementary but more sophisticated sieve method. However (16.39) cannot be weakened, given the current state of knowlege concerning exceptional zeros. For this reason the constant implied in (16.41) is not effectively computable nor is the one in our main Theorem 1. All the other arguments in this paper give effective results, such as Theorem 2.

The formula (16.15) is valid in any range of the moduli. Therefore one could attempt to use this also in the middle range to give a direct treatment of $V(\beta)$ by an appeal to the theory of the large sieve rather than by the synthesis of the three methods from Sections 11, 12, 13. However such an approach fails because it ignores the intrinsic nature of the factor $\mathcal{J}(\chi)$ (which is a partial character sum over the subgroup of rational residue classes given by (16.13)). The change of the argument of $\mathcal{J}(\chi)$ plays an important role. Note also that $|\mathcal{J}(\chi)|$ can be as large as $\varphi(d)$ if $\chi$ is a real character (one of those constructed in Section 19), while the average value is $\sqrt{\varphi(d)}$; precisely we have

(16.42) $$\frac{1}{\Phi(d)} \sum_{\chi \ (\mathrm{mod}\, d)} |\mathcal{J}(\chi)|^2 = \varphi(d) .$$

## 17. Transformations of $W(\beta)$

It remains to estimate the partial sum $W(\beta)$ of $T(\beta)$ over $d > |\Delta|/4X$. In addition to this lower bound, $4d$ divides the determinant $\Delta$ and, since $d$ is extremely large it is useful to switch $d$ into the complementary divisor of $|\Delta|/4$. We get (see (10.13) and (10.17))

$$W(\beta) = 2 \sum_d f^*(4d) \sum\sum_{\substack{(z_1,z_2)=1 \\ \Delta(z_1,z_2)\equiv 0(4d)}} \beta_{z_1}\bar{\beta}_{z_2} \left(\frac{z_2/z_1}{|\Delta|d}\right) \frac{4d}{|\Delta|} \varphi\left(\frac{|\Delta|}{4d}\right) \log 2 \left| \frac{z_1 z_2}{\Delta} \right|$$

because $\left(\frac{z_2/z_1}{|\Delta|/4d}\right) = \left(\frac{z_2/z_1}{|\Delta|d}\right)$. Here we write

$$\frac{4d}{|\Delta|}\varphi\left(\frac{|\Delta|}{4d}\right) = \sum_{4bd|\Delta} \frac{\mu(b)}{b},$$

so that $W(\beta)$ becomes

(17.1) $\quad 2 \sum_b \frac{\mu(b)}{b} \sum_d f^*(4d) \sum\sum_{\substack{(z_1,z_2)=1 \\ \Delta(z_1,z_2)\equiv 0(4bd)}} \beta_{z_1}\bar{\beta}_{z_2} \left(\frac{z_2/z_1}{|\Delta|d}\right) \log 2 \left|\frac{z_1 z_2}{\Delta}\right|.$

Note that from our restrictions for the support of $\beta_z$ (see Section 5) and the above summation condition we have

(17.2) $\qquad\qquad (z_1, \bar{z}_1) = (z_2, \bar{z}_2) = (z_1, z_2) = 1,$

(17.3) $\qquad\qquad z_1 \equiv z_2 \pmod{8},$

(17.4) $\qquad\qquad 0 < r_1 r_2 \equiv 1 \pmod{8}.$

The positivity of $r_1 r_2$ means that the points $z_1, z_2$ lie in the same half-plane, but in fact they even both belong to the same fixed narrow sector (5.12).

LEMMA 17.1. *Assuming (17.2)–(17.4) we have*

(17.5) $\qquad\qquad \left(\frac{z_2/z_1}{|\Delta|}\right) = \left(\frac{s_1}{|r_1|}\right)\left(\frac{s_2}{|r_2|}\right)$

*where* $\Delta = \Delta(z_1, z_2) = \operatorname{Im} \bar{z}_1 z_2 = r_1 s_2 - r_2 s_1$ *for* $z_1 = r_1 + is_1$ *and* $z_2 = r_2 + is_2$.

Note that $s_1, s_2$ are even and $(r_1, s_1) = (r_2, s_2) = 1$.

*Proof.* We shall appeal to the quadratic reciprocity law

(17.6) $\qquad\qquad \left(\frac{a}{|b|}\right)\left(\frac{b}{|a|}\right) = (-1)^{\frac{a-1}{2}\frac{b-1}{2}} (a,b)_\infty$

for $a, b$ odd and co-prime, where $\left(\frac{c}{d}\right)$ is the Jacobi symbol and $(a,b)_\infty$ is the Hilbert symbol at $\infty$ defined by

$$(17.7) \qquad (a,b)_\infty = \begin{cases} -1 & \text{if } a < 0, b < 0 \\ 1 & \text{otherwise.} \end{cases}$$

We have $z_2/z_1 = (r_1 r_2 + s_1 s_2 + i\Delta)|z_1|^{-2}$, whence

$$\left(\frac{z_2/z_1}{|\Delta|}\right) = \left(\frac{|z_1|^2}{|\Delta|}\right)\left(\frac{r_1 r_2 + s_1 s_2}{|r_1 s_2 - r_2 s_1|}\right).$$

Put $r_1 = \rho u_1$, $r_2 = \rho u_2$ with $(u_1, u_2) = 1$, so $(u_1 u_2, s_1 s_2) = 1$, $u_1 u_2 \equiv 1 \pmod 8$ and $u_1 u_2 > 0$. Putting $|z_1|^2 = r_1^2 + s_1^2 = \rho^2 u_1^2 + s_1^2$, we have

$$(\rho^2 u_1 u_2 + s_1 s_2) u_1^2 = u_1 u_2 |z_1|^2 + u_1 s_1 (u_1 s_2 - u_2 s_1),$$

and so

$$\left(\frac{r_1 r_2 + s_1 s_2}{|r_1 s_2 - r_2 s_1|}\right) = \left(\frac{s_1 s_2}{\rho}\right)\left(\frac{\rho^2 u_1 u_2 + s_1 s_2}{|u_1 s_2 - u_2 s_1|}\right) = \left(\frac{s_1 s_2}{\rho}\right)\left(\frac{u_1 u_2 |z_1|^2}{|u_1 s_2 - u_2 s_1|}\right).$$

Since $|z_1|^2$ is a square modulo $\rho$ we get

$$\left(\frac{z_2/z_1}{|\Delta|}\right) = \left(\frac{s_1 s_2}{\rho}\right)\left(\frac{u_1 u_2}{|u_1 s_2 - u_2 s_1|}\right) = \left(\frac{s_1 s_2}{\rho}\right)\left(\frac{u_1 s_2 - u_2 s_1}{u_1 u_2}\right)$$

by (17.6) (actually one has to pull out from the determinant the whole power of 2 before application of (17.6) and install it back after at no cost because of the convention (8.15) and the congruence $u_1 u_2 \equiv 1 \pmod 8$). Hence

$$\left(\frac{z_2/z_1}{|\Delta|}\right) = \left(\frac{s_1 s_2}{\rho}\right)\left(\frac{-u_2 s_1}{|u_1|}\right)\left(\frac{u_1 s_2}{|u_2|}\right) = \varepsilon \left(\frac{s_1}{|r_1|}\right)\left(\frac{s_2}{|r_2|}\right),$$

where

$$\varepsilon = \left(\frac{-u_2}{|u_1|}\right)\left(\frac{u_1}{|u_2|}\right) = 1$$

by (17.6) using $0 < u_1 u_2 \equiv 1 \pmod 8$. This completes the proof of (17.5). $\square$

Inserting (17.5) in (17.1) we arrive at

$$W(\beta) = 2 \sum_b \frac{\mu(b)}{b} \sum_d f^*(4d) \sum\sum_{\substack{(z_1, z_2) = 1 \\ \Delta(z_1, z_2) \equiv 0 (4bd)}} \beta'_{z_1} \bar{\beta}'_{z_2} \left(\frac{z_2/z_1}{d}\right) \log 2 \left|\frac{z_1 z_2}{\Delta}\right|$$

where

$$(17.8) \qquad \beta'_z = i^{\frac{r-1}{2}} \left(\frac{s}{|r|}\right) \beta_z \qquad \text{if } z = r + is.$$

Here we could introduce the factor $i^{\frac{r-1}{2}}$ because $r_1 \equiv r_2 \pmod 8$ by virtue of (17.3). This factor will turn out to be quite natural. Since the condition $\Delta(z_1, z_2) \equiv 0 \pmod{4bd}$ is equivalent to $z_1 \equiv \omega z_2 \pmod{4bd}$ for some rational $\omega \pmod{4bd}$ we have

$$(17.9) \qquad W(\beta) = 2 \sum_b \frac{\mu(b)}{b} \sum_d f^*(4d) \sum_{\omega \pmod{4bd}} \left(\frac{\omega}{d}\right)$$

$$\sum\sum_{\substack{(z_1, z_2)=1 \\ z_1 \equiv \omega z_2 \pmod{4bd}}} \beta'_{z_1} \bar{\beta}'_{z_2} \log 2 \left| \frac{z_1 z_2}{\Delta} \right| .$$

We reduce the summation to $b \leqslant X$ and estimate the remaining sum trivially by $O(X^{-1} N^2)$. Then we remove the condition $(z_1, z_2) = 1$ at the cost of $O(P^{-1} N^2)$. Next we mollify the logarithmic factor $\log 2 |z_1 z_2/\Delta|$ as in (16.7) estimating the contribution from terms near the diagonal by $O(H^{-1} N^2)$ and we insert the truncated Fourier series (16.9) at the cost of $O(K^{-1} H^2 N^2)$. We detect the congruence $z_1 \equiv \omega z_2 \pmod{4bd}$ by characters of $(\mathbb{Z}[i]/4bd\mathbb{Z}[i])^*$ getting

$$(17.10) \quad W(\beta) = 2 \sum_{|k| \leqslant K} \hat{u}(k) \sum_{b \leqslant X} \frac{\mu(b)}{b} \sum_d \frac{f^*(4d)}{\Phi(bd)} \sum_{\chi \pmod{4bd}} \mathcal{J}(\chi) |S_\chi^k(\beta')|^2$$

$$+ O\left( (X^{-1} + P^{-1} + H^{-1} + K^{-1} H^2) N^2 \log N \right)$$

where $\Phi(bd)$ is the order of the group $(\mathbb{Z}[i]/4bd\mathbb{Z}[i])^*$ (see the formula (16.12)), $\mathcal{J}(\chi)$ is the character sum over rational classes

$$(17.11) \qquad \mathcal{J}(\chi) = \sum_{\omega \pmod{4bd}} \chi(\omega) \left(\frac{\omega}{d}\right) ,$$

and

$$(17.12) \qquad S_\chi^k(\beta') = \sum_z \beta'_z \chi(z) (z/|z|)^k .$$

The formula (17.10) closely resembles (16.15), the only important difference being between the character sums (17.12) and (16.14) wherein $\beta$ is turned into $\beta'$. This extra spin of $\beta'$ will be vital in producing cancellation in $S_\chi^k(\beta')$. Recall that $\beta_z$ is given by (5.13) and it carries the support conditions $z \equiv 1 \pmod{2(1+i)}$ and $(z, \Pi) = 1$ (remember that (5.7) was detected by characters and (5.9) was waived in Section 10). We no longer need the angle mollifier $q(\alpha)$ in $\beta_z$. Thus we replace it by the truncated Fourier series (16.27) which changes $S_\chi^k(\beta')$ to

$$\sum_{|\ell| \leqslant L} \hat{q}(\ell) S_\chi^{k+\ell}(\beta') + O(L^{-1} N) .$$

Here the error term $O(L^{-1}N)$ gives by (16.42) a contribution to $W(\beta)$ which is at most $O(KX^2L^{-2}N^2\log X)$. We choose $H = X$ and $K = L = X^3$ getting

$$(17.13) \quad W(\beta) \ll (\log N) \sum_{|k|\leqslant 2X^3} \sum_{d\leqslant X^2} \sum_{\chi \pmod{4d}} |S_\chi^k(\beta')|^2 \\ + (X^{-1} + P^{-1})N^2 \log XN \ .$$

Now we have

$$(17.14) \quad S_\chi^k(\beta') = \sum_{(z,\Pi)=1}^{\wedge} \beta_z [z]\chi(z)\bigl(z/|z|\bigr)^k$$

where $\sum^{\wedge}$ restricts to primary numbers, the coefficients $\beta_z$ being given by

$$(17.15) \quad \beta_z = p(n)\mu(n) \sum_{c|n,\ c\leqslant C} \mu(c)$$

with $n = z\bar{z}$ and the symbol $[z]$ being defined by

$$(17.16) \quad [z] = i^{\frac{r-1}{2}} \left(\frac{s}{|r|}\right) \quad \text{if } z = r + is \ .$$

We shall estimate every sum $S_\chi^k(\beta')$ separately for each Grossencharacter

$$(17.17) \quad \psi(z) = \chi(z)\bigl(z/|z|\bigr)^k \ .$$

This time the source of cancellation is the symbol $[z]$ which is attached to $\beta_z$, not the character $\psi(z)$ nor the Möbius function $\mu(|z|^2)$. Therefore the conductor of $\psi$ plays a minor role. Our goal will be:

PROPOSITION 17.2.   *For every character* (17.17) *we have*

$$(17.18) \quad S_\chi^k(\beta') \ll N^{1-\delta}$$

*uniformly in* $d(|k|+1) \leqslant N^\delta$ *for some positive constant $\delta$.*

By Proposition 17.2 and (17.10) we at once derive

PROPOSITION 17.3.   *For $1 \leqslant X < N^\varepsilon$ we have*

$$(17.19) \quad W(\beta) \ll \bigl(P^{-1} + X^{-1}\bigr) N^2 \log N \ .$$

## 18. Proof of main theorem

First, assuming Proposition 17.2, we prove Proposition 10.2. Recall by (10.14) the decomposition $T(\beta) = U(\beta) + V(\beta) + W(\beta)$ in accordance with the size of the divisors of the determinant. These three parts were estimated respectively in Sections 16, 15, 17. Inserting the results (15.8), (16.41), (17.19) we obtain
$$T(\beta) \ll \left(P^{-1} + X^{-\frac{1}{3}}\right) N^2 \log N$$
uniformly in $X \leqslant (\log N)^\nu$ for any $\nu$ with an implied constant depending on $\nu$. Choosing $X = (\log N)^{3\sigma+3}$ we complete the proof of Proposition 10.2.

Assuming only Proposition 17.2 we have all the pieces finally in place for the proof of Theorem 1. Recall that this had already, in the early sections of the paper been reduced to the proof of the bound (4.23) for the bilinear form $\mathcal{B}^*(M, N)$ and to this end we need only the results up to Proposition 10.2. We re-trace the path we followed from (4.23) to Proposition 10.2 and the conditions imposed on the parameters $A$, $A'$, $B$, $\tau$, and $P$. We required that

(4.4) $$(\log \log x)^2 \leqslant \log P \leqslant (\log x)(\log \log x)^{-2}$$

and

(4.11) $$\tau \geqslant (\log x)^{A+124}.$$

From $\mathcal{B}^*(M, N)$ we passed to $\mathcal{B}(M, N)$ which introduced the additional restrictions

(5.11) $$B \geqslant \tfrac{3}{2}A + A'$$

and $P \geqslant (\log x)^{A+A'}$, the latter following already from (4.4).

From $\mathcal{B}(M, N)$ we passed by Cauchy's inequality to $\mathcal{D}(M, N)$ and then to $\mathcal{D}^*(M, N)$, the latter step introducing the new restrictions

(5.23) $$\tau \leqslant x^{\frac{1}{3}\eta},$$

and

(5.24) $$P \geqslant \tau^2 (\log x)^{2A+4A'+508}.$$

Next, from $\mathcal{D}^*(M, N)$ we passed to $\mathcal{D}_0(M, N)$ encountering the additional requirement

(9.12) $$\tau \leqslant (\log x)^{\frac{1}{2}B - \frac{7}{4}A - 2A' - \frac{1}{2}b}$$

where $b$ depends only on $\eta$ in Proposition 4.1. Note that this condition supersedes (5.23) and that, together with (4.4), it implies (5.24).

Finally $\mathcal{D}_0(M, N)$ is given by (10.10) where the error term requires the conditions

(10.11) $$\tau \geqslant (\log x)^{A+2^{20}},$$

which supersedes (4.11), and

(10.12) $$A' \geqslant 2A + 2^{20}.$$

The main term in (10.10) introduces no new conditions due to Proposition 10.2 and (4.4).

It remains to show the consistency of this collection of conditions. To this end we may for example take $A' = 2A + 2^{20}$ and $\tau = (\log x)^{A+2^{20}}$. If we take $B$ sufficiently large we then ensure the remaining conditions (5.11) and (9.12). This completes the proof of Theorem 1.

We still need to prove Proposition 17.2. We postpone this task to establish first a larger background so as to be able to include a proof of the more general Theorem $2^\psi$ stated in Section 26. This itself does not depend on the results already established. Some of the arguments of the proof have independent interest so we present these in considerable generality. The proof of Proposition 17.2 has essentially the same ingredients as that of Theorem $2^\psi$ plus a new combinatorial identity (see Proposition 24.2) which saves one third of the work by allowing us to make an appeal to Theorem $2^\psi$ itself.

## 19. Real characters in the Gaussian domain

Let $q$ be an odd integer and $\omega \pmod{q}$ a root of

$$(19.1) \qquad \omega^2 + 1 \equiv 0 \pmod{q}.$$

Therefore every prime factor of $q$ is congruent to 1 $(\bmod\, 4)$, so it splits in $\mathbb{Z}[i]$. The number of roots $\omega \pmod{q}$ is equal to $2^t$, where $t$ is the number of distinct prime factors of $q$. Given $\omega \pmod{q}$ we define a function on $\mathbb{Z}[i]$ by

$$(19.2) \qquad \xi(z) = \left(\frac{r + \omega s}{q}\right) \quad \text{if } z = r + is.$$

Clearly $\xi$ is periodic of period $q$, and it is multiplicative as well since

$$(r_1 + \omega s_1)(r_2 + \omega s_2) \equiv r_1 r_2 - s_1 s_2 + \omega(r_1 s_2 + r_2 s_1) \pmod{q}.$$

Therefore $\xi$ is a character on Gaussian integers modulo $q$; it is an extension of the Jacobi symbol on the rational integers

$$(19.3) \qquad \xi(r) = \left(\frac{r}{q}\right) \quad \text{if } r \in \mathbb{Z}.$$

When $z$ is a unit we have

$$(19.4) \qquad \xi(\pm 1) = 1$$

and

$$(19.5) \qquad \xi(\pm i) = \begin{cases} 1 & \text{if } q \equiv 1 \pmod{8} \\ -1 & \text{if } q \equiv 5 \pmod{8}. \end{cases}$$

Thus $\xi(z)$ is an even character and, if $q \equiv 1 \pmod 8$ then $\xi(z)$ is a function on ideals.

The roots of (19.1) correspond to the representations

(19.6) $$q = u^2 + v^2 \quad \text{with } (u, v) = 1,\ v \text{ even}.$$

These are given by

(19.7) $$\omega \equiv -v\bar{u} \pmod q.$$

Since

(19.8) $$\left(\frac{u}{q}\right) = \left(\frac{|u|}{q}\right) = \left(\frac{q}{|u|}\right) = \left(\frac{v^2}{|u|}\right) = 1$$

we also have

(19.9) $$\xi(z) = \left(\frac{ur - vs}{q}\right) \quad \text{if } z = r + is.$$

If we require $w = u + iv$ to be primary and primitive then the correspondence (19.7) between the roots of (19.1) and the representations (19.6) is one-to-one and the characters (19.2) can be written as

(19.10) $$\xi(z) = \left(\frac{\operatorname{Re} wz}{|w|^2}\right).$$

For $q = |w|^2$ prime these characters were considered by Dirichlet [Di]. We shall also denote (19.10) by

(19.11) $$\xi_w(z) = \left(\frac{z}{w}\right)$$

and we call this the Dirichlet symbol even when $w$ is not prime. No confusion should arise between the Jacobi symbol $\chi_q(r) = \left(\frac{r}{q}\right)$ and the Dirichlet symbol $\xi_w(z) = \left(\frac{z}{w}\right)$ because the latter is defined and will be used exclusively for $w$ which is primary and primitive (a rational integer $q$ is not primary primitive unless $q = 1$ in which case the two symbols coincide). Note that

$$\left(\frac{z}{w}\right) = 0 \quad \text{if and only if } (w, \bar{z}) \neq 1.$$

If both $w$ and $z$ are primary primitive then

(19.12) $$\left(\frac{z}{w}\right) = \left(\frac{w}{z}\right).$$

Indeed, if $wz$ is primitive then by (19.8)

$$\left(\frac{z}{w}\right)\left(\frac{w}{z}\right) = \left(\frac{\operatorname{Re} wz}{|wz|^2}\right) = 1,$$

(because the Jacobi symbol is multiplicative in the lower entry) whereas if not then $(w, \bar{z}) \neq 1$ in which case both sides of (19.12) vanish.

Applying the reciprocity formula (19.12) we deduce that for any $z, w_1, w_2$ primary primitive

$$\left(\frac{z}{w_1}\right)\left(\frac{z}{w_2}\right) = \left(\frac{w_1 w_2}{z}\right).$$

Note that $w_1 w_2$ is not necessarily primitive since $w_1, \overline{w}_2$ may have a nontrivial common factor. Let $e$ be the primary associate of $(w_1, \overline{w}_2)$. Pulling out the divisor $d = e\bar{e}$ we obtain the primary primitive number $w_1 w_2/d$. Then applying the reciprocity law (17.6) and (19.12) we infer that

$$(19.13) \qquad \left(\frac{z}{w_1}\right)\left(\frac{z}{w_2}\right) = \left(\frac{|z|^2}{d}\right)\left(\frac{z}{w_1 w_2/d}\right).$$

This is true for all $z$, as seen by reducing to the case of $z$ primary primitive. In particular, taking $w_1 = w_2 = w$ primary primitive so that $d = (w, \overline{w}) = 1$, the formula (19.13) gives

$$(19.14) \qquad \xi_w \xi_{\bar{w}} = \chi_q \circ N$$

where $q = w\overline{w}$, $\chi_q$ is the Jacobi symbol and $N$ is the norm map. Now using (19.14) for $q = d$ we can write (19.13) as

$$(19.15) \qquad \xi_{w_1} \xi_{w_2} = \xi_e \xi_{\bar{e}} \xi_{w_1 w_2/e\bar{e}}$$

where $e$ is the primary associate of $(w_1, \overline{w}_2)$. This is the multiplicativity law in the lower entry for the Dirichlet symbol.

Having $\xi(z) = \xi_w(z)$ as a function in $z \in \mathbb{Z}[i]$ one can define a character on odd ideals $\mathfrak{a} \subset \mathbb{Z}[i]$ by setting $\xi(\mathfrak{a}) = \xi(z)$ if $\mathfrak{a} = (z)$ with $z$ primary. Note that if $w \equiv 1 \pmod 4$ then $q = w\overline{w} \equiv 1 \pmod 8$, and the choice of the primary generator is not necessary. With $\xi$ we associate the $L$-function

$$(19.16) \qquad L(s, \xi) = \sum_{\mathfrak{a} \text{ odd}} \xi(\mathfrak{a})(N\mathfrak{a})^{-s}.$$

Since $\xi$ is completely multiplicative we have the Euler product

$$\begin{aligned} L(s, \xi) &= \prod_{\mathfrak{p} \text{ odd}} \left(1 - \xi(\mathfrak{p})(N\mathfrak{p})^{-s}\right)^{-1} \\ &= \prod_{p > 2} \left\{1 - \lambda(p)p^{-s} + \chi_q(p)p^{-2s}\right\}^{-1} \end{aligned}$$

where

$$\lambda(p) = \begin{cases} \xi(\pi) + \xi(\bar{\pi}) & \text{if } p = \pi\bar{\pi}, \pi \text{ primary,} \\ 0 & \text{if } p \equiv -1 \pmod 4. \end{cases}$$

If the character $\xi$ is nontrivial then the corresponding theta series

$$\theta(z, \xi) = \sum_{\mathfrak{a} \text{ odd}} \xi(\mathfrak{a}) e(zN\mathfrak{a})$$

is a cusp form of weight one, level $4q$ and character $\chi_4\chi_q$. Hence $L(s, \xi)$ is an entire function. Moreover if $\xi$ is primitive we have the functional equation

$$(19.17) \qquad \left(\frac{\sqrt{q}}{\pi}\right)^s \Gamma(s)L(s,\xi) = W\left(\frac{\sqrt{q}}{\pi}\right)^{1-s}\Gamma(1-s)L(1-s,\xi)$$

where $W$ is the normalized Gauss sum. We shall not use any property of $L(s,\xi)$ in this paper although doing so would strengthen some of the results.

## 20. Jacobi-Kubota symbol

In Section 17 we arrived at the sums of $\beta_z$ twisted by Grossencharacters $\psi(z)$ and the symbol $[z]$ over $z$ primary and primitive. However, it makes sense to define

$$(20.1) \qquad [z] = i^{\frac{r-1}{2}} \left(\frac{s}{|r|}\right)$$

for all $z = r + is \equiv 1 \pmod 2$. Note that $[z]$ vanishes if $z$ is not primitive. Here the factor $i^{\frac{r-1}{2}}$ is attached since it simplifies forthcoming relations. One could further refine $[z]$ by including the Hilbert symbol

$$(r,s)_\infty = \begin{cases} -1 & \text{if } r,s < 0 \\ 1 & \text{otherwise,} \end{cases}$$

which amounts to extending the Jacobi symbol to all odd moduli (positive or negative) by setting

$$(20.2) \qquad \left(\frac{s}{r}\right) = \left(\frac{s}{|r|}\right)(r,s)_\infty \, .$$

However, being unable to take advantage of such a refinement, in this paper we stay with (20.1).

We shall be interested in the multiplicative structure of $[z]$ where, "multiplicative" refers to that within the Gaussian ring rather than rational multiplication relative to the coordinates $r, s$. For this reason we refer to $[z]$ as the Jacobi-Kubota symbol rather than simply the Jacobi symbol.

Of course, $[z]$ is not multiplicative per se, yet it is nearly so, up to a factor which is the Dirichlet symbol.

LEMMA 20.1.  *If $w$ is primary primitive and $z \equiv 1 \pmod 2$ then*

$$(20.3) \qquad [wz] = \varepsilon[w][z]\left(\frac{z}{w}\right)$$

*with $\varepsilon = \pm 1$ depending only on the quadrants in which $w, z, wz$ are located. Precisely, if $w = u + iv$ and $z = r + is$ then (20.3) reads as*

$$(20.4) \qquad i^{\frac{a-1}{2}}\left(\frac{us+vr}{|ur-vs|}\right) = \varepsilon i^{\frac{u-1}{2}+\frac{r-1}{2}}\left(\frac{v}{|u|}\right)\left(\frac{s}{|r|}\right)\left(\frac{ur-vs}{u^2+v^2}\right)$$

where $a = \operatorname{Re} wz = ur - vs$ and $\varepsilon = \varepsilon(w, z)$ is given by

$$(20.5) \qquad \varepsilon = (u, v)_\infty (r, -v)_\infty \qquad \text{if } ur > vs,$$

$$(20.6) \qquad \varepsilon = (u, v)_\infty (-r, v)_\infty \qquad \text{if } ur < vs.$$

*Remarks.* The formula (20.3) is reminiscent of a similar rule for the theta-multiplier (cf. [Sh]). The presence of the Dirichlet symbol in (20.3) destroys the multiplicativity but this is just what enables successful estimation of the relevant bilinear forms (see the next three sections).

*Proof.* We give a direct proof of (20.3) using the quadratic reciprocity law (17.6) together with its supplement

$$(20.7) \qquad \left(\frac{2}{d}\right) = i^{\frac{d^2-1}{4}} \qquad \text{if } 2 \nmid d.$$

We also need the following easy properties of the Hilbert symbol

$$(20.8) \qquad (x, y) = (y, x)$$
$$(20.9) \qquad (v, xy) = (v, x)(v, y)$$
$$(20.10) \qquad (x, -y) = (x, y)\operatorname{sign} x$$
$$(20.11) \qquad (-x, -y) = -(x, y)\operatorname{sign} xy.$$

Put $wz = ur - vs + i(us + vr) = a + ib$, say. First, we consider the case $a > 0$. We can assume that $(u, v) = (r, s) = 1$, or else both sides of (20.4) vanish. Put $v = \rho v_1$ and $r = \rho r_1$ with $(v_1, r_1) = 1$. Then

$$\left(\frac{b}{a}\right) = \left(\frac{us + vr}{ur - vs}\right) = \left(\frac{us}{\rho}\right)\left(\frac{us + vr}{ur_1 - v_1 s}\right) = \left(\frac{us}{\rho}\right)\left(\frac{u^2 r_1/v_1 + vr}{ur_1 - v_1 s}\right)$$
$$= \left(\frac{us}{\rho}\right)\left(\frac{v_1 r_1}{ur_1 - v_1 s}\right)\left(\frac{u^2 + v^2}{ur_1 - v_1 s}\right) = \left(\frac{us}{\rho}\right)\left(\frac{v_1 r_1}{ur_1 - v_1 s}\right)\left(\frac{u^2 + v^2}{ur - vs}\right)$$

because $\left(\frac{u^2+v^2}{\rho}\right) = \left(\frac{u^2}{\rho}\right) = 1$. Since $0 < u^2 + v^2 \equiv 1 \pmod 4$ we have

$$\left(\frac{u^2 + v^2}{ur - vs}\right) = \left(\frac{ur - vs}{u^2 + v^2}\right) = \left(\frac{z}{w}\right)$$

by the reciprocity law, giving the last factor in (20.4). Therefore

$$\left(\frac{b}{a}\right) = \left(\frac{us}{\rho}\right)\left(\frac{v_1 r_1}{ur_1 - v_1 s}\right)\left(\frac{z}{w}\right)$$

and we need to compute the factors on the right. First, by the reciprocity (17.6) we find that

$$\left(\frac{r_1}{ur_1 - v_1 s}\right) = (-1)^{\frac{r_1-1}{2}\frac{ur_1-1}{2}}\left(\frac{ur_1 - v_1 s}{|r_1|}\right) = (\operatorname{sign} r)(-1)^{\frac{r_1-1}{2}\frac{u-1}{2}}\left(\frac{v_1 s}{|r_1|}\right).$$

Similarly, writing $v_1 = 2^\alpha v'$ with $v'$ odd, we compute by twice applying reciprocity that

$$\left(\frac{v_1}{ur_1 - v_1 s}\right) = (-1)^{\frac{v'-1}{2}\frac{ur_1-1}{2}}\left(\frac{ur_1}{|v'|}\right)\left(\frac{2^\alpha}{ur_1 - v_1 s}\right)$$

$$= (ur, v)_\infty \left(\frac{v'}{|ur_1|}\right)\left(\frac{2^\alpha}{ur_1 - v_1 s}\right).$$

In the last symbol we change $v'$ back to $v_1$ and we evaluate the resulting Jacobi symbol for $2^\alpha$ by means of (20.7) getting

$$\left(\frac{2^\alpha}{ur_1 - v_1 s}\right)\left(\frac{2^\alpha}{|ur_1|}\right) = (-1)^{\frac{\alpha vs}{4}} = (-1)^{\frac{vs}{4}}$$

because $(ur_1)^2(ur_1 - v_1 s)^2 - 1 \equiv 2v_1 s \equiv 2vs \pmod{16}$. Note also that we have $(\operatorname{sign} r)(ur, v)_\infty = (\operatorname{sign} r)(r, v)_\infty(u, v)_\infty = (r, -v)_\infty(u, v)_\infty = \varepsilon$. Gathering these results we obtain

$$\left(\frac{b}{a}\right) = \varepsilon(-1)^{\frac{r_1-1}{2}\frac{u-1}{2}+\frac{vs}{4}}\left(\frac{us}{\rho}\right)\left(\frac{v_1}{|u|}\right)\left(\frac{s}{|r_1|}\right)\left(\frac{z}{w}\right),$$

because the symbol $\left(\frac{v_1}{|r_1|}\right)$ appears twice and therefore annihilates itself. Now, again by the reciprocity (17.6) we transform

$$\left(\frac{us}{\rho}\right)\left(\frac{v_1}{|u|}\right)\left(\frac{s}{|r_1|}\right) = \left(\frac{u}{\rho}\right)\left(\frac{\rho}{|u|}\right)\left(\frac{v}{|u|}\right)\left(\frac{s}{|r|}\right)$$

$$= (-1)^{\frac{\rho-1}{2}\frac{u-1}{2}}\left(\frac{v}{|u|}\right)\left(\frac{s}{|r|}\right).$$

Hence

$$\left(\frac{b}{a}\right) = \varepsilon(-1)^\nu\left(\frac{v}{|u|}\right)\left(\frac{s}{|r|}\right)\left(\frac{z}{w}\right)$$

where the parity of $\nu$ is given by

$$\nu \equiv \frac{r_1-1}{2}\frac{u-1}{2} + \frac{\rho-1}{2}\frac{u-1}{2} + \frac{vs}{4} \equiv \left(\frac{\rho r - 1}{2} + \frac{\rho-1}{2}\right)\frac{u-1}{2} + \frac{vs}{4}$$

$$\equiv \left(\rho\frac{r+1}{2} - 1\right)\frac{u-1}{2} + \frac{vs}{4} \equiv \frac{r-1}{2}\frac{u-1}{2} + \frac{vs}{4}$$

$$\equiv -\frac{r-1}{2}\frac{u-1}{2} + \frac{vs}{4} = \frac{u-1}{4} + \frac{r-1}{4} - \frac{a-1}{4}.$$

This completes the proof of (20.4) in the case $a = ur - vs > 0$.

If $a < 0$ we apply (20.4) with $(r, s)_\infty$ changed to $(-r, -s)_\infty$. The left side changes by the factor

$$i^{-a}\left(\frac{-1}{|a|}\right) = i^{-a-|a|+1} = i$$

while the right side changes by the factor (see (20.10), (20.11))

$$i^{-r}\left(\frac{-1}{|r|}\right)(r, -v)_\infty(-r, -v)_\infty = i(r, -v)_\infty(-r, v)_\infty .$$

These changes lead to (20.4) with $\varepsilon$ switched from (20.5) to (20.6). This completes the proof of Lemma 20.1. □

We shall have a problem to separate $w, z$ in $\varepsilon(w, z)$. In this connection one can use the formula

(20.12) $$2(x, y)_\infty = 1 + \text{sign}\, x + \text{sign}\, y - \text{sign}\, xy$$

to write

(20.13) $$2(r, -v)_\infty = 1 + \text{sign}\, vr + \text{sign}\, r - \text{sign}\, v$$

(20.14) $$2(-r, v)_\infty = 1 + \text{sign}\, vr - \text{sign}\, r + \text{sign}\, v .$$

Hence one can write both cases (20.5), (20.6) in a single form:

(20.15) $$2\varepsilon(w, z)(u, v)_\infty = 1 + \text{sign}(vr) - (\text{sign}\, v - \text{sign}\, r)\text{sign}(\text{Re}\, wz) .$$

*Remarks.* If $w, z$ are primary primitive we know by (20.3) and (19.12) that $\varepsilon(w, z)$ is symmetric, however this property is not so apparent from (20.15).

The formula (20.15) reduces the problem of the separation of variables in $\varepsilon(w, z)$ to that in $\text{sign}(\text{Re}\, wz)$, and the latter requires an application of harmonic analysis. One can proceed in a number of ways. We shall use the characters $e^{ik \arg z}$ since these fit well the analysis previously employed. Note that $\text{sign}(\text{Re}\, z) = 1$ is equivalent to $|\arg z| < \frac{\pi}{2}$, whence

(20.16) $$\text{sign}(\text{Re}\, wz) = t(\arg wz)$$

where $t(\alpha)$ is the periodic function of period $2\pi$ which is even, takes value 1 if $0 \leqslant \alpha < \frac{\pi}{2}$ and $-1$ if $\frac{\pi}{2} < \alpha \leqslant \pi$. Since $wz$ is not purely imaginary because of $wz \equiv 1 \pmod 2$, we can smooth $t(\alpha)$ slightly at $\alpha = \frac{\pi}{2}$ without altering the values (20.16) but getting the Fourier expansion

(20.17) $$t(\alpha) = \sum_k \hat{t}(k) e^{ik\alpha}$$

which converges absolutely and has $\ell_1$-norm of its coefficients almost bounded. Precisely, if $w, z$ are restricted by $|wz| \leqslant R$ and $wz \equiv 1 \pmod 2$ then $wz$ stays away from the imaginary axis by an angle $R^{-1}$ and this is sufficient room for the modification of $t(\alpha)$ so that

$$\text{(20.18)} \qquad \sum_k |\hat{t}(k)| \ll \log 2R \ .$$

By (20.16) and (20.17) we write

$$\text{(20.19)} \qquad \text{sign}(\text{Re } wz) = \sum_k \hat{t}(k) \left(\frac{wz}{|wz|}\right)^k \ .$$

Inserting (20.19) into (20.15) we obtain an expression for $\varepsilon(w, z)$ in which the variables $w$, $z$ appear separately.

## 21. Bilinear forms in Dirichlet symbols

In this section we establish (among other things) estimates for general sums of type

$$\text{(21.1)} \qquad Q(M, N) = \sum_w^* \sum_z \alpha_w \beta_z \left(\frac{z}{w}\right)$$

where $\alpha_w$, $\beta_z$ are complex coefficients supported in the discs

$$\text{(21.2)} \qquad |w|^2 \leqslant M$$
$$\text{(21.3)} \qquad |z|^2 \leqslant N$$

and the $*$ restricts the summation over $w$ to the primary primitive numbers. To simplify the arguments we assume, as we may, that

$$\text{(21.4)} \qquad |\alpha_w| \leqslant 1$$
$$\text{(21.5)} \qquad |\beta_z| \leqslant 1 \ .$$

Our aim is to improve on the trivial bound

$$\text{(21.6)} \qquad Q(M, N) \ll MN \ .$$

For the sake of simplicity we do not attempt to show the strongest results. The method is robust and it deserves a more precise consideration in a separate project. In particular one could employ the theory of the $L$-function (19.16) but we choose direct arguments which are sufficiently powerful. We begin by the following

LEMMA 21.1. *Let $w_1, w_2$ be primary primitive. Put $q = |w_1 w_2|^2$ and $d = |(w_1, \overline{w}_2)|^2$, so $d^2 \mid q$. Then we have*

$$\text{(21.7)} \qquad \sum_{\zeta \,(\text{mod } q)} \left(\frac{\zeta}{w_1}\right) \left(\frac{\zeta}{w_2}\right) = \begin{cases} q\varphi(d)\varphi(q/d) & \textit{if } q \textit{ and } d \textit{ are squares} \\ 0 & \textit{otherwise.} \end{cases}$$

*Proof.* By (19.13) we have

$$\left(\frac{\zeta}{w_1}\right)\left(\frac{\zeta}{w_2}\right) = \left(\frac{|\zeta|^2}{d}\right)\left(\frac{\zeta}{w_1 w_2/d}\right).$$

Hence the sum is equal to

$$\sum\sum_{r,s \;(\mathrm{mod}\, q)} \left(\frac{r^2+s^2}{d}\right)\left(\frac{r+\omega s}{q/d^2}\right) = \sum\sum_{r,s \;(\mathrm{mod}\, q)} \left(\frac{r-\omega s}{d}\right)\left(\frac{r+\omega s}{q/d}\right)$$

for some $\omega^2 + 1 \equiv 0 \;(\mathrm{mod}\, q)$. We change the variables $r, s$ into $r - \omega s \equiv x$ and $r + \omega s \equiv y$ getting

$$\sum_{x \;(\mathrm{mod}\, q)} \left(\frac{x}{d}\right) \sum_{y \;(\mathrm{mod}\, q)} \left(\frac{y}{q/d}\right) = q \sum_{x \;(\mathrm{mod}\, d)} \left(\frac{x}{d}\right) \sum_{y \;(\mathrm{mod}\, q/d)} \left(\frac{y}{q/d}\right) = 0$$

unless $d$ and $q/d$ are squares in which case the sum equals to $q\varphi(d)\varphi(q/d)$, completing the proof of (21.7). □

LEMMA 21.2. *We have*

$$(21.8) \qquad Q(M,N) \ll \left\{M^3 N^{\frac{1}{2}} + M^2 N^{\frac{3}{4}} + M^{\frac{1}{2}} N\right\}(MN)^\varepsilon.$$

*Proof.* By Cauchy's inequality

$$|Q(M,N)|^2 \leq \|\beta\|^2 \sum_z \Big| \sum_w{}^* \alpha_w \left(\frac{z}{w}\right) \Big|^2$$

$$= \|\beta\|^2 \sum_{w_1}{}^* \sum_{w_2}{}^* \alpha_{w_1} \bar{\alpha}_{w_2} \sum_z \left(\frac{z}{w_1}\right)\left(\frac{z}{w_2}\right).$$

Here and below $z$ runs over Gaussian integers in the disc (21.3). Splitting the inner summation into residue classes $\zeta \;(\mathrm{mod}\, q)$ with $q = |w_1 w_2|^2$ we get by elementary counting that

$$\sum_z \left(\frac{z}{w_1}\right)\left(\frac{z}{w_2}\right) = \sum_{\zeta \;(\mathrm{mod}\, q)} \left(\frac{\zeta}{w_1}\right)\left(\frac{\zeta}{w_2}\right) \left\{\frac{\pi N}{q^2} + O\left(\frac{\sqrt{N}}{q} + 1\right)\right\}.$$

Hence by Lemma 21.1 we get

$$Q(M,N)^2 \ll N^2 \sum\sum_{\substack{m_1, m_2 \leq M \\ m_1 m_2 = \square}} \tau(m_1 m_2) + NM^4(\sqrt{N} + M^2)$$

which yields (21.8). □

The estimate (21.8) is trivial if $N < M^4$. We shall refine this by exploiting the multiplicativity of $\xi_w(z)$ in $z$. Applying Hölder's inequality we get

$$Q^k(M,N) \ll M^{k-1} \sum_w{}^* \Big| \sum_z \beta_z \left(\frac{z}{w}\right) \Big|^k = M^{k-1} \tilde{Q}(M, N^k)$$

where $\tilde{Q}$ is the bilinear form of type (21.1) with coefficients
$$\tilde{\beta}_z = \sum_{z_1\ldots z_k=z} \beta_{z_1}\ldots\beta_{z_k}$$
and some $\tilde{\alpha}_w$ with $|\tilde{\alpha}_w| = |\alpha_w|$. Since $\tilde{\beta}_z \ll N^\varepsilon$ we get by applying (21.8) that
$$Q^k(M,N) \ll M^{k-1}\left\{M^3 N^{\frac{k}{2}} + M^2 N^{\frac{3k}{4}} + M^{\frac{1}{2}} N^k\right\}(MN)^\varepsilon ,$$
whence
$$Q(M,N) \ll \left\{M^{1+\frac{2}{k}}N^{\frac{1}{2}} + M^{1+\frac{1}{k}}N^{\frac{3}{4}} + M^{1-\frac{1}{2k}}N\right\}(MN)^\varepsilon .$$
This holds for any integer $k \geq 1$ and any $\varepsilon > 0$ with the implied constant depending only on $k$ and $\varepsilon$. Since $k$ can be arbitrarily large we have already established a nontrivial bound for $Q(M,N)$ whenever $M$ and $N$ have the same order of magnitude in the logarithmic scale. The above result would be sufficient for a proof of Theorem 1 but not of Theorem 2. However, it is possible to do better by exploiting a symmetry of $Q(M,N)$. With this aim in mind we choose $k = 6$, getting

COROLLARY.    *We have*
$$(21.9) \qquad Q(M,N) \ll \left\{M^{\frac{4}{3}}N^{\frac{1}{2}} + M^{\frac{7}{6}}N^{\frac{3}{4}} + M^{\frac{11}{12}}N\right\}(MN)^\varepsilon$$
*where the implied constant depends only on $\varepsilon$.*

Still (21.9) is trivial if $N < M^{\frac{2}{3}}$. Now we refine (21.9) by exploiting the reciprocity law (19.12). This requires both $w$ and $z$ to be primary primitive. Let $Q^*(M,N)$ denote the form (21.1) restricted to primary primitive $w$ and $z$. Of course, (21.9) holds for $Q^*(M,N)$. Interchanging $w$ with $z$ we switch $M$ with $N$ in (21.9). Then taking the minimum of the two bounds (21.9) we deduce that
$$Q^*(M,N) \ll \left(MN^{\frac{11}{12}} + NM^{\frac{11}{12}}\right)(MN)^\varepsilon .$$
This estimate extends to $Q(M,N)$ through the following arrangement:
$$Q(M,N) = \sum_d \sum_w^* \sum_z^* \alpha_w \beta_{dz}\left(\frac{d}{w\overline{w}}\right)\left(\frac{z}{w}\right) = \sum_d Q^*\left(M, Nd^{-2}\right) ,$$
say, where $Q^*(M, Nd^{-2})$ has coefficients depending on $d$ and bounded by 1. Hence we conclude

PROPOSITION 21.3.    *For any complex coefficients $\alpha_w, \beta_z$ supported in the discs (21.2), (21.3) and satisfying (21.4), (21.5) respectively we have*
$$(21.10) \qquad Q(M,N) \ll (M+N)^{\frac{1}{12}}(MN)^{\frac{11}{12}+\varepsilon}$$
*where the implied constant depends only on $\varepsilon$.*

Note that the bound (21.10) is nontrivial whenever $M$ and $N$ have the same order of magnitude in the logarithmic scale.

Our actual goal is to estimate bilinear forms in the Jacobi-Kubota symbol $[wz]$ (see (20.1)) rather than in the Dirichlet symbol $\left(\frac{z}{w}\right)$ (see (19.10)). However, by virtue of the multiplier rule (20.3) together with the formula (20.15) and the Fourier series (20.19) these problems are essentially equivalent as long as the coefficients $\alpha_w, \beta_z$ are arbitrary but bounded. We denote

$$(21.11) \qquad \mathcal{K}(M,N) = \sum_w \sum_z \alpha_w \beta_z [wz]$$

where $\alpha_w, \beta_z$ are complex coefficients supported on primary numbers in the discs (21.2), (21.3) respectively. We shall also need the restricted bilinear form

$$(21.12) \qquad \mathcal{K}^*(M,N) = \sum\sum_{(w,z)=1} \alpha_w \beta_z [wz] \;.$$

PROPOSITION 21.4. *For any complex coefficients $\alpha_w, \beta_z$ supported on primary numbers in the discs* (21.2), (21.3) *and satisfying* (21.4), (21.5) *respectively we have*

$$(21.13) \qquad \mathcal{K}(M,N) \ll (M+N)^{\frac{1}{12}} (MN)^{\frac{11}{12}+\varepsilon}$$

*where the implied constant depends only on $\varepsilon$. The same result holds true for the restricted bilinear form $\mathcal{K}^*(M,N)$.*

*Proof.* By Proposition 21.3 (without loss of generality one can restrict $w$ and $z$ to primitive numbers because otherwise $[wz]$ vanishes) we at once obtain (21.13). It remains to derive the same bound for $\mathcal{K}^*(M,N)$. Detecting the condition $(w,z)=1$ by Möbius inversion we obtain

$$\mathcal{K}^*(M,N) = \sum_e \mu(e) \sum_w \sum_z \alpha_{ew} \beta_{ez} [e^2 wz] \;.$$

Here we can assume that $e^2$ is primitive since otherwise $[e^2 wz]$ vanishes. Then by (20.3) we write

$$[e^2 wz] = \varepsilon [e^2][wz] \left(\frac{w}{e^2}\right) \left(\frac{z}{e^2}\right)$$

where $\varepsilon = \pm 1$ depends only on the quadrants in which $e^2 wz$, $e^2$, $wz$ are located. For each $e$ we can separate the variables $w, z$ in the $\varepsilon$–factor by applying (20.15) and (20.19) which will cost us $\log N$ by virtue of (20.18). In this way we obtain bilinear forms of type $\mathcal{K}(M|e|^{-2}, N|e|^{-2})$. These are estimated by $O\left(|e|^{-\frac{23}{6}} (M+N)^{\frac{1}{12}} (MN)^{\frac{11}{12}+\varepsilon}\right)$ due to (21.13). Summing this bound over $e$ we conclude that (21.13) holds for $\mathcal{K}^*(M,N)$. □

## 22. Linear forms in Jacobi-Kubota symbols

In addition to the general bilinear forms $\mathcal{K}(M,N)$ and $\mathcal{K}^*(M,N)$ we shall need bounds for the special linear forms

$$\mathcal{K}(N) = \sideset{}{^\wedge}\sum_{z \in \mathfrak{B}} \psi(z)[wz] \tag{22.1}$$

$$\mathcal{K}^*(N) = \sideset{}{^\wedge}\sum_{\substack{z \in \mathfrak{B} \\ (z,w)=1}} \psi(z)[wz] \,. \tag{22.2}$$

where $\wedge$ restricts to the primary numbers, $\mathfrak{B}$ is a polar box contained in the disc $|z|^2 \leq N$, $w$ is a fixed primary number and $\psi$ is a Hecke character given by (17.17), namely

$$\psi(z) = \chi(z)\bigl(z/|z|\bigr)^k \,.$$

Here $\chi$ is a character on residue classes to modulus $4d$ and $k$ is a rational integer.

Note the trivial bound $\mathcal{K}(N) \ll N$. We seek a bound which is nontrivial uniformly in large ranges of the relevant parameters $d$, $|k|$, $|w|$. We succeeded to establish

PROPOSITION 22.1.    *Given $\psi$ and $w$ as above we have*

$$\mathcal{K}(N) \ll d(|k|+1)|w|N^{\frac{3}{4}} \log|w|N, \tag{22.3}$$

$$\mathcal{K}^*(N) \ll d(|k|+1)|w|\tau(|w|^2)N^{\frac{3}{4}} \log|w|N \,. \tag{22.4}$$

*where the implied constant is absolute.*

*Proof.* We use the Polyá-Vinogradov theorem which asserts that

$$\sum_{n \leq N} \chi(n) \ll \sqrt{q}\,\log q \tag{22.5}$$

for any nontrivial Dirichlet character $\chi \pmod{q}$ where the implied constant is absolute.

For the proof of (22.3) we can assume that $w$ is primitive or else $[wz]$ vanishes for all $z$, hence so does $\mathcal{K}(N)$. By Lemma 20.1 we write

$$\mathcal{K}(N) = \varepsilon[w] \sideset{}{^\wedge}\sum_{z \in \mathfrak{B}} \psi(z)[z]\xi_w(z) \,. \tag{22.6}$$

Here we have pulled out the factor $\varepsilon = \pm 1$ which is permissible since it can be made constant by splitting $\mathfrak{B}$ into sixteen sectors (if necessary) to keep the arguments of $z$ and $wz$ in fixed quadrants. Write the real character $\xi_w$ as (see (19.2))

$$(22.7) \qquad \xi_w(z) = \left(\frac{r+\omega s}{q}\right) \quad \text{if } z = r+is$$

where $q = |w|^2$ and $\omega^2 + 1 \equiv 0 \pmod{q}$. Inserting (17.17), (20.1) and (22.7) in (22.6) we get

$$\mathcal{K}(N) = \varepsilon[w] \sideset{}{^\wedge}\sum_{r+is\in\mathfrak{B}} i^{\frac{r-1}{2}} \chi(r+is)\left(\frac{r+is}{|r+is|}\right)^k \left(\frac{r+\omega s}{q}\right) \left(\frac{s}{|r|}\right).$$

Observe that $\left(\frac{r+\omega s}{q}\right) = \left(\frac{\omega}{q}\right)\left(\frac{s-\omega r}{q}\right)$. For given $r$ we translate $s$ by $\omega r$ getting

$$|\mathcal{K}(N)| \leqslant \sum_{\substack{|r|<\sqrt{N} \\ r\equiv 1(\text{mod }2)}} \left| \sum_{s\in I(r)} \chi(s+\omega r-ir)\left(\frac{s+\omega r-ir}{|s+\omega r-ir|}\right)^k \left(\frac{s}{q|r|}\right) \right|$$

where $I(r)$ is a segment having length $< 2\sqrt{N}$ of integers in the progression $s \equiv r-1-\omega r \pmod{4}$. The inner sum satisfies

$$\sum_{s\in I(r)} \ll d(|k|+1)\sqrt{q|r|}\log q|r|$$

by the Polyá-Vinogradov estimate (22.5) provided that $q|r|/(d,q|r|)$ is not a square. Here the factor $d$ is lost by splitting $s$ into residue classes modulo $4d$ which is necessary to fix the values of $\chi$, whereas the factor $|k|+1$ is lost from an application of partial summation which is needed to remove the sector character. The condition that $q|r|/(d,q|r|)$ is not a square ensures that the relevant Jacobi symbol is not the trivial character. If it is a square we simply use the trivial bound $2\sqrt{N}$. Summing these bounds over $r$ we obtain (22.3). $\square$

Now we apply (22.3) to estimate the reduced linear form $\mathcal{K}^*(N)$. Detecting the condition $(z,w) = 1$ by Möbius inversion we get

$$\mathcal{K}^*(N) = \sideset{}{^\wedge}\sum_{e|w} \mu(e)\psi(e) \sideset{}{^\wedge}\sum_{z\in\mathfrak{B}/e} \psi(z)[ewz].$$

Hence we obtain (22.4).

*Remarks.* The bound (22.3) improves the trivial bound $\mathcal{K}(N) \ll N$ if

$$(22.8) \qquad d(|k|+1)|w| \ll N^{\frac{1}{4}}(\log N)^{-1}.$$

One can improve the range of uniformity (22.8) considerably by applying the well known estimate of Burgess [Bu] in place of (22.5), but we do not need such a strong result here (we did use Burgess' estimate in an earlier version of this work).

## 23. Linear and bilinear forms in quadratic eigenvalues

To the Hecke character (17.17) we associate the "quadratic eigenvalues"

$$\lambda(n) = \sum_{z\bar{z}=n}^{\wedge} \psi(z)[z] \tag{23.1}$$

where, as before, $\wedge$ restricts the summation to primary numbers. Writing $z = r + is$ this becomes

$$\lambda(n) = \sum_{r^2+s^2=n}^{\wedge} i^{\frac{r-1}{2}} \psi(r+is) \left(\frac{s}{|r|}\right) \tag{23.2}$$

where naturally enough, in this sum $\wedge$ restricts to pairs of integers satisfying (5.4) and (5.5). Note that $\lambda(n)$ vanishes if $n$ has prime factors other than $p \equiv 1 \pmod{4}$. We suspect, but have not examined thoroughly, that $\lambda(n)$ are related to the Fourier coefficients of some kind of metaplectic Eisenstein series or a cusp form, by analogy to the Hecke eigenvalues (16.30) which generate a modular form of integral weight. In the quadratic eigenvalue (23.1) we focus on the symbol $[z]$, yet the presence of the Hecke character $\psi$ offers welcome flexibility by means of which one can create simpler objects such as

$$\lambda_0(n) = \sum_{r^2+s^2=n} \left(\frac{s}{r}\right) \tag{23.3}$$

where $r, s$ are both positive and $r$ is odd. Here one can also put $r, s$ into a prescribed sector and one can require these to be in fixed residue classes to a given modulus. Theorem 2 concerns the eigenvalue (23.1) stripped to (23.3).

By analogy to $\mathcal{K}(M, N)$ we construct the bilinear form

$$\mathcal{L}(M, N) = \sum_m \sum_n \alpha(m)\beta(n)\lambda(cmn) \tag{23.4}$$

where $\alpha(m), \beta(n)$ are complex coefficients with

$$|\alpha(m)| \leqslant 1 \quad \text{for } 1 \leqslant m \leqslant M \tag{23.5}$$
$$|\beta(n)| \leqslant 1 \quad \text{for } 1 \leqslant n \leqslant N \tag{23.6}$$

and $c$ is a positive integer. Have in mind that $\lambda$ is not multiplicative so the introduction of $c$ makes the sum (23.4) more general and in practice this offers some extra flexibility. We shall also consider the restricted bilinear form

$$\mathcal{L}^*(M, N) = \sum\sum_{(m,n)=1} \alpha(m)\beta(n)\lambda(cmn) \ . \tag{23.7}$$

From Proposition 21.4 we derive

PROPOSITION 23.1.  *For any complex coefficients $\alpha(m), \beta(n)$ satisfying* (23.5), (23.6) *respectively, and for any $c \geqslant 1$ we have*

$$(23.8) \qquad \mathcal{L}(M,N) \ll \tau(c)(M+N)^{\frac{1}{12}}(MN)^{\frac{11}{12}+\varepsilon}$$

*where the implied constant depends only on $\varepsilon$. The same bound holds true for* $\mathcal{L}^*(M,N)$.

*Proof.* First we establish (23.8) for the reduced form $\mathcal{L}^*(M,N)$ with the coefficients $\alpha(m), \beta(n)$ supported on numbers prime to $c$. Since $c, m, n$ are mutually co-prime we have

$$\lambda(cmn) = \sideset{}{^\wedge}\sum_{e\bar{e}=c} \psi(e) \sideset{}{^\wedge}\sum_{w\bar{w}=m} \psi(w) \sideset{}{^\wedge}\sum_{z\bar{z}=n} \psi(z) [ewz] \ .$$

Hence

$$\mathcal{L}^*(M,N) = \sideset{}{^\wedge}\sum_{e\bar{e}=c} \psi(e) \sideset{}{^\wedge}\sum\sideset{}{^\wedge}\sum_{(w\bar{w}, z\bar{z})=1} \alpha_w \beta_z \psi(wz)[ewz]$$

where $\alpha_w = \alpha(w\bar{w})$ and $\beta_z = \beta(z\bar{z})$. We can assume that $ewz$ is primitive or else the symbol $[ewz]$ vanishes. Thus $e$ is primitive and $(wz, \overline{wz}) = 1$ in which case the condition $(w\bar{w}, z\bar{z}) = 1$ reduces to $(w, z) = 1$. Applying (20.3) we get

$$\mathcal{L}^*(M,N) = \sideset{}{^*}\sum_{e\bar{e}=c} \psi(e)[e] \sideset{}{^\wedge}\sum\sideset{}{^\wedge}\sum_{(w,z)=1} \varepsilon \alpha_w \beta_z \psi(wz)\left(\frac{wz}{e}\right)[wz]$$

where $\varepsilon = \pm 1$ depends only on the quadrants in which $e, wz$, and $ewz$ are located. For each $e$ we can separate the variables $w, z$ in the $\varepsilon$–factor by the use of (20.15) and (20.19). Moreover we write

$$\psi(wz)\left(\frac{wz}{e}\right) = \psi(w)\left(\frac{w}{e}\right)\psi(z)\left(\frac{z}{e}\right)$$

getting bilinear forms of type $\mathcal{K}^*(M,N)$ with their coefficients $\alpha_w, \beta_z$ twisted by characters. These forms satisfy the bound (21.13). Summing over $e$ we conclude that $\mathcal{L}^*(M,N)$ satisfies (23.8). Now we can remove the condition that $(c, mn) = 1$ as follows

$$\begin{aligned}
\mathcal{L}^*(M,N) &= \sum\sum_{\substack{ab|c^\infty \\ (a,b)=1}} \sum\sum_{\substack{(m,n)=1 \\ (mn,c)=1}} \alpha(am)\beta(bn)\lambda(abcmn) \\
&\ll \sum\sum_{\substack{ab|c^\infty \\ (a,b)=1}} \tau(abc)\left(a^{-1}M + b^{-1}N\right)^{\frac{1}{12}}\left((ab)^{-1}MN\right)^{\frac{11}{12}+\varepsilon} \\
&\ll \tau(c)(M+N)^{\frac{1}{12}}(MN)^{\frac{11}{12}+\varepsilon} \ .
\end{aligned}$$

Finally we remove the condition $(m,n) = 1$ as follows

$$\begin{aligned}\mathcal{L}(M,N) &= \sum_d \sum_{(m,n)=1} \sum \alpha(dm)\beta(dn)\lambda(cd^2 mn) \\ &\ll \sum_d \tau(cd^2) d^{-\frac{23}{12}} (M+N)^{\frac{1}{12}} (MN)^{\frac{11}{12}+\varepsilon} \\ &\ll \tau(c)(M+N)^{\frac{1}{12}} (MN)^{\frac{11}{12}+\varepsilon} .\end{aligned}$$

This completes the proof of Proposition 23.1. $\square$

By analogy to $\mathcal{K}(N)$ and $\mathcal{K}^*(N)$ we construct the special linear forms

(23.9) $$\mathcal{L}(N) = \sum_{n \leqslant N} \lambda(mn) ,$$

(23.10) $$\mathcal{L}^*(N) = \sum_{\substack{n \leqslant N \\ (n,m)=1}} \lambda(mn) .$$

PROPOSITION 23.2. *Given a Hecke character $\psi$ defined by (17.17) and a positive integer $m$ we have the bounds*

(23.11) $$\mathcal{L}(N) \ll d(|k|+1)\tau(m)^4 \sqrt{m} N^{\frac{3}{4}} \log mN ,$$
(23.12) $$\mathcal{L}^*(N) \ll d(|k|+1)\tau(m)^2 \sqrt{m} N^{\frac{3}{4}} \log mN ,$$

*where the implied constant is absolute.*

*Proof.* We have

$$\mathcal{L}^*(N) = \sideset{}{^\wedge}\sum_{w\overline{w}=m} \psi(w) \sum_{\substack{z\overline{z} \leqslant N \\ (z,w)=1}} \psi(z)[wz]$$

so (23.12) follows from (22.4). Next we have

$$\mathcal{L}(N) = \sum_{e|m} \sum_{c|e^\infty} \sum_{\substack{n \leqslant N/c \\ (n,m)=1}} \lambda(cmn)$$

so (23.11) follows from (23.12). $\square$

## 24. Combinatorial identities for sums of arithmetic functions

Given a nice arithmetic function $f : \mathbb{N} \to \mathbb{C}$ we are interested in estimating sums of $f(n)$ twisted by the quadratic eigenvalue $\lambda(n)$. Since $\lambda(n)$ changes sign at random, one should expect a lot of cancellation so that a bound

(24.1) $$\sum_{n \leqslant x} f(n)\lambda(n) \ll x^{1-\delta}$$

for some $\delta > 0$ is not out of the question for many interesting functions normalized as to satisfy $f(n) \ll n^\varepsilon$. In this section we develop some combinatorial identities by means of which one can reduce this problem to estimates for general bilinear forms $\mathcal{L}^*(M, N)$ and for special linear forms over primes.

Fix $r \geqslant 2$. For every squarefree number, say
$$\ell = p_1 p_2 \ldots, \quad \text{with } p_1 > p_2 > \ldots,$$
we define the divisor $\mathfrak{d} = \mathfrak{d}(\ell)$ by setting

(24.2) $$\mathfrak{d} = p_1 \ldots p_r p_{2r} p_{3r} \ldots .$$

Here and throughout this section the product of primes terminates when it runs out of primes of $\ell$. As usual the empty product is defined to be one. Note that

(24.3) $$\mathfrak{d} \leqslant \ell^{1/r} p_1^{r-1} .$$

We shall call $\mathfrak{d} = \mathfrak{d}(\ell)$ the "separation" divisor of $\ell$ for reasons soon to be clear. Then we define the divisors $m = m(\ell)$, $n = n(\ell)$ by setting

(24.4) $$m = (p_{r+1} \ldots p_{2r-1})(p_{3r+1} \ldots p_{4r-1}) \ldots ,$$

(24.5) $$n = (p_{2r+1} \ldots p_{3r-1})(p_{4r+1} \ldots p_{5r-1}) \ldots .$$

Note that

(24.6) $$n \leqslant m \leqslant \mathfrak{d} n .$$

One can characterize $m, n$ solely in terms of the separation divisor. Indeed, writing
$$\mathfrak{d} = \pi_1 \pi_2 \ldots, \quad \text{with } \pi_1 > \pi_2 > \ldots ,$$
that is $\pi_j = p_j$ if $j \leqslant r$ and $\pi_{k+r-1} = p_{kr}$ if $k \geqslant 1$, we see that $m$ has the following properties:

($\mathfrak{d}^+$) the first largest $r - 1$ prime divisors of $m$ are in $(\pi_{r+1}, \pi_r)$,
   the second largest $r - 1$ prime divisors of $m$ are in $(\pi_{r+3}, \pi_{r+2})$,
   the third ..., etc.

Similarly $n$ has the following properties:

($\mathfrak{d}^-$) the first largest $r - 1$ prime divisors of $n$ are in $(\pi_{r+2}, \pi_{r+1})$,
   the second largest $r - 1$ prime divisors of $n$ are in $(\pi_{r+4}, \pi_{r+3})$,
   the third ..., etc.

Hence we obtain the factorization

(24.7) $$\ell = \mathfrak{d} m n .$$

and every squarefree $\ell$ has unique factorization of this type. By (24.6) and (24.7) it follows that

(24.8) $$m, n \leqslant \sqrt{\ell} .$$

Suppose $\ell \leqslant x$ and that $\ell$ is free of prime divisors larger than

(24.9) $$z = x^{1/r^2} .$$

Then it follows by (24.3) that the separation divisor is quite small, namely $\mathfrak{d} \leqslant D$ with

(24.10) $$D = x^{2/r} .$$

By the above combinatorics we arrive at

LEMMA 24.1. *Let $f(\ell)$ be an arithmetic function supported on squarefree numbers $\ell \leqslant x$. Let $r \geqslant 2$. Then we have*

(24.11) $$\sum_{\ell | P(z)} f(\ell) = \sum_{\mathfrak{d} \leqslant D} \sum\sum_{m,n \leqslant \sqrt{x}} \gamma_{\mathfrak{d}}^+(m) \gamma_{\mathfrak{d}}^-(n) f(\mathfrak{d} m n)$$

*where $\gamma_{\mathfrak{d}}^+, \gamma_{\mathfrak{d}}^-$ are the characteristic functions of those integers having the properties $(\mathfrak{d}^+)$, $(\mathfrak{d}^-)$ respectively, and $z$, $D$ are defined by (24.9), (24.10).*

In this lemma $P(z)$ denotes the product of all primes $\leqslant z$. To remove the restriction $\ell \mid P(z)$ we appeal to another combinatorial formula

(24.12) $$\sum_{\ell} f(\ell) = \sum_{\ell | P(z)} f(\ell) + \sum_{p>z} \sum_{q} f(pq) \nu(pq, z)^{-1}$$

where $\nu(\ell, z)$ denotes the number of prime factors of $\ell$ which are $> z$. Note here that $q$ runs over integers. We can write $\nu(pq, z)^{-1} = \gamma(q)$ since it does not depend on $p$. If $q \leqslant z$ then $\gamma(q) = 1$. We single out this part of (24.12) and insert (24.11) getting

PROPOSITION 24.2. *Let $f(\ell)$ be a function supported on squarefree numbers $\ell \leqslant x$. Then we have*

(24.13) $$\sum_{\ell} f(\ell) = \sum_{\mathfrak{d} \leqslant D} \sum\sum_{m,n \leqslant \sqrt{x}} \gamma_{\mathfrak{d}}^+(m) \gamma_{\mathfrak{d}}^-(n) f(\mathfrak{d} m n)$$
$$+ \sum\sum_{p,q>z} \gamma(q) f(pq) + \sum\sum_{p>z \geqslant q} f(pq)$$

*for any $r \geqslant 2$ where $\gamma_{\mathfrak{d}}^+, \gamma_{\mathfrak{d}}^-$ are the characteristic functions of the integers having the properties $(\mathfrak{d}^+)$, $(\mathfrak{d}^-)$ respectively, $\gamma(q) = (1 + \nu(q,z))^{-1}$, and $z$, $D$ are defined by (24.9), (24.10).*

*Remarks.* On the right-hand side of (24.13) we have double sums over $m, n \leqslant \sqrt{x}$ and $p, q > z$ each of which is a bilinear form in $f$. The last sum in (24.13) will be treated for each $q$ individually (since $q$ is relatively small) as a special linear form over primes $\sum_{p} f(pq)$.

## 25. Estimation of $S_\chi^k(\beta')$

In Section 18 we have completed the proof of Theorem 1 by an appeal to Proposition 17.2 which is yet to be established. In this section we reduce the proof of Proposition 17.2 to Theorem 2, actually in a somewhat more general form, and to the latter we devote the next (and last!) section.

Recalling (17.14), (17.15) and (23.1) we write

$$(25.1) \qquad S_\chi^k(\beta') = \sum_{(m,\Pi)=1} g(m)\mu(m)\gamma(m)\lambda(m)$$

where $g(m)$ is a smooth function supported on $N' \leqslant m \leqslant (1+\theta)N'$ with $N < N' < 2N$ and satisfying $g^{(j)} \ll (\theta N)^{-j}$ for $j = 0, 1, 2$. In fact $g(m)$ is the function $p(n)$ in (4.13) but we have changed the notation to avoid ambiguity with the variable in sums over primes. We put

$$(25.2) \qquad \gamma(m) = \sum_{c|m,\ c \leqslant C} \mu(c)$$

with $1 \leqslant C \leqslant N^{1-\eta}$ for a small $\eta > 0$. Therefore, changing the order of summation,

$$S_\chi^k(\beta') = \sum_{\substack{c \leqslant C \\ (c,\Pi)=1}}^\flat \sum_{(\ell,c\Pi)=1} \mu(\ell)g(c\ell)\lambda(c\ell) \ .$$

First we can assume that

$$(25.3) \qquad 1 \leqslant C \leqslant N^\eta$$

because the remaining partial sum of $S_\chi^k(\beta')$ over $m = c\ell$ with $N^\eta < c \leqslant N^{1-\eta}$ is a bilinear form of type $\mathcal{L}^*(A,B)$ with $A, B \leqslant N^{1-\eta}$, $AB < 2N$ for which Proposition 23.1 gives the bound (after splitting into dyadic boxes, normalizing the coefficients, and separating the variables in $g(c\ell)$ via the Mellin transform)

$$(25.4) \qquad \mathcal{L}^*(A,B) \leqslant N^{1-\frac{\eta}{12}+\varepsilon},$$

and this is sufficient for (17.18).

Now we decompose the inner sum over $\ell$ according to the formula (24.13) with $x = N$ getting

$$S_\chi^k(\beta') = \sum_{\substack{c \leqslant C \\ (c,\Pi)=1}}^\flat \left\{ \sum_{\mathfrak{d} \leqslant N^{2/r}} T(c,\mathfrak{d}) + T(c) - S(c) \right\},$$

where $T(c, \mathfrak{d})$ and $T(c)$ are bilinear forms

$$T(c, \mathfrak{d}) = \sum_{\substack{m,n \leqslant \sqrt{N} \\ (mn, c\Pi) = 1}} \gamma_{\mathfrak{d}}^+(m) \gamma_{\mathfrak{d}}^-(n) \mu(mn) g(cmn) \lambda(cmn)$$

$$T(c) = \sum\sum_{p,q > N^{1/r^2}} \gamma(q) \mu(pq) g(cpq) \lambda(cpq)$$

and $S(c)$ involves a long sum over primes with a small parameter, namely it is given by

$$S(c) = \sum\sum_{q \leqslant N^{1/r^2} < p} \mu(q) g(cpq) \lambda(cpq).$$

The first double sum is a bilinear form of type (23.7) for which, after separating the variables $m, n$ in $g(cmn)$, Proposition 23.1 gives

$$T(c, \mathfrak{d}) \ll N^{\frac{23}{24} + \varepsilon} .$$

Similarly for the second double sum Proposition 23.1 also gives

$$T(c) \ll N^{1 - 1/12r^2 + \varepsilon} .$$

The inner sum over primes in $S(c)$ may be estimated using Theorem $2^\psi$, to be proved in the final section. This gives

$$S(c) \ll \sum_{q \leqslant N^{1/r^2}} cd(|k| + 1) q N^{\frac{76}{77}} \ll cd(|k| + 1) N^{\frac{76}{77} + \frac{2}{r^2}}$$

uniformly in $c, d, |k|$, the parameters of the Grossencharacter involved in $\lambda$. Adding these bounds, summing over $c$, and choosing $r$ sufficiently large we obtain Proposition 17.2.

It remains to prove Theorem $2^\psi$, a task we have postponed to the final section since it may generate an interest on its own.

## 26. Sums of quadratic eigenvalues at primes

Recall that

(26.1) $$\lambda(n) = \sum_{z\bar{z}=n}^{\wedge} \psi(z)[z]$$

where $\psi(z)$ is the Grossencharacter (17.17) and $[z]$ is the Jacobi-Kubota symbol. We shall establish

THEOREM $2^\psi$.   *For any $c \geqslant 1$ we have*

(26.2) $$\sum_{n \leqslant x} \Lambda(n) \lambda(cn) \ll c\mathfrak{f} x^{\frac{76}{77}}$$

*where $\mathfrak{f} = d(|k| + 1)$ and the implied constant is absolute.*

We begin by stating a result, see [DFI, Lemma 9] where also the proof may be found, which is useful for separating integral variables $m$, $n$ constrained by an inequality $mn \leqslant x$.

LEMMA 26.1.    *For $x \geqslant 1$ there exists a function $h(t)$ such that*

$$\int_{-\infty}^{\infty} |h(t)| dt < \log 6x$$

*and for every positive integer $k$*

$$\int_{-\infty}^{\infty} h(t) k^{it} dt = \begin{cases} 1 & \text{if } k \leqslant x \\ 0 & \text{otherwise.} \end{cases}$$

For the proof of Theorem $2^\psi$ we split the sum according to the formula of R.C. Vaughan

$$(26.3) \qquad \Lambda(n) = \sum_{\substack{a|n \\ a \leqslant y}} \mu(a) \log \frac{n}{a} - \sum_{\substack{ab|n \\ a,b \leqslant y}} \mu(a)\Lambda(b) + \sum_{\substack{ab|n \\ a,b > y}} \mu(a)\Lambda(b) \ .$$

This identity is valid for $n > y$. If $n \leqslant y$ then the right side vanishes. Hence the sum (26.2) splits into $S_0 + S_1 - S_2 + S_3$ where

$$\begin{aligned} S_0 &= \sum_{n \leqslant y} \Lambda(n)\lambda(cn) \ll \tau(c)y \ , \\ S_1 &= \sum_{a \leqslant y} \mu(a) \sum_{m \leqslant x/a} \lambda(acm) \log m \ , \\ S_2 &= \sum_{a,b \leqslant y} \mu(a)\Lambda(b) \sum_{m \leqslant x/ab} \lambda(abcm) \ , \\ S_3 &= \sum_{\substack{ab \leqslant x \\ a,b > y}} \mu(a)\Lambda(b) \sum_{m \leqslant x/ab} \lambda(abcm) \ . \end{aligned}$$

In $S_1$ we first insert $\log m = \int_1^m t^{-1} dt$ to avoid partial summation and, only then, apply Proposition 23.2 to infer that

$$S_1 \ll c^{\frac{1}{2}} \mathfrak{f} \sum_{a \leqslant y} \tau(ac)^4 a^{-\frac{1}{4}} x^{\frac{3}{4}} (\log cx)^2 \ll c\mathfrak{f} y^{\frac{3}{4}} x^{\frac{3}{4}+\varepsilon}.$$

For $S_2$ we directly apply Proposition 23.2, this time getting

$$S_2 \ll c\mathfrak{f} y^{\frac{3}{2}} x^{\frac{3}{4}+\varepsilon} \ .$$

By Proposition 23.1 we infer that (after splitting the variables into dyadic segments and separating them by the aid of Lemma 26.1)

$$S_3 \ll cy^{-\frac{1}{12}} x^{1+\varepsilon} \ .$$

We sum up the above bounds and choose $y = x^{\frac{3}{19}}$ getting $c\mathfrak{f}x^{\frac{75}{76}+\varepsilon}$ which is slightly better than (26.2).

*Remarks.* Theorem $2^{\psi}$ should be compared with the bound in [DI] for the corresponding sum over primes of the Fourier coefficients of a cusp form with respect to the theta multiplier. However the arguments of [DI] are rather different from these used here. Of course, we have not made an attempt to get the best bound by the available technology. For example if we had treated the part of $S_2$ corresponding to $ab \leqslant y$ as a sum of type $S_1$ and the remaining part as a sum of type $S_3$ the choice $y = x^{3/10}$ would give (26.2) with exponent $\frac{39}{40} + \varepsilon$. We challenge the reader (if she/he is not yet burnt out of energy as we are) to further reduce this exponent to a single digits fraction. It could be also interesting to improve the dependence of the bound (26.2) on the involved parameters $c$ and $\mathfrak{f}$ or simply to get the nontrivial bound

$$\sum_{n \leqslant x} \Lambda(n)\lambda(cn) \ll x^{1-\delta}$$

but uniformly in $c\mathfrak{f}$ as large as possible.

University of Toronto, Toronto, Canada
*E-mail address*: frdlndr@math.toronto.edu

Rutgers University
*E-mail address*: iwaniec@math.rutgers.edu

## References

[Bo1]  E. Bombieri, The asymptotic sieve, *Rend. Accad. Naz. dei XL*, **1/2** (1977), 243–269.
[Bo2]  ———, Le Grand Crible dans la Théorie Analytique des Nombres, *Astérisque*, **18** Soc. Math. France (Paris), 1987.
[BD]   E. Bombieri and H. Davenport, On the large sieve method, *Abhandlungen aus Zahlentheorie und Analysis Zur Erinnerung an Edmund Landau* (1877–1938), Deutscher Verlag der Wissenschaften, pp 11–22, (Berlin), 1968.
[BFI]  E. Bombieri, J. Friedlander, and H. Iwaniec, Primes in arithmetic progressions to large moduli, *Acta Math.* **156** (1986), 203–251.
[Bu]   D.A. Burgess, On character sums and $L$-series, *Proc. London Math. Soc.* **12** (1962), 193–206.
[DH]   H. Davenport and H. Halberstam, The values of a trigonometric polynomial at well spaced points, *Mathematika* **13** (1966), 91–96. ibid. **14** (1967), 229–232.
[Di]   P.G.L. Dirichlet, Démonstration d'une propriété analogue à la loi de reciprocité qui existe entre deux nombres premiers quelconques, *J. Reine Angew. Math.* **9** (1832), 379–389.
[DFI]  W. Duke, J. Friedlander, and H. Iwaniec, Bilinear forms with Kloosterman fractions, *Invent. Math.* **128** (1997), 23–43.
[DI]   W. Duke and H. Iwaniec, Bilinear forms in the Fourier coefficients of half-integral weight cusp forms and sums over primes, *Math. Ann.* **286** (1990), 783–802.
[FI]   E. Fouvry and H. Iwaniec, Gaussian primes, *Acta Arith.* **79** (1997), 249–287.

[FI1]   J. FRIEDLANDER AND H. IWANIEC, Bombieri's sieve, *Analytic Number Theory, Proc. Halberstam Conf.* ed. B.C. Berndt et al., Vol. 1 pp. 411–430, Birkhäuser (Boston), 1996.
[FI2]   ______, Using a parity-sensitive sieve to count prime values of a polynomial, Proc. Nat. Acad. Sci. USA **94** (1997), 1054–1058.
[FI3]   ______, Asymptotic sieve for primes, Ann. of Math., **148** (1998), 1041–1065.
[Go]    D. GOLDFELD, A simple proof of Siegel's theorem, Proc. Nat. Acad. Sci. USA **71** (1974), 1055.
[GM]    F. GOUVEA AND B. MAZUR, The square-free sieve and the rank of elliptic curves, J. Amer. Math. Soc. **4** (1991), 1–23.
[GR]    I.S. GRADSHTEYN AND I.M. RYZHIK, Table of Integrals, Series, and Products, Academic Press, Orlando, 1980.
[He]    E. HECKE, Über eine neue Art von Zetafunktionen, Math. Zeit. **6** (1920), 11–51.
[Ho]    C. HOOLEY, On the Barban-Davenport-Halberstam theorem I, J. Reine Angew. Math. **274/275** (1975), 206–223.
[Iw]    H. IWANIEC, Primes represented by quadratic polynomials in two variables, Acta Arith. **24** (1974), 435–459.
[Sh]    G. SHIMURA, On modular forms of half integral weight, Ann. of Math. **97** (1973), 440–481.
[Se]    A. SELBERG, On elementary methods in primenumber-theory and their limitations Cong. Math. Scand. Trondheim **11** (1949), 13–22.
[ST]    C.L. STEWART AND J. TOP, On ranks of twists of elliptic curves and power-free values of binary forms, J. Amer. Math. Soc. **8** (1995), 943–973.